\documentclass[oneside, a4paper,11pt]{amsart}

\usepackage[utf8]{inputenc}

\DeclareFontFamily{U}{wncy}{}
\DeclareFontShape{U}{wncy}{m}{n}{<->wncyr10}{}
\DeclareSymbolFont{mcy}{U}{wncy}{m}{n}
\DeclareMathSymbol{\sha}{\mathord}{mcy}{"58} 
\DeclareMathSymbol{\yer}{\mathord}{mcy}{"5E} 
\DeclareMathSymbol{\zem}{\mathord}{mcy}{"5A} 
\newcommand{\sym}{\sha}

\usepackage{geometry}

\usepackage{amsmath}
\usepackage{amsfonts}
\usepackage{amsthm}
\usepackage{amssymb}
\usepackage{graphicx}

\usepackage{tikz}
\usetikzlibrary{arrows}

\usepackage{aliascnt}
\usepackage{xcolor}
\usepackage[hyperfootnotes=false]{hyperref}
\hypersetup{ 
	linkbordercolor=blue,
    linktoc=none, 
}

\usepackage{array} 
\usepackage{enumerate}

\allowdisplaybreaks

\usepackage{mathtools}

\theoremstyle{plain}
\newtheorem{theorem}{Theorem}[section]

\newtheorem*{acknowledgementX}{Acknowledgment}

\newtheorem{conjecture}[theorem]{Conjecture}

\newtheorem{lemma}[theorem]{Lemma}

\newtheorem{maintheorem}[theorem]{Main Theorem}

\newtheorem{definition}[theorem]{Definition}

\newtheorem{example}[theorem]{Example}

\newtheorem{problem}[theorem]{Problem}
\newtheorem*{problemX}{Problem}
\newtheorem{fact}[theorem]{Fact}

\newtheorem{remark}[theorem]{Remark}

\newcommand{\id }{\mathrm{id}}

\newcommand{\B}{\mathcal{B}} 	
\newcommand{\g}{\mathfrak{g}} 	
\renewcommand{\sl}{\mathfrak{sl}}	
\newcommand{\Z}{\mathbb{Z}}  	
\newcommand{\N}{\mathbb{N}}  	
\newcommand{\R}{\mathbb{R}}  	
\newcommand{\C}{\mathbb{C}}  	
\renewcommand{\S}{\mathbb{S}}  	

\newcommand{\V}{\mathcal{V}}	
\renewcommand{\H}{\mathcal{H}}

\newcommand{\Fp}{\mathrm{F_-}}
\newcommand{\Fm}{\mathrm{F_+}}
\newcommand{\Fpm}{\mathrm{F_\pm}}
\newcommand{\rFp}{\mathrm{\tilde{F}_-}}

\newcommand{\rFpm}{\mathrm{\tilde{F}_\pm}}

\newcommand{\ch}{\mathrm{ch}}
\renewcommand{\exp}[1]{{\mathrm e}^{#1}}

\renewcommand{\i}{{\mathrm i}}

\newcommand{\Y}{{\mathrm Y}}
\newcommand{\resY}{\mathrm{ResY}}

\newcommand{\W}{{\mathcal W}}

\newcommand{\Beta}{\mathrm{B}}
\newcommand{\Sel}{\mathrm{Sel}}
\renewcommand{\d}{\mathrm{d}}

\renewcommand{\Re}{\mathfrak{Re}}

\renewcommand{\mod}{\mbox{mod}}

\DeclareMathOperator{\ord}{ord}

\newcommand{\sgn}{\mathrm{sgn}}

\catcode`,\active

\catcode`\,12

\newcommand{\md}{\text{-}}
\newcommand{\Res}{\mathrm{Res}}

\newcommand{\res}[1]{\mathrm{Res}\left(#1\right)}

\newcommand{\hamburger}[4] 
{
  \thispagestyle{empty}
  \vspace*{-2cm}
  \begin{flushright}
    ZMP-HH / #2 \\
    Hamburger Beitr{\"a}ge zur Mathematik Nr. #3 \\
  \end{flushright}
  \vspace{0.5cm}
  \begin{center}
    \Large \bf
    #1
  \end{center}
  \vspace{0.2cm}
  \begin{center}	
    Simon D. Lentner \\
    Algebra and Number Theory, 
    University Hamburg,\\
    Bundesstra{\ss}e 55, D-20146 Hamburg \\
    {simon.lentner@uni-hamburg.de} \\
  
\end{center}
}

\begin{document}
\hamburger{Quantum groups and Nichols algebras \\ acting on conformal field theories}{17-6}{647}{}
\begin{abstract}
We prove that certain screening operators in conformal field theory obey the algebra relations of a corresponding Nichols algebra with diagonal braiding. Our result proves in particular a long-standing expectation that the Borel parts of  small quantum groups appear as the algebra of screening operators. 

The proof is based on a novel, intimate relation between Hopf algebras, vertex algebras and a class of multivalued analytic special functions, which are generalizations of Selberg integrals. We prove that the zeroes of these special functions correspond to the algebra relations of the respective Nichols algebra, by proving an analytical quantum symmetrizer formula for the functions. Moreover, certain poles of the functions encode module extensions and a Weyl group action. At other poles, the quantum {symmetrizer} formula fails and the screening operators generate an extension of the Nichols algebra.

The intended application of our result is the conjectural logarithmic Kazhdan-Lusztig correspondence. More generally, our result seems to suggest that non-local screening operators in an arbitrary vertex algebra should be described by appropriate Nichols algebras, just as local screening operators can be described by Lie algebras.


\end{abstract}
\title[]{}
\maketitle


\thispagestyle{empty}
\enlargethispage{.1cm}
\setcounter{tocdepth}{2}
\tableofcontents

\section{Introduction}
\subsection{Background}
\enlargethispage{1.6cm}
A \emph{vertex operator algebra} is, very roughly speaking, a commutative algebra that depends analytically on a complex variable $z$. More precisely, a vertex operator algebra~$\V$ is an infinite-dimensional graded vector space with a linear map  
$$\Y:\;\V\otimes_\C \V\to \V[[z,z^{-1}]],$$
where $\V[[z,z^{-1}]]$ denotes Laurent series in a formal variable $z$  with coefficients in $\V$. The axioms of a vertex operator algebra include a version of commutativity or locality, which relates $\Y(a,z)\Y(b,w)$ and $\Y(b,w)\Y(a,z)$ for $z,\;w,\;z-w\neq 0$. As an implication, one also has a version of associativity, which relates these two expressions to $\Y(\Y(a,z-w)b,w)$. 
An additional axiom requires that conformal transformations of the variable $z$ in $\Y(a,z)$ are compatible with an action of the Virasoro algebra on $\V$, which is part of the data. Standard mathematical textbooks on vertex operator algebras include \cite{Kac98,FB04}. Vertex operator algebras are motivated by physics, where they describe the holomorphic (chiral) part of a 2-dimensional quantum field theory with conformal symmetry.\\ 


Vertex operator algebras have a rich representation theory, as anticipated in physics in \cite{MS89} and proven in \cite{Huang08,HLZ-VIII}:
If the category of representations of a vertex operator algebra $\V$ is semisimple and finite, then $\V$ is called \emph{rational}.
In this case, the category of representations possesses a natural tensor product, dualities and a nondegenerate braiding; it is a \emph{modular tensor category}. 
Any modular tensor category gives rise to 3-dimensional topological invariants and to actions of mapping class groups \cite{Tur}, which again matches the expectation from physics. But for the category of representations of a vertex operator algebra $\V$, this topological data is accompanied by analytical data, which is associated to $2$-dimensional surfaces with complex structure, and which transforms nicely under the mapping class groups \cite{Zhu96}. For example, attached to any vertex algebra module $\mathcal{M}$ is an analytic function called \emph{graded dimension}~$\ch_{\mathcal{M}}(\tau)$. The graded dimensions of all simple $\V$-modules piece together to a vector-valued modular form, under the action of the mapping class group of the torus $\mathrm{SL}_2(\Z)$. \\

The most basic example of a vertex operator algebra is the \emph{Heisenberg vertex algebra} $\H^r$ of rank $r$. The simple modules over $\H^r$ will be called $\V_\lambda$ and they are parametrized by $\lambda\in\C^r$. The tensor product discussed above is here $\V_\lambda\boxtimes \V_\mu=\V_{\lambda+\mu}$, the braiding is given by $e^{\pi\i(\lambda,\mu)}$, and the graded dimension of $\V_\lambda$ for any $\lambda$ is essentially the inverse~of the Dedekind eta function. Physically, this vertex operator algebra describes $r$ free bosons. 

Evidently, the representation category of $\H^r$ is infinite, but it can be used to build an important class of rational vertex operator algebras: Let $L$ be an even integral lattice of rank $r$, embedded into $\C^r$. Then the direct sum of $\H^r$-modules $\bigoplus_{\mu\in L} \V_\mu$ can be endowed with the structure of a vertex operator algebra $\V_L$, called \emph{lattice vertex algebra}. The simple modules over $\V_L$ will be called $\V_{[\lambda]}$ and they have the form $\bigoplus_{\mu \in [\lambda]} \V_\mu$ for any coset $[\lambda]\in L^*/L$ in the dual lattice~$L^*$. Overall, the representation category of $\V_L$ is equivalent to the category of vector spaces graded by the finite abelian group $L^*/L$, and with braiding and associator given via the quadratic form $e^{\pi\i(\lambda,\lambda)}$. The graded dimension attached to the simple module $\V_{[\lambda]}$ is essentially the theta function associated to the lattice coset $[\lambda]$, and together these functions transform as a vector valued modular form under $\mathrm{SL}_2(\Z)$, with respect to the Weil representation of $\mathrm{SL}_2(\Z)$ associated to $L^*/L$ and $e^{\pi\i(\lambda,\lambda)}$.\\

\enlargethispage{2.1cm}
A similar behaviour is widely expected if the category of representations of a vertex operator algebra $\V$ is finite, but not semisimple. In this case, we call $\V$ \emph{logarithmic}, for reasons that become clear later.
Understanding and constructing logarithmic vertex operator algebras is a fundamental question, but only few examples are known so far. Some have applications to spin chains or to critical percolation, see the survey \cite{GJRSV13}. 

\enlargethispage{.6cm}
By results in \cite{HLZ-VIII}, the representation category of a logarithmic vertex operator algebra is a (non-semisimple) braided tensor category. Conjecturally, it is a non-semisimple modular tensor category in the sense of \cite{KL01,Shim16} or in a slightly weaker sense.\footnote{Added: The weaker duality structure is now discussed in \cite{ALSW21}. This includes the example of a lattice vertex algebra $\V_L$ with a shifted Virasoro action that is relevant to this article, see \ref{fact_Vir}} As such, it still produces topological invariants and mapping class group actions, but not a full topological field theory. Moreover, the graded dimensions of the simple modules do not form a vector-valued modular form; one has to add additional so-called pseudo-characters associated to projective covers, as described for example in the surveys \cite{M14,CG17}. This is analogous to the fact that the space of class functions of a non-semisimple Hopf algebra has a larger dimension than the number of simple modules.  \\ 

On the algebra side, the main sources for non-semisimple modular tensor categories are the representation categories of \emph{small quantum groups} $u_q(\g)$, which we now discuss. Drinfeld and Jimbo have constructed the quantum group $U_q(\g)$ associated to a finite-dimensional semisimple complex Lie algebra $\g$ as a deformation of the universal enveloping algebra $U(\g)$ by a formal parameter $q$. The quantum group $U_q(\g)$ is a Hopf algebra, which means that its category of representations has a tensor product. 

Lusztig observed in \cite{Lusz90a} that if $q$ is specialized to a root of unity of order $\ell$, then the representation category changes drastically in comparison to the representation category of $\g$, and he discovered the existence of a finite-dimensional Hopf algebra quotient~$u_q(\g)$. The representation category of $u_q(\g)$ is non-semisimple and, roughly speaking, all simple representations of $\g$ with weight larger than $\ell$ become indecomposables. Lusztig conjectured that this category is related to representations of $\g$ in finite characteristic $\ell$ (if prime) and related to representations of the affine Lie algebra $\hat{\g}$ at shifted negative level.
The first conjecture was proven in \cite{AJS94}, the second conjecture was proven by Kazhdan and Lusztig in a series of papers culminating in \cite{KL94}. \\

\noindent
Our main motivation for the present article is the following long-term goal:

\begin{problemX}[Logarithmic Kazhdan-Lusztig Correspondence]
Construct a logarithmic $\quad$ vertex operator algebra $\W(\g,\ell)$, whose representation category is a non-semisimple $\quad$ modular tensor category, which is equivalent to the representation category of the small quantum group $u_q(\g)$ with $q$ a primitive $\ell$-th root of unity (or the slight variation below)
\end{problemX}

It is widely believed among experts that such a vertex operator algebra $\W(\g,\ell)$ can be constructed by using  a strategy called free-field realization. Previously, this strategy was applied very successfully in \cite{Wak86,FF88} to realize affine Lie algebras $\hat{\g}$ and their so-called Hamiltonian reductions $\mathrm{W}(\g,\ell)$, and this construction is nowadays a cornerstone of the geometric Langlands correspondence, see \cite{F95}.     

The idea of free-field realization, in our case, is to construct $\W(\g,\ell)$ for even $\ell=2p$  as a subalgebra of the lattice vertex algebra $\V_L$, where $L$ is  the root lattice of $\g$, rescaled by~$p$. It is expected that certain linear maps, the so-called \emph{short screening operators}, give an action of the Borel part $u_q(\g)^+$ of the small quantum group on the sum of simple modules of $\V_L$, and that the kernel of this action in $\V_L$ is the vertex operator algebra $\W(\g,\ell)$ with the conjectured representation category. Moreover, it is expected that  the so-called \emph{long screening operators} give an action of the Borel part of the dual Lie algebra~$\g^\vee$, which extends on $\W(\g,\ell)$ to an action of all of $\g^\vee$, and whose kernel contains the Hamiltonian reduction $\mathrm{W}(\g,\ell)$. These two actions should combine to an action of parts of Lusztig's infinite-dimensional quantum group of divided powers $U_q^\mathcal{L}(\g)$.

 In the smallest example $\sl_2$, the Hamiltonian reduction $\mathrm{W}(\mathfrak{sl}_2,\ell)$ gives a vertex operator algebra associated to the Virasoro algebra. In turn, $\W(\mathfrak{sl}_2,\ell)$ is an extension by three generators $W^+,W^0,W^-$, which form the adjoint representation of $\mathfrak{sl}_2$, and this extension turns out to be the widely studied logarithmic triplet algebra $\mathcal{W}(p)$, first constructed in  \cite{Kausch91}. In terms of representation theory, this extension merges infinitely many simple modules of the Virasoro algebra to finitely many simple modules of $\mathcal{W}(p)$. These modules can be identified with the simple modules of the small quantum group $u_q(\mathfrak{sl}_2)$.\\

The main result in this article, which we discuss in the next subsection, is the proof that indeed the short screening operators give an action of $u_q(\g)^+$. It also relates certain poles on the analysis side to non-split extensions of modules over $u_q(\g)^+$ and to Weyl reflections.
While this result is very far from establishing the logarithmic Kazhdan-Lusztig correspondence, it may be considered to be a necessary first step and it shows that indeed the small quantum group appears on the vertex operator algebra side.
 
More generally, our proof produces actions of arbitrary diagonal Nichols algebras (finite or infinite dimensional) via screening operators on corresponding lattice vertex algebras. In fact, we would expect much more generally that non-local screening operators on any vertex operator algebra produce an action of a corresponding Nichols algebra. At this point, one could speculate about scenarios, where the role of $u_q(\g)$ in the conjectures above is taken by other pointed Hopf algebras and beyond. Overall, we are interested in the following long-term goal, which today is surely out of reach:

\begin{problemX}
	For every non-semisimple modular tensor category $\mathcal{C}$, construct a logarithmic vertex operator algebra $\W$, whose representation category is equivalent to $\mathcal{C}$.
\end{problemX}	
   
   \enlargethispage{1.3cm}
\noindent
We finally discuss the current state of the logarithmic Kazhdan-Lusztig correspondence:

The case $\g=\mathfrak{sl}_2$ with $\mathcal{W}(p)=\mathcal{W}(\mathfrak{sl}_2,\ell)$ and  $\ell=2p$ is today quite well understood. It initially appears in \cite{FHST,FGSTsl2}. The best understood case is $p=2$, in which case $\mathcal{W}(2)$ is isomorphic to a known vertex operator algebra called \emph{symplectic fermions}. Using this, the equivalence of abelian categories $\mathcal{W}(2)\md\mod\cong u_q(\mathfrak{sl}_2)\md\mod$ was proven in \cite{FGSTsl2}, and later \cite{GR} have constructed a quasi-Hopf algebra $\tilde{u}_q(\mathfrak{sl}_2)$ similar to ${u}_q(\mathfrak{sl}_2)$, and an equivalence of modular tensor categories $\mathcal{W}\md\mod\cong \tilde{u}_q(\mathfrak{sl}_2)\md\mod$, up to an open conjecture about symplectic fermions. For general $p$, the abelian category of representations $\mathcal{W}(p)\md\mod$ was largely determined in \cite{AM08}. The action of $u_q(\mathfrak{sl}_2)^+$ on the free field realization via screening operators is given in \cite{TW}, and their work is based on a special case of the Selberg integrals appearing in our article, see Section~\ref{sec_Equal}. The equivalence of abelian categories  for $\sl_2$ has been proven in \cite{TN}. A factorizable quasi-Hopf algebra  $\tilde{u}_q(\mathfrak{sl}_2)$ was constructed in \cite{CGR20}, and we would expect an equivalence of modular tensor categories between representations of $\W(\sl_2,\ell)$ and $\tilde{u}_q(\sl_2)$.
\enlargethispage{1cm}
\footnote{Added: In the meanwhile this has been proven in \cite{CLR21} for $p=2$ and in \cite{GN21} for arbitrary $p$.}

For $\g$ simply-laced of higher rank, the reader is referred to \cite{AM14} and to an intriguing unpublished paper by Feigin and Tipunin \cite{FT} that sketches a construction of the modules of $\W(\g,\ell)$ as cohomologies of a vertex algebra bundle over the flag variety. This is a direct analog to the Borel-Weil-Bott theorem for Lie algebras, and it predicts an explicit formula for the graded dimension of these modules, using the Euler characteristic.

On the other hand, in \cite{GLO18} we have constructed a quasi-Hopf algebra variant~$\tilde{u}_q(\g)$ of the quantum groups at roots of unity of even order, such that the representation category is~a non-semisimple modular tensor category. The construction specializes for $\sl_2$ to the quasi-Hopf algebra $\tilde{u}_q(\mathfrak{sl}_2)$ in \cite{GR,CGR20} mentioned above. Note~that~for quantum groups at even order roots of unity the representation category is rarely a modular tensor category, but as we show this can be remedied at the cost of introducing~a nontrivial associator.\footnote{Added: For a different view on this phenomenon from the perspective of algebraic groups see \cite{Ne21}} We expect $\tilde{u}_q(\g)$ to be the correct algebra side of the logarithmic Kazhdan Lusztig correspondence. The associator is readily visible on the vertex algebra~side.

Moreover, there was a remarkable attempt by Semikhatov and Tipunin to generalize the program from quantum groups to diagonal Nichols algebras \cite{ST}, such as  calculations for the super-Lie algebra $\mathfrak{sl}(2|1)$ in \cite{STsuper} and conjectural  central charges for each rank 2 Nichols algebra in \cite{SCharges}. This was the author's initial contact to the question.


\subsection{Main Results}

In this article we prove that non-local screening operators built from lattice intertwining algebras, as defined below, fulfill the relations of the Nichols algebra associated to the diagonal braiding that corresponds to the given non-locality. In particular, this proves the long-standing expectation that short screening operators for rescaled root lattices give an action of the Borel part of the small quantum group~$u_q(\g)^+$. It suggests that Nichols algebras are the natural algebra structures generated by non-local screening operators and their braidings, just as Lie algebras are the natural algebra structures generated by local screening operators. We strongly expect similar statements to be true beyond lattice intertwiners and beyond Nichols algebras of diagonal braidings.

Mathematically, our proof shows an intimate relation between Nichols algebras and multivalued analytic functions, which we call \emph{generalized Selberg integrals}: We prove that such a function has a zero, wherever we have a relation in the Nichols algebra of the diagonal braiding that is related to the monodromy of the multivalued function. This relation between Nichols algebras and complex analysis might be of independent interest. \\

\enlargethispage{2.5cm}
We now discuss these results more precisely: A crucial notion in the theory of vertex algebra modules is that of an  \emph{intertwining operator} $\Y_{\mathcal{M},\mathcal{N},\mathcal{L}}$ between any three given modules~ $\mathcal{M},\mathcal{N},\mathcal{L}$. Intertwining operators could be viewed as the vertex algebra analogs of $V$-balanced bilinear maps between modules $M,N,L$ over a commutative ring~$V$. In full analogy to this case, one defines the tensor product of vertex~algebra modules $\mathcal{M}\boxtimes\mathcal{N}$ to be the universal object that admits an intertwining operator of type ${\mathcal{M},\mathcal{N},\mathcal{M}\boxtimes\mathcal{N}}$. 

For example, for the Heisenberg vertex algebra $\H^r$ the tensor product of two simple modules $\V_\lambda,\V_\mu$ for $\lambda,\mu\in \C^r$ turns out to be ${\V_\lambda\boxtimes \V_\mu=\V_{\lambda+\mu}}$. This is due to the existence and essentially uniqueness of an intertwining operator, given on the generators $e^{\phi_\lambda},e^{\phi_\mu}$~by
\begin{align*}
\Y_{\V_\lambda,\V_\mu,\V_{\lambda+\mu}}:\;\V_\lambda\otimes_\C \V_\mu &\longrightarrow \V_{\lambda+\mu}\{z\}\\
e^{\phi_\lambda}\otimes e^{\phi_\mu} &\longmapsto e^{\phi_{\lambda+\mu}}\;z^{(\lambda,\mu)}+\text{higher terms}
\end{align*}  
where $\C\{z\}$ denotes power series in $z^m$ with complex exponents $m\in\C$, and where $(\lambda,\mu)$ denotes the standard scalar product in $\C^r$.

In \cite{Huang08} and in \cite{HLZ-VIII}, where intertwining operators may also involve $\log(z)$, the braiding $\mathcal{M}\boxtimes\mathcal{N}\to \mathcal{N}\boxtimes \mathcal{M}$  is defined by analytically continuing the intertwining operators counterclockwise from $z$ to $-z$. So the double braiding expresses, to which extend the intertwining operators are multivalued functions. For the Heisenberg algebra they involve  $z^{(\lambda,\mu)}$, accordingly the braiding between $\V_\lambda,\V_\mu$ is the scalar $e^{\pi\i (\lambda,\mu)}$.\\

Let us now fix some element in some vertex algebra module. In the present article, we fix some vector $\alpha\in\C^r$ and fix in the respective $\H^r$-module $\V_{\alpha}$ the generator $\exp{\phi_\alpha}$. We then consider left-multiplication of $\exp{\phi_\alpha}$ on any module $\V_{\lambda}$ in the sense that we use the intertwining operator that we have by definition of the tensor product $\V_{\alpha}\boxtimes \V_{\lambda}= \V_{\alpha+\lambda}$:
\begin{align*}
\Y_{\V_\alpha,\V_\lambda,\V_{\alpha+\lambda}}(e^{\phi_\alpha},z):\;
\V_{\lambda} &\longrightarrow \V_{\alpha+\lambda}\{z\}
\end{align*}
We then define the (non-local) \emph{screening operator} associated to  $\exp{\phi_\alpha}$ as the linear map
$$\zem_\alpha:=\Res(\Y(e^{\phi_\alpha},z)):\;\V_{\lambda}\longrightarrow \overline{\V}_{\alpha+\lambda}$$
where by the residue $\Res(f(z))$ of a multivalued function $f(z)$ on $\C^\times$ we mean the contour integral around a lift of the unit circle, which is not a closed path anymore. As a consequence, not just the coefficient of $z^{-1}$ but of all non-integral $z$-powers contribute to the screening operator, and the screening operator actually maps to the algebraic closure $\overline{\V}_{\alpha+\lambda}$. If we evaluate a product of two such screening operators, the resulting coefficients are infinite series that have to be checked for convergence. In our case we get generalized hypergeometric series, instead of the binomial coefficients appearing for local screenings. \\

\enlargethispage{1cm}
We ask: What are the algebra relations between such screening operators $\zem_{\alpha_1},\ldots,\zem_{\alpha_r}$? Our main result is that they are the relations of the Nichols algebra associated to the diagonal braiding $e^{\pi\i(\alpha_i,\alpha_j)}$. To prove this, we must prove that the hypergeometric series above, which could be regarded as the structure constants of the algebra of screening operators, have a zero wherever the Nichols algebra has a relation. 

It seems astonishing that the zeroes of a family of analytic functions should match the relations of a corresponding Nichols algebra, which is in general very complicated. To prove this fact, we prove an analytic quantum symmetrizer formula by writing our hypergeometric function as a contour integral of a multivalued function in variables $z_1,\ldots, z_n$ and rewriting this as a quantum symmetrizer of another contour integral. The braiding $e^{\pi\i(\alpha_i,\alpha_j)}$ comes from monodromies around branch points $z_i=z_j$. If this second integral is non-singular, then this proves for the first integral that any formal linear combination in the kernel of the quantum symmetrizer becomes zero. But the kernel of the quantum symmetrizer consists by definition of the relations in the respective Nichols algebra.\\
%

\newpage
\subsection{Annotated Table of Content}~\\

\noindent
\hyperref[sec_Nicholsalgebra]{\sc 2. Quantum Symmetrizers in Quantum Groups and Nichols algebras.}\;    First, we
 briefly review one of the equivalent definitions of the Nichols algebra associated to a vector space with braiding: It is the quotient of a free algebra by the kernel of the \emph{quantum symmetrizer} $\sha$. Yet, it is difficult to obtain defining relations. The Borel part of the small quantum group $u_q(\g)^+$ is such a Nichols algebra, see Example \ref{exm_quantumgroup}. Thus, our goal is to prove that a linear combination of products of screening operators vanishes, if the formal linear combination of products lies in the kernel of the quantum symmetrizer. \\

\noindent
\hyperref[sec_LatticeIntertwiners]{\sc 3. Calculating in Vertex Algebras using Hopf Algebras.}\;   \\

\noindent
\hyperref[subsec_HopfIntertwiners]{\sc 3.1. Some Hopf Algebra Preliminaries.}\;   
We introduce an approach developed by the author in \cite{Len07} that encapsulates much of the vertex algebra power series combinatorics into certain Hopf algebra structures. However, we wish to emphasize that using this framework is optional for the proof of our main results. 

The reader unfamiliar with vertex operator algebras may simply accept Definition \ref{def_HopfIntertwiner} as a source of interesting linear maps on a vector space $\V$ endowed with certain Hopf algebra structures, which is designed to ultimately fulfill the identity in Theorem \ref{thm_associativity}. The reader familiar with vertex algebras may skip to the explicit examples in the next section.\\

\noindent
\hyperref[subsec_LatticeIntertwiners]{\sc 3.2. The Lattice Intertwiner Algebra.}\;    
We introduce two examples, both well-studied, that are of main interest to our article: The \emph{Heisenberg vertex algebra}~$\H^r$ in Definition \ref{def_HeisenbergVOA} and the \emph{lattice intertwining algebra} $\V_\Lambda$ in Definition \ref{def_latticeVOA}, which is associated to a lattice $\Lambda$,  not necessarily integral. The lattice intertwining algebra is a non-local generalization of a lattice vertex algebra and it can be viewed as an algebra in the braided category of representations over $\H^r$ or over some underlying lattice vertex algebra $\V_L$.\\

\noindent 
\hyperref[subsec_QFT]{\sc 3.3. Quantum Field Theory in a Nutshell.}\;    This section is optional and contains a short introduction to the role of vertex operator algebras in quantum field theory. \\
 
 \noindent
 \hyperref[subsec_Residues]{\sc 3.4. Mode- and Residue-Operators.}\;    For each element $a\in \V_\Lambda$ we turn the intertwining operator $\Y(a,z)$ into a linear map by taking a formal residue, and call the resulting linear map $\resY(a)$ the (non-local) \emph{screening operator}. Here, the residue of a multivalued function is defined by lifting a circle to a path in the multivalued covering. As a consequence, to be mathematically precise, the screenings operator $\resY(a)$ is an infinite sum and it maps into the algebraic closure $\overline{\V}_\Lambda$. When we compose such operators, we have to prove convergence of the resulting infinite sums of coefficients (Lemma \ref{lm_ConditionalConvergence}).
 
 The reader familiar with vertex algebras should be warned that such non-local screening operators are in general badly behaved. For example, when acting on modules they do not admit a commutator formula
 and are not compatible with the action of the Virasoro algebra or the kernel of the screening (but suitable powers are, see {Section}~\ref{subsec_EqualWeyl}). 
 Nevertheless, non-local screening operators are crucial for certain constructions.\\
 
\noindent
\hyperref[subsec_Screenings]{\sc 3.5. Screening Operators.}\;    We concentrate on screening operators~$\resY(a)$ for two specific families of elements $a$ in the lattice intertwining algebra $\V_\Lambda$, namely~$\partial\phi_\alpha$ and $\exp{\phi_\alpha}$, which lead to screening operators $\yer_\alpha$ and $\zem_\alpha$ for any $\alpha\in \Lambda$. The operator $\yer_\alpha$ is simply dual to the  $\Lambda$-grading on $\V_\Lambda$. The operator $\zem_\alpha$ is the main subject of our article. \\

\noindent
\hyperref[subsec_TrivialLevel]{\sc 3.6. Example: Root Lattices of Lie Algebras.}\;    We give a trivial but instructive example for the main motivation of this article: If we take $\Lambda$ to be the (even integral) root lattice of a simply-laced semisimple Lie algebra $\g$ with simple roots $\alpha_i$, then the formula derived in the previous subsection suffices to prove that $\yer_{\alpha_i},\zem_{\alpha_i}$ give an action of the generators $H_{\alpha_i},E_{\alpha_i}$ of the Borel part of $U(\g)$ on $\V_\Lambda$. In fact, this $\V_\Lambda$ is isomorphic to the affine Lie algebra $\hat{\g}$ at level $1$, which is called the Frenkel-Kac-Segal construction \\


\noindent
\hyperref[sec_Associativity]{\sc 4. The Non-Local Associativity Formula.}\;    
 We apply our Hopf algebra framework to simplify products of intertwining operators, and subsequently simplify products of screening operators $\prod_{i=1}^n\resY(a_i)$ of arbitrary elements $a_1,\ldots, a_n$. In Theorem~\ref{thm_associativity} we derive a general formula that involves two types of terms: First, it involves certain Hopf algebra expressions that are essentially invariant under permutation and even rebracketing of the $\Y(a_i)$. Second, it involves certain complex numbers arising from generalized hypergeometric series, which we call the \emph{quantum monodromy numbers}
$$\Fpm((m_i,m_{ij})_{ij})=\hspace{-.2cm}\sum_{(k_{ij})_{ij}}
\prod_{i}\frac{(e^{2\pi\i(m_i+\sum_{i<j} m_{ij})}-1)/2\pi\i}{
	1+m_i+\sum_{i<j}(m_{ij}-k_{ij})+\sum_{j<i}k_{ji}}
\prod_{i<j} (\pm 1)^{k_{ij}}{m_{ij} \choose k_{ij}}$$
depending on parameters $m_{ij},m_i\in\C$ and summing over $k_{ij}\in\N_0$ for ${1\leq i<j \leq n}$. Below we will derive a more reasonable expression for them in terms of contour integrals.

\enlargethispage{1.6cm}
Our main example of interest is the lattice intertwining algebra $\V_\Lambda$ with $a_i=\exp{\phi_{\alpha_i}}$ for some fixed $\alpha_i\in\Lambda$. For example, if we apply a product of screening operators to an element of the form $v=\exp{\phi_\lambda}$ in $\V_\lambda$, then the result of Theorem \ref{thm_associativity} reads
\begin{align*}
\left(\prod_{i=1}^n\zem_{\alpha_i}\right)\exp{\phi_\lambda}
=\left(\prod_{i=1}^n\resY(a_i)\right)v
&=\sum_{k_1,\ldots,k_n=0}^\infty\left(\exp{\phi_\lambda}\prod_{i=1}^n\frac{\partial^{k_i}}{k_i!}\exp{\phi_{\alpha_i}}\right)\cdot \Fp((m_i+k_i,m_{ij})_{ij}) \\
&\vphantom{\prod_{i=1}^n}=\exp{\phi_{\alpha_1+\ldots+\alpha_n+\lambda}}\cdot \Fp((m_i,m_{ij})_{ij}) +\text{ higher terms}
\end{align*}
with the following set of parameters for $\Fp$
$$m_{ij}:=(\alpha_i,\alpha_j),\qquad m_i:=(\alpha_i,\lambda)$$

One should interpret Theorem \ref{thm_associativity} as the ``first half'' of a non-local operator product expansion formula: Locality, associativity and the commutator formula now only depend on properties of the analytic expressions $\Fpm$. For local vertex operators, correspondingly for $m_{ij}\in 2\Z$, these $\Fpm$ are simply binomial expressions, and they lead to the usual formula for the commutator of two screening operators. However, in the non-local setting we do not find such easy relations. Our goal is to prove that the relations between the~$\Fp$ match the relations in the Nichols algebra for diagonal braiding $q_{ij}:=e^{\pi\i m_{ij}}$.\\

\noindent
\hyperref[sec_QuantumSymmetrizer]{\sc 5. The Quantum Symmetrizer Formula.}~\\

\noindent
\hyperref[subsec_QuantumMonodromyNumbers]{\sc 5.1. Quantum Monodromy Numbers $\Fpm$ as Series.}\;    We define our central notion, the quantum monodromy numbers $\Fpm((m_i,m_{ij})_{ij})$, as series. For the smallest cases we can give closed expressions as follows, where $\Beta(a,b)$ is the Euler Beta function:  
\begin{align*}
\Fp(m_1)
&=\frac{(e^{2\pi\i m_1}-1)/2\pi\i}{m_1+1}\\
\Fp(m_1,m_2,m_{12})
&=\frac{e^{2\pi\i m_2}-1}{2\pi\i}\frac{e^{2\pi\i m_1+2\pi\i m_{12}}-1}{2\pi\i}
\;\frac{1}{m_1+m_2+m_{12}+2}\cdot\\
&\cdot \left(\Beta(m_2+1,m_{12}+1)+\frac{\sin\pi m_1}{\sin\pi(m_1+m_{12})}\Beta(m_1+1,m_{12}+1)\right)\\
\end{align*}

\noindent
\hyperref[subsec_smallness]{\sc 5.2. A Regularity Condition: Subpolar $(m_{ij})_{ij}$.}\;    
We call a set of complex para\-meters $(m_{ij})_{1\leq i<j\leq n}$ \emph{subpolar} if, roughly, for any index subset $J$ the sums of $m_{ij}$ over all $i,j\in J$ has a not too negative real part. Sometimes we restrict ourselves to intervals~$J$. 

As a relevant case, Lemma \ref{lm_smallness} establishes that $m_{ij}=(\alpha_i,\alpha_j)$ is subpolar if the $\alpha_i$ are elements in a positive definite euclidean vector space and all $(\alpha_i,\alpha_i)<1$. 

The main consequence of this notion is Lemma \ref{lm_Christian2}, which requires some combinatorial effort. It states that $\vert\prod_{1\leq i<j\leq n} (z_i-z_j)^{m_{ij}}\vert$ is bounded by a linear combination of terms $\prod_{1\leq i<n} |z_i-z_{i+1}|^{\hat{{m}}_{i,i+1}}$, where all $\hat{{m}}_{i,i+1}>-1$. So ``subpolar''means that the singularities of  $\prod_{1\leq i<j\leq n} (z_i-z_j)^{m_{ij}}$ at $z_i=z_j$ are integrable, and in this way being subpolar is a frequent requirement for convergence of $n$-fold integrals and sums in this article.   \\

\noindent
\hyperref[subsec_Selberg]{\sc 5.3. The Generalized Selberg Integrals $\Sel$.}\;    We study the integral 
 \begin{align*}
\Sel((m_i,\bar{m}_i,m_{ij})_{ij})
&:=\idotsint\displaylimits_{1> z_1>\ldots >z_{n}> 0} 
\prod_i z_i^{m_i}\prod_i (1-z_i)^{\bar{m}_i}\prod_{i<j} (z_i-z_j)^{m_{ij}}
\;\d z_1\cdots \d z_{n}
\end{align*}
For example $\Sel(m_1,\bar{m}_1)=\Beta(m_1+1,\bar{m}_1+1)$. If all $m_i,m_{i},m_{ij}$ are equal, then the famous Selberg integral \cite{Sel44} gives an expression in terms of Gamma functions. If all $m_i,\bar{m}_i$ are equal modulo $\Z$ and all $m_{ij}$ are equal modulo $2\Z$, then the Selberg integral can be expressed in terms of Jack polynomials; in our context this is used already in~\cite{TW}.  

Lemma \ref{lm_SelbergConvergence} shows that the generalized Selberg integral is absolutely convergent if $(m_{ij})_{ij}$ is subpolar and if the real parts of $(m_i)_i$ are not too negative.
It will be highly relevant in Section \ref{subsec_EqualWeyl} that $\Sel$ has many singularities depending on  $(m_i)_i$ and $(m_{ij})_{ij}$, while $\Fp$ and $\rFp$ discussed below have few singularities that depend only on $(m_{ij})_{ij}$.\\

\noindent
\hyperref[subsec_redF]{\sc 5.4. The Reduced Quantum Monodromy Numbers $\rFpm$.}\;    We study the integral
$$\rFpm((m_i,m_{ij})_{ij}):=\frac{1}{(2\pi\i)^n}\int_{\triangle}
	\prod_i z_i^{m_i}\prod_{i<j} (z_i\pm z_j)^{m_{ij}}\;\d z_1\cdots \d z_n$$ 
where we integrate over the following image of an $n$-simplex in the multivalued covering
$$\triangle:=\{(e^{2\pi\i t_1},\ldots,e^{2\pi\i t_n}) \mid 0<t_1<\cdots <t_n<1\}$$
This integral $\rFp$ converges if $(m_{ij})_{ij}$ is subpolar on intervals, and it can be expressed in terms of  generalized Selberg integrals, assuming they converge, which depends on $(m_i)_i$.  

It would be very interesting to know an analytical continuation of $\rFp$. We expect poles to be located   wherever being subpolar fails with one $m_{ij}$-sum being a negative \emph{integer}.\\

\noindent
\hyperref[subsec_Fintegral]{\sc 5.4. The Quantum Monodromy Numbers $\Fpm$ as Integrals.}\;    We study the integral 
	$$\Fpm((m_i,m_{ij})_{ij})=\frac{1}{(2\pi\i)^n}\int_{\Box}
	\prod_i z_i^{m_i}\prod_{i<j} (z_i\pm z_j)^{m_{ij}}\;\d z_1\cdots \d z_n$$ 
	where we integrate over $\Box$, a certain lift  of the $n$-torus to the multivalued covering, on which the integrand has branch points. Technically, we thus work with the limit $(r_i)_i\to 1$ of the integral over the following domain for $r_1>\ldots>r_n$ in the multivalued cover:  
	$$\Box_{r_1,\ldots r_n}:=\left\{(r_1e^{2\pi\i t_1},\ldots,r_ne^{2\pi\i t_n}) \mid 0<t_1,\cdots,t_n<1\right\}$$

This integral converges if $(m_{ij})_{ij}$ is subpolar. As intended, we find that the integral $\Fpm$ above and the series $\Fpm$ in Section \ref{subsec_QuantumMonodromyNumbers} are expressions for the same analytic function. 

It will be highly relevant in Section \ref{subsec_EqualLiouville} that there are parameters $(m_{ij})_{ij}$, which are not subpolar and $\rFp$ has a singularity, yet $\Fp$ still converges.\\

\enlargethispage{1cm}
\noindent
\hyperref[subsec_Main]{\sc 5.4. Main Theorem.}\;    The \emph{Analytic Quantum Symmetrizer Formula} in Theorem \ref{thm_QuantumSymmetrizer}, for $(m_{ij})_{ij}$ subpolar, is the main result of this section. It states that
 $$
\Fp((m_i,m_{ij})_{ij})=\sum_{\sigma\in \S_n} q(\sigma)\;\rFp((m_{\sigma^{-1}(i)},m_{\sigma^{-1}(i)\sigma^{-1}(j)})_{ij})=:\sym_{n,q}\;\rFp((m_i,m_{ij})_{ij})
$$
 for the braiding matrix $q_{ij}=e^{\pi\i\;m_{ij}}$ and the associated braiding factors $q(\sigma),\;\sigma\in\S_n$. The main idea of proof is to subdivide the integral $\Fp$ over $\Box$ into $n!$ integrals $\rFp$ over~$\triangle$:
 $$\Box=\bigcup_{\sigma\in\S_n} \sigma\triangle$$
 
  The key observation in this article is that, as a consequence of writing $\Fp$ as quantum symmetrizer of some $\rFp$, we have proven that any linear combination of $\Fp$ vanishes, whenever it is in the kernel of the quantum symmetrizer. This holds independently of the explicit expression for $\rFp$, as long as $(m_{ij})_{ij}$ is subpolar and thus $\rFp$ is convergent.
 
 In Example \ref{exm_F2} we give for $n=2$ again explicitly the quantum symmetrizer formula and discuss in this case the role of the singularities of both $\rFp$ and $\Sel$.
 
As another consequence, we prove in Lemma \ref{lm_ConditionalConvergence} that also the series for $\Fp$ is conditionally convergent if $(m_{ij})_{ij}$ is subpolar, and even slightly beyond.\\

\noindent
\hyperref[sec_Equal]{\sc 6. Example: The Case of Equal $m_i,m_{ij}$.}\;   
 This section treats the important and instructional case, where all $m_i$ are all equal to some $m_{\alpha\lambda}$ and where all $m_{ij}$ are equal to some $m_{\alpha\alpha}$. In some assertions we only require equality modulo $\Z$ and $2\Z$, respectively. This appears in particular for a power of one screening $(\zem_\alpha)^n$ acting on a module $\V_\lambda$. \\

\noindent
\hyperref[subsec_EqualPrefactors]{\sc 6.1. Calculating the Combinatorial Prefactors.}\;  We calculate in the case of equal $m_i,m_{ij}$ modulo $\Z,2\Z$ certain combinatorial prefactors that appear when passing from $\Sel$ to $\rFp$ and to $\Fp$. These prefactors nicely factorize into $q$-numbers depending on $q_{\alpha\lambda}^2=e^{2\pi\i m_{i}},q_{\alpha\alpha}=e^{\pi\i m_{ij}}$. Their zeroes match the poles we encounter below. \\

\hyperref[subsec_EqualCalc]{\sc 6.2. Calculating $\Fp,\rFp,\Sel$.}\;  For equal $m_i,m_{ij}$ we derive explicit formulae for the functions $\Sel$, $\rFp$ and $\Fp$, using the Selberg integral formula and the prefactors calculated in the previous subsection. If we just have an equality modulo~$\Z,2\Z$, then we encounter  Jack polynomials, as they appear in \cite{TW}. 

We exhibit three families of poles interacting with the zeroes of the prefactors: The first poles appear in $\Sel$, depending on $(m_i)_i$, and they disappear in $\rFp$, where they cause exceptionally non-zero values. The second poles in $\Sel,\rFp$, appear for non-subpolar~$(m_{ij})_{ij}$. In mild cases $\Fp$ remains convergent, see Lemma \ref{lm_ConditionalConvergence}, but the divergency of $\rFp$ causes exceptionally non-zero values of $\Fp$ that violate the Nichols algebra relation. The remaining poles of $\Sel,\rFp,\Fp$ lead to infinities in the composition of screening operators.\\

\enlargethispage{2.5cm}
\noindent
\hyperref[subsec_EqualWeyl]{\sc 6.3. The Pole of $\Sel$ at $m_{\alpha\lambda}+(n-1)m_{\alpha\alpha}/2=-1$: Reflection operators.}\;    
 We discuss the first type of poles. They cause exceptionally non-zero terms in screening operator powers $(\zem_\alpha)^n$, depending on the acted-upon module $\V_{[\lambda]}$, and we shall call the resulting operators  \emph{reflection operators}. In rank $1$, this reproduces the powers of a screening appearing in the Felder complex \cite{Fel89}. In higher rank, we show on the level of weights that reflection operators correspond to the reflections defined in the theory of Nichols algebras. For example, for quantum groups we recover the Weyl group dot action. We conjecture that the reflection operators in fact give local operators, in the sense that they commute with the action of vertex algebra elements in the kernel of the screening operators. This could ultimately explain, why indecomposable modules on the vertex algebra side correspond to indecomposable modules on the Nichols algebra side. \\
 
 \noindent
 \hyperref[subsec_EqualLiouville]{\sc 6.4. The Poles of $\rFp$ at $m_{\alpha\alpha}=-2/n$: Violated Nichols Algebra Relations.} \;   
 We discuss the second type of poles. They serve as counterexamples to our main theorem, because now~$(m_{ij})_{ij}$ is non-subpolar and in turn a Nichols algebra relation fails, here~$(\zem_\alpha)^n\neq 0$ for $n=\ord(q_{\alpha\alpha})$. We can prove a remarkable occurrence, namely that the exceptionally non-zero operator $(\zem_\alpha)^n$ is equal to the  local operator $\zem_{n\alpha}$. The reasoning in the proof of this fact applies for any Nichols algebra relation if the parameters $(m_{ij})_{ij}$ are on the boundary of the subpolar region. Then, roughly speaking, many degrees of freedom differ only by terms inside the subpolar region and hence these terms vanish, so we get strong invariance properties. In the case at hand, the nonzero quantum monodromy numbers only depend on $\sum_i m_i$, not on the individual $m_i$. This invariance property then allows us to group all terms in the associativity formula in Section \ref{sec_Associativity} together and prove that the exceptionally non-zero power $(\zem_\alpha)^n$ is equal  $\zem_{n\alpha}$.

In particular, for a negative-definite rescaling of a root lattice, where $(m_{ij})$ is not subpolar, we expect that the algebra of screening operators $\zem_{\alpha_i}$ is the positive part of the quantum group of Kac-DeConcini-Procesi with an infinite-dimensional center, which is generated by the non-zero powers $\zem_{n\alpha_i}$ and similar powers for non-simple roots. From~a physics perspective, we expect that the kernel of these screenings describes Liouville- and Toda theories \cite{FF93} at rational coupling constants, where they become non-semisimple.\\

 \noindent
\hyperref[sec_KazhdanLusztig]{\sc 7. Applications to the Logarithmic Kazhdan-Lusztig Correspondence.}~\\

\noindent
\hyperref[subsec_NicholsAlgebraAction]{\sc 7.1. Nichols Algebra Action.}\;    We apply the quantum symmetrizer formula to prove in Theorem \ref{thm_NicholsAlgebra} our main result: For basis elements $\alpha_1,\ldots,\alpha_{r}$ in a positive definite lattice satisfying $\Vert\alpha_i\Vert^2\leq 1$, the screening operators $\zem_{\alpha_i}$ fulfill the relations of the Nichols algebra associated to the diagonal braiding matrix $q_{ij}=e^{\pi\i(\alpha_i,\alpha_j)}$. 

The assumption  $\Vert\alpha_i\Vert^2\leq 1$ assures that for every list $(\alpha_{i(x)})_{x\in\{1,\ldots,n\}}$ the associated  matrix $m_{xy}:=(\alpha_{i(x)},\alpha_{i(y)})$ is subpolar, so $\rFp$ is convergent. This is no formality: If the assumption fails, then the Nichols algebra relations are violated and catch additional terms due to poles of $\rFp$. For example, in the situation of Section \ref{subsec_TrivialLevel} we have $\Vert\alpha_i\Vert^2=2$, so the Nichols algebra is trivial, but the comutators of screenings give the Lie algebra $\g$. For Liouville models we have  $\Vert\alpha_i\Vert^2<0$ and the powers of screenings become nonzero.  \\

\enlargethispage{2.3cm}

\noindent
\hyperref[subsec_OutlookKL]{\sc 7.2. The Logarithmic Kazhdan-Lusztig Correspondence.}\;    We discuss this program and try to formulate the conjectures that are implicitly assumed from the experts, in sufficiently large  generality. The author also wants to sketch some intermediate steps, that he would find useful to approach this problem further, and some suggestions to construct further examples beyond lattice vertex algebras and quantum groups.\\

\noindent
\hyperref[subsec_ExampleKL]{\sc 7.3 Example $\mathfrak{sl}_2$ at $\ell=4$.}\;    
We discuss the smallest example in light of the program.\\

\section{Quantum Symmetrizers in Quantum Groups and Nichols Algebras}\label{sec_Nicholsalgebra}

A Nichols algebra is a universal algebra associated to a vector spaces with {braiding}. Its first and most prominent appearance is in \cite{Lusz93} as a characterization of the positive part of the quantum group $U_q(\g)^+$ and the small quantum group $u_q(\g)^+$. In fact, any Hopf algebra contains naturally in its graded algebra a Nichols algebra. In this way,~Nichols algebras are a key component in the classification of pointed Hopf algebras in \cite{AS}.

\begin{definition}[Braided vector space]
	A braided vector space $(M,c)$ is a vector space~$M$ together with a bijective linear map $c:M\otimes M \to M\otimes M$ called \emph{braiding}, such that on~$M\otimes M\otimes M$ the \emph{braid relation} holds:
	$$(\id_M\otimes c)(c\otimes \id_M)(\id_M\otimes c)=(c\otimes \id_M)(\id_M\otimes c)(c\otimes \id_M)$$ 
\end{definition}
\noindent
For a braided vector space $(M,c)$ we have an action $\rho_n$ of the \emph{braid group} $\mathbb{B}_n$ on $M^{\otimes n}$, where the overcrossing of the $i$-th and $i+1$-th strand, which we denote by $\tau_{i,i+1}$, acts by
$$\rho_n(\tau_{i,i+1}):=\mathrm{id}_M^{\otimes (i-1)}\otimes c \otimes \mathrm{id}_M^{\otimes (n-i-1)}$$
The braidings that appear in this article are of the following form:
\begin{example}[Diagonal braiding]
	Let $M$ be complex vector space with a choice of a basis~$x_1,\ldots, x_r$ and let $(q_{ij})_{1\leq i,j\leq r}$ be an arbitrary matrix of invertible scalars $q_{ij}\in \mathbb{C}^\times$. Then we define the \emph{diagonal braiding} on $M$ to be the following map $M\otimes M\to M\otimes M$:
	$$c:\;(x_i\otimes x_j)\mapsto q_{ij}(x_j\otimes x_i)$$
\end{example}
\noindent

We now define the Nichols algebra of a braided vector space $(M,c)$. We remark that this definition applies more generally to an object $M$ in a braided tensor category.

\begin{definition}[Quantum symmetrizer and Nichols algebra]\label{def_QuantumSymmetrizer}~
	
\noindent Let $(M, c)$ be a braided vector space, and $\rho_n$ the action of the braid group $\mathbb{B}_n$ on $M^{\otimes n}$.
 \begin{itemize}
 	\item The \emph{Matsumoto section} $s:\mathbb{S}_n \rightarrow  \mathbb{B}_n$ is a set-theoretic section of the canonical surjection $\mathbb{B}_n\twoheadrightarrow \mathbb{S}_n$. The map $s$ is uniquely defined by mapping the transposition~$(i,i+1)$ to the braid $\tau_{i,i+1}$, and by mapping any product of such transpositions of minimal length (a \emph{reduced expression}) to the respective product of braids. 
 	
 	The map $s$ is well-defined, because a theorem of Matsumoto states that all reduced expressions for an element in $\S_n$ are connected by braid relations. Note that $s$ is not a group-theoretic section, because $(i,i+1)^2=1$, whereas $\tau_{i,i+1}^2\neq 1$.
 	\item  The \emph{quantum symmetrizer} is the following linear map\footnote{This Cyrillic symbol is called \emph{sha}.} $M^{\otimes n}\to M^{\otimes n}$ 
 	\begin{align*}
 	\sym_{n}: =\sum _{\sigma\in \mathbb{S}_n} \rho_n(s(\sigma))
 	\end{align*}
 	\item The \emph{Nichols algebra} associated to $(M,c)$ is the  tensor algebra (or free algebra) generated by $M$ modulo the kernel of the quantum symmetrizer in each degree 
 	\[\B(M) :=\bigoplus _{n\geq 0} M^{\otimes n}/ \ker (\sym_{n}).\]
 \end{itemize}
 \end{definition}

\begin{lemma}\label{lm_braidingfactor}
	For a diagonal braiding $q_{ij}$ we have explicitly
	$$\sym_{n,q}: =\sum _{\sigma\in \mathbb{S}_n} q(\sigma)\sigma,\qquad
	q(\sigma)=\prod_{i<j,\;\sigma(i)>\sigma(j)} q_{ij}$$ 
\end{lemma}	
\begin{proof}
	The formula is easily proven by induction: For $\sigma'=(k,k+1)\sigma$ with $\mathrm{length}(\sigma')=\mathrm{length}(\sigma)+1$ we have by definition $\rho_n(s(\sigma'))=\rho_n(\tau_{k,k+1})\rho_n(s(\sigma))$. On one hand this causes an additional factor $q_{\sigma^{-1}(k),\sigma^{-1}(k+1)}$ in $\sym_{n,q}$. On the other hand, the assumption of $\mathrm{length}(\sigma)+1$  means that the set of inversions $(i,j)$ increases by this new element:
\[\{(i,j)\mid i<j,\;\sigma'(i)>\sigma'(j)\}
=\{(\sigma^{-1}(k),\sigma^{-1}(k+1))\}\;\cup\;\{(i,j)\mid i<j,\;\sigma(i)>\sigma(j)\} \qedhere\] 
\end{proof}

This particular definition of Nichols algebras is due to Woronowicz \cite{Wor89} and Rosso \cite{Rosso97}, and in this way Nichols algebras will appear in the present article. 

Note however, that the kernels of the maps $\sym_{n,q}$ are in general hard to calculate in explicit terms. So this definition of $\B(M)$ does not mean that the relations are known.
In fact, a Nichols algebra comes with the structure of a Hopf algebra in a braided sense, and it can be defined as being the Hopf algebra fulfilling several equivalent universal properties. Most structural results discussed below are derived from this perspective.

\begin{example}[Rank $r=1$ \cite{Nichols78}]\label{exm_rank1}
    Let $M=x\mathbb{C}$ be a $1$-dimensional vector space with diagonal braiding given by $q_{11}=q\in\mathbb{C}^\times$. Then $M^{\otimes n}=x^n\mathbb{C}$ is $1$-dimensional and
    $$\sym_{n,q}=\sum_{\tau\in\S_n}q^{\mathrm{length}(\tau)}=\prod_{k=1}^n\frac{1-q^k}{1-q}=:[n]_q!$$
    From the factorization we see that this polynomial has zeroes at all roots of unity $q$ of order $k$ for $1<k\leq n$. Hence the Nichols algebra of the braided vector space $(x\mathbb{C},q)$ is
    \[\B(M)=\begin{cases}
		\mathbb{C}[x]/(x^\ell),\quad & \text{if $q\neq 1$ is a primitive $\ell$-th root of unity} \\
		\mathbb{C}[x], & \text{else}
             \end{cases}\]
\end{example}

\begin{example}[Quantum group]\label{exm_quantumgroup}
      Let $\g$ be a finite-dimensional complex semisimple Lie algebra of rank $r$ with root system $\Phi$, simple roots $\alpha_1,\ldots,\alpha_r$ and Killing form $(\alpha_i,\alpha_j)$, normalized so that  $(\alpha_i,\alpha_i)\in\{2,4,6\}$. Let $q$ be a primitive $\ell$-th root of unity. Consider the $r$-dimensional vector space $M$ with diagonal braiding
      $q_{ij}:=q^{(\alpha_i,\alpha_j)}$.
      Then the Nichols algebra $\B(M)$ is isomorphic to the positive part of the small quantum group $u_q(\g)^+$.
\end{example}

In fact, a general Nichols algebra possesses many structures familiar from Lie algebras: It has a \emph{generalized root system}, a \emph{Weyl groupoid} and a \emph{Poincare-Birkhoff-Witt (PBW) basis}, see \cite{Heck09,HS10,AHS10}. A generalized root system is a set of hyperplanes, one for each root, that fulfills some crystallographic integrality condition. However, in contrast to root systems of Lie algebras, not all choices of simple roots (Weyl chambers), have the same braiding matrix or even the same Cartan matrix. This behaviour is already familiar from Lie superalgebras, which can be seen as special cases of Nichols algebras:

\begin{example}[Quantum group for $\mathfrak{sl(2|1)})$]
     Let $q$ be a primitive $\ell$-th root of unity. Consider the braided vector spaces $M',M''$ with the following different braiding matrices:
     \[q_{ij}'=\begin{pmatrix}
               -1  & q \\
               q & -1
              \end{pmatrix}
    \qquad
    q_{ij}''=\begin{pmatrix}
               -1  & q^{-1} \\
               q^{-1} & q^{2}
              \end{pmatrix}\]
    The associated root system is the standard Lie algebra root system of type $A_2$ with three positive roots. However, two of the roots are decorated by self-braiding $q_{\alpha\alpha}=-1$ and one root has self-braiding $q_{\alpha\alpha}=q^2$, hence there are two different types of Weyl chambers. 
    
    More precisely, in the first type of chamber with braiding matrix $q_{ij}'$ the simple roots $\alpha_1',\alpha_2'$ are decorated with $q'_{11}=q'_{22}=-1$ and the third root $\alpha_1'+\alpha_2'$ is decorated with the bimultiplicative combination $q'_{11}q'_{12}q'_{11}q'_{21}q'_{22}=q^{2}$. After reflection on $\alpha_1'$ we have simple roots $-\alpha_1',\alpha_1'+\alpha_2'$, which are decorated with $-1,q^{2}$ and the braiding matrix is $q_{ij}''$. 
    
    The Nichols algebras $\B(M'),\B(M'')$ are not isomorphic, but have the same dimension $2\cdot 2\cdot \mathrm{ord}(q^2)$. They are different Borel parts of the small quantum supergroup $u_q(\mathfrak{sl(2|1)}))$.
\end{example}
\noindent
 In this case, both set of simple roots give rise to the same Cartan matrix and Dynkin diagram of type $A_2$. We also give an example where different Cartan matrices appear:
\begin{example}[Quantum group for $\mathrm{D(2|1;\alpha)}$]
	 For roots of unity $q,s,r$ with $qrs=1$, we consider the following different braiding matrices of rank $3$:
	\[\begin{pmatrix}
	q^2  & q^{-1} & 1 \\
	q^{-1}  & -1 & r^{-1} \\
	1 & r^{-1} & r^2 
	\end{pmatrix},\quad
	\begin{pmatrix}
	q^2  & q^{-1} & 1 \\
	q^{-1}  & -1 & s^{-1} \\
	1 & s^{-1} & s^2 
	\end{pmatrix},\quad
	\begin{pmatrix}
	r^2  & r^{-1} & 1 \\
	r^{-1}  & -1 & s^{-1} \\
	1 & s^{-1} & s^2 
	\end{pmatrix},\quad
	\begin{pmatrix}
	-1  & q & s \\
	q  & -1 & r \\
	s & r & -1 
	\end{pmatrix}\]
	The associated generalized root system has $7$ positive roots. The corresponding arrangement of $7$ hyperplanes $\alpha^\perp$ in $\R^3$ forms a subdivided Archimedean solid as follows:   	
	\medskip
	\begin{center}
		\includegraphics[scale=.38, angle=0]{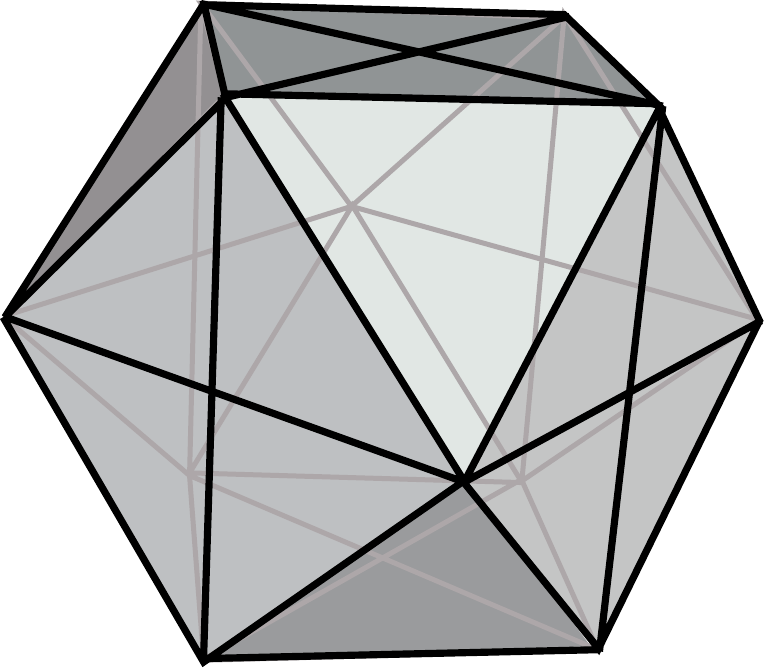}
	\end{center}
	\medskip
	The triangles are the different Weyl chambers, i.e. the choices of a basis of simple roots. Corresponding to the first three braiding matrices, there are $24$ right triangles, whose~Cartan matrix and Dynkin diagram is of type $A_3$. Corresponding to the fourth braiding matrix, there are $8$ equilateral triangles, whose Cartan matrix and Dynkin diagram resembles the affine type $A_2^{(1)}$. Reflecting an equilateral triangle on any side gives a right triangle. 
\end{example}

Finite dimensional Nichols algebras of diagonal type over fields of characteristic $0$ were classified by I. Heckenberger~\cite{Heck09}. Finite generalized root systems and Weyl groupoids are classified in \cite{p-CH10}: There are the familiar series $A_n,B_n,C_n,D_n$, the series $D_{m,n}$ familiar from Lie superalgebras, and 74 exceptional generalized root systems.  \\

\section{Vertex Algebras and Screening Operators}\label{sec_LatticeIntertwiners}

\subsection{Some Hopf Algebra Preliminaries}\label{subsec_HopfIntertwiners}

We introduce some Hopf algebra structures on infinite-dimensional graded vector spaces, from which we are able to construct a large class of vertex operators and intertwining operators. This is previous work of the author \cite{Len07} and will be useful during calculations, in particular with intertwining operators, since it encapsulates most power series combinatorics into finite algebra structures, which leaves only the relevant complex-analytic functions. This is largely done in Section~\ref{sec_Associativity}.

However, we wish to emphasize that using this general framework is optional for our main results. The reader familiar with vertex algebras may skip to the explicit lattice intertwining algebra $\V_\Lambda$ in the next subsection, and may in this case rewrite our algebraic expressions in terms of normally ordered products, creators and annihilators.\\

\enlargethispage{.5cm}

Let $H$ be a Hopf algebra with counit $\varepsilon:H\to \C$ and coproduct $\Delta:H\to H\otimes H$, which we denote in Sweedler notation $\Delta(h)=h^{(1)}\otimes h^{(2)}$. An $H$-module algebra $R$ is an algebra with an action of $H$, such that $h.(r'r'')=(h^{(1)}.r')(h^{(2)}.r'')$ and $h(1)=\varepsilon(h)$.
Similarly, an $H$-module coalgebra $V$ is a coalgebra with an action of $H$, such that $(h.v)^{(1)}\otimes (h.v)^{(2)}=h^{(1)}v^{(1)}\otimes h^{(2)}v^{(2)}$ and $\varepsilon(h.v)=\varepsilon(h)\varepsilon(v)$. 
If $H$ is a cocommutative Hopf algebra, then an $H$-module bialgebra or Hopf algebra $V$ is a bialgebra or Hopf algebra with an $H$-action that is compatible with both the algebra and coalgebra structure of~$V$. 
   
\begin{definition}
	We fix from now on the Hopf algebra $H$ to be the polynomial algebra $H:=\C[\partial]$, where the generator $\partial$ has coproduct $\partial^{(1)}\otimes \partial^{(2)}=1\otimes \partial+\partial\otimes 1$ and counit $\varepsilon(\partial)=0$. A module algebra over $H$ is simply an algebra together with a derivation $\partial$.
	
	As $H$-module algebra $R$ we fix from now on the algebra of finite linear combinations of $z^m$ with complex exponents $m\in \C$, where $\partial$ acts by derivation $-\frac{\partial}{\partial z}$. We also denote by $\C\{x\}$  the space of infinite linear combinations of $z^m$, as in \cite{HLZ-I}. 
\end{definition}
\noindent
The algebra $R$ contains the algebra of Laurent polynomials $\C[z,z^{-1}]$ and the space $\C\{z\}$ contains the space of Laurent series. Note that the series appearing later-on will be, by construction, bounded from below in the sense that the exponents come from a finite set translated by $\N$.
In more general non-semisimple situations, we would follow \cite{HLZ-II} and also include polynomials in $\log(z)$ into the definition of $R$. 
 Recall the subtlety that expansions with complex exponents are not unique, see \cite{HLZ-V} Def. 7.5.
\begin{definition}[\cite{Len07}]
	Let $\V$ be an $H$-module Hopf algebra. Then a \emph{Hopf pairing with coefficients} in the $H$-module algebra $R$ is a map 
	$\V\otimes \V\to R$ fulfilling
	\begin{align*}
	\langle a,bc \rangle
	&=\langle a^{(1)},b \rangle\langle a^{(2)},c \rangle
	\hspace{2cm} \langle a,\partial.b \rangle
	=-\frac{\partial}{\partial z}\langle a,b \rangle\\
	\langle ab,c \rangle
	&=\langle a,c^{(2)} \rangle\langle b,c^{(1)} \rangle
	\hspace{2.05cm}
	\langle \partial.a,b \rangle
	=\frac{\partial}{\partial z}\langle a,b \rangle
	\end{align*}
\end{definition}

The main purpose of introducing this algebraic data is that it can be used to write down vertex operators and intertwining operators.\footnote{The intended advantages of this approach compared to the standard one are as follows: First, the pairing $\langle a,b\rangle$ only contains finitely many terms, while $\Y$ contains infinitely many more terms, but these are all a consequence of the translation axiom. Second, associativity and locality with delta-functions appear automatically in our~$\Y$, see Section \ref{sec_Associativity}. Third, the first two axioms for $\langle a,b\rangle$ encapsulate a combinatorial version of Wick's theorem. Its standard form is recovered for $a,b$ products of primitive elements.}   

\begin{definition}\label{def_HopfIntertwiner}
	Let $\V$ be an $H$-module Hopf algebra and $\langle-,-\rangle(z)$ a Hopf pairing with coefficients in $R$. Then we define the associated \emph{intertwining operator} as follows:
	\begin{align*}
	\Y:\;\V\otimes_\C\V &\longrightarrow \V\{z\}\\
	a\otimes b&\longmapsto  \Y(a,z)b \; :=\; \sum_{k\geq 0}\langle a^{(1)},b^{(1)}\rangle(z) 
	\cdot b^{(2)}\cdot\frac{z^k}{k!}\partial^k a^{(2)}
	\end{align*}
	
\end{definition}

In \cite{Len07} the author has proven that this definition satisfies for general Hopf algebra data a generalized form of vertex algebra associativity and twisted locality. We discuss this in Section \ref{sec_Associativity}. For integral $z$-powers these properties reduce by \cite{Len07} Sec. 5.1 to the familiar vertex algebra axioms, under some additional assumptions. \\

\begin{example}
To illustrate the definition in general, we compute the following vertex operators involving the vacuum $1\in\V$: 
\begin{align*}
\Y(1,z)\;a
&={\langle 1, a^{(1)}\rangle}
\cdot a^{(2)} \cdot 
\sum_{k\geq 0}\frac{z^k}{k!}\;
\partial^k.1
= z^0\cdot a\\
\Y(a,z)\;1
&={\langle a^{(1)}, 1 \rangle}
\cdot 1 \cdot 
\sum_{k\geq 0}\frac{z^k}{k!}\;
\partial^k.a^{(2)}
=\sum_{k\geq 0}\frac{z^k}{k!} \cdot 
\partial^k.a\\
\end{align*}
\end{example}

\subsection{The Lattice Intertwining Algebra}\label{subsec_LatticeIntertwiners}

For the theory of vertex operator algebras the reader is referred to the textbooks \cite{Kac98,FB04}. We will concentrate on the main object of study in this article, the lattice intertwining algebra $\V_\Lambda$ \cite{DL93}, which can be viewed as a non-local generalization of the lattice vertex algebra to non-integral lattices. 

\begin{definition}[Heisenberg Vertex Algebra $\H^r$]\label{def_HeisenbergVOA}
Let $\C^r$ be a vector space with inner product $(,)$. We define the associated Heisenberg vertex algebra $\H^r$ as follows: 
\begin{itemize}	
\item $\H^r$ is the commutative algebra generated by formal symbols $\partial^k \phi_\alpha$ for all $k>0$ and all $\alpha\in \C^r$, together with the linear relations 
$$ a\;\partial^k\phi_\alpha+b\;\partial^k\phi_\beta
=\partial^k\phi_{a\alpha+b\beta},\qquad \forall a,b\in\C$$
If $\alpha_1,\ldots, \alpha_r$ is a basis of $\C^r$, then $\H^r$ is the polynomial ring in the variables~$\partial^k \phi_{\alpha_i}$.
\item $\H^r$ has an $\N_0$-grading via $|\partial^k\phi_\alpha|=k$. The first homogeneous components are: 
\smallskip
\begin{center}
	\begin{tabular}{c|cc}
		$\N_0$-degree & $\quad$dimension & basis \\
		\hline
		0 & 1 & 1 \\
		1 & $r$ & $\partial\phi_{\alpha_i}$  \\
		2 & $r+\frac{r(r+1)}{2}$ & $\partial^2\phi_{\alpha_i},\;\partial\phi_{\alpha_i}\partial\phi_{\alpha_j}$\\
		3 & $\cdots$ & $\cdots$ 	
\end{tabular}
\enlargethispage{1cm}
\end{center}
\item $\H^r$ becomes a module algebra over $H=\C[\partial]$ by setting
$$\partial.\left(\partial^k\phi_\alpha\right)=\partial^{k+1}\phi_\alpha$$ 
and by using the Leibniz rule on products of generators in $\H^r$.
\item $\H^r$ can be endowed with the structure of a vertex algebra. Usually, this is done by defining $\Y(\partial\phi_\alpha,z)$ in terms of the Heisenberg Lie algebra acting on $\H^r$ and then extending this to products in the argument by normally ordered products. 

In the context of this article, we can alternatively obtain $\Y(a,z)$ from Definition \ref{def_HopfIntertwiner}, by turning $\H^r$ into a Hopf algebra with $\Delta(\partial^k\phi_\alpha):=1\otimes \partial^k\phi_\alpha+\partial^k\phi_\alpha\otimes 1$ and by defining on $\H^r$ a Hopf pairings with coefficients, given on generators by 
$$\langle \partial \phi_\alpha, \partial\phi_\beta \rangle
:=(\alpha,\beta)z^{-2}$$
\end{itemize}
\end{definition}
As an example we compute 
\begin{align*}
\Y(\partial\phi_\alpha,z)\;\partial\phi_\beta
&=\langle \left(\partial\phi_\alpha\right)^{(1)}, 
\left(\partial\phi_\beta\right)^{(1)}\rangle
\cdot \left(\partial\phi_\beta \right)^{(2)}\cdot 
\sum_{k\geq 0}\frac{z^k}{k!}
\partial^k.\left(\partial\phi_\alpha\right)^{(2)}\\
&=\underbrace{\langle\partial\phi_\alpha,\partial\phi_\beta\rangle}_{=
	(\alpha,\beta)z^{-2}}\cdot 1 
\cdot \underbrace{\sum_{k\geq 0} \frac{z^k}{k!}\partial^k.1}_{=z^0\cdot 1}
+ \underbrace{\langle \partial\phi_\alpha, 1\rangle}_{=0} 
\cdot \partial\phi_\beta 
\cdot \underbrace{\sum_{k\geq 0} \frac{z^k}{k!}\partial^k.1}_{=z^0\cdot 1}\\
&+ \underbrace{\langle 1,\partial\phi_\beta\rangle}_{=0} \cdot 1 
\cdot \sum_{k\geq 0}\frac{z^k}{k!}\partial^k.\partial\phi_\alpha
+ \underbrace{\langle 1,1\rangle}_{=z^0} \cdot \partial\phi_\beta
\cdot \sum_{k\geq 0}\frac{z^k}{k!}\partial^k.\partial\phi_\alpha\\
&=(\alpha,\beta)z^{-2}\cdot 1+ \sum_{k\geq 0} \frac{z^k}{k!}\cdot
\partial\phi_\beta
\cdot\partial^{k+1}\phi_\alpha
\end{align*}

\enlargethispage{-.2cm}
\begin{definition}[Lattice Intertwiner Algebra $\V_\Lambda$]\label{def_latticeVOA}
	Let $\Lambda\subset \C^r$ be a lattice of rank $r$, with an inner product~$(,):\Lambda\times\Lambda\to \C$, which is not necessarily integer-valued. Then we define the \emph{lattice intertwining algebra} as follows:
	\begin{itemize}  
		\item $\V_\Lambda$ as a commutative algebra is defined as
		$$\V_\Lambda=\H^{r}\otimes \C[\Lambda]$$	 
		and we denote the basis of the group ring $\C[\Lambda]$ by formal symbols $\exp{\phi_\alpha},\alpha\in\Lambda$, so 
		$$\exp{\phi_\alpha}\exp{\phi_\beta}=\exp{\phi_{\alpha+\beta}}
		\hspace*{.4cm}$$
		\item $\V_\Lambda$ becomes a module algebra over $H=\C[\partial]$ by setting for all $\alpha\in \Lambda$
		$$\partial.\exp{\phi_\alpha}=\partial\phi_\alpha\exp{\phi_\alpha}$$ 
		\item $\V_\Lambda$ has a compatible grading by the lattice $\Lambda$
		$$\V_\Lambda=\bigoplus_{\lambda\in\Lambda} \V_\lambda$$ 
		where $\V_\lambda$ consists of all elements $u\exp{\phi_\lambda}$ with $u\in \H^r$, i.e. a differential polynomial.\newline
		As we will see, the $\V_\lambda$ are simple modules over the vertex subalgebra $\V_0=\H^r$.
		\item $\V_\Lambda$ can be endowed with the structure of an intertwining algebra. Usually this is done by directly defining $\Y(\exp{\phi_\alpha},z)$. In the context of this article, we can alternatively obtain $\Y$ from Definition \ref{def_HopfIntertwiner}, by turning $\V_\Lambda=\H^r\otimes \C[\Lambda]$ into a Hopf algebra with grouplike elements $\Delta(\exp{\phi_\alpha}):=\exp{\phi_\alpha}\otimes\exp{\phi_\alpha}$ and by defining on $\V_\Lambda$ the Hopf pairings with coefficients on the algebra generators by 
	  \begin{align*}
		\langle \exp{\phi_\alpha},\exp{\phi_\beta} \rangle
		&=z^{(\alpha,\beta)}\\
		\langle \exp{\phi_\alpha},\partial \phi_\beta \rangle
		&=-(\alpha,\beta)z^{-1}\\
		\langle \partial \phi_\alpha, \exp{\phi_\beta} \rangle
		&=(\alpha,\beta)z^{-1}\\
		\langle \partial \phi_\alpha, \partial\phi_\beta \rangle
		&=(\alpha,\beta)z^{-2}
	  \end{align*}
	  This assignment on algebra generators defines a Hopf pairing with coefficients, since the definition $\partial.\exp{\phi_\alpha}=\partial\phi_\alpha\exp{\phi_\alpha}$ is compatible with the  pairing's product rule. 
	\end{itemize}   
For most of this article we do not fix an action of the Virasoro algebra. As discussed in Section \ref{fact_Vir}, there is a $1$-parameter family of Virasoro actions on $\V_\Lambda$, and there we make~a choice compatible with a choice of screenings resp. a Weyl chamber in the {Nichols~algebra}.
\end{definition}

\noindent
For later use we define the differential polynomials $P_{\alpha,k}\in \H^r$ by
	$$\frac{1}{k!}\partial^k.\exp{\phi_\alpha}=P_{\alpha,k}\;\exp{\phi_\alpha}\;$$
The first few of these polynomials $P_{\alpha,k}$ are
$$1,\;\;
	\partial\phi_\alpha,\;\;
	\frac{1}{2!}\left(\partial\phi_\alpha\partial\phi_\alpha
	+\partial^2\phi_\alpha\right),\;\;
	\frac{1}{3!}\left(\partial\phi_\alpha\partial\phi_\alpha\partial\phi_\alpha
	+3\;\partial\phi_\alpha\partial^2\phi_\alpha
	+\partial^3\phi_\alpha\right),\;\;\ldots$$
 The coproduct of $P_{k,\alpha}$ follows from the fact that $\V_\Lambda$ is a $\C[\partial]$-module Hopf algebra:
	$$\Delta(P_{k,\alpha})
	=\sum_{k_1+k_2=k} P_{k_1,\alpha} \otimes P_{k_2,\alpha}$$

\noindent
As examples, we calculate the following intertwining operators in $\V_\Lambda$ involving $\exp{\phi_\alpha}$:
\begin{align*}
\Y(\exp{\phi_\alpha},z)\partial\phi_\beta
&=\left\langle \left(\exp{\phi_\alpha}\right)^{(1)}, 
\left(\partial\phi_\beta\right)^{(1)}\right\rangle
\cdot \left(\partial\phi_\beta\right)^{(2)}\cdot 
\sum_{k\geq0} \frac{z^k}{k!}\partial^k.\left(\exp{\phi_\alpha}\right)^{(2)}\\
&=\langle \exp{\phi_\alpha},\partial\phi_\beta \rangle
\cdot 1 \cdot 
\sum_{k\geq0} \frac{z^k}{k!}\partial^k.\exp{\phi_\alpha}
+\langle \exp{\phi_\alpha}, 1\rangle \cdot \partial\phi_\beta \cdot 
\sum_{k\geq0} \frac{z^k}{k!}\partial^k.\exp{\phi_\alpha}   \\
&=\left(-(\alpha,\beta)z^{-1}\cdot \exp{\phi_\alpha}+ z^0 \cdot \partial\phi_\beta \exp{\phi_\alpha} \right) 
\cdot\sum_{k\geq 0} z^k\;P_{k,\alpha}\\
\Y(\exp{\phi_\alpha},z)\exp{\phi_\beta}
&=\left\langle \left(\exp{\phi_\alpha}\right)^{(1)}, 
\left(\exp{\phi_\beta}\right)^{(1)}\right\rangle
\cdot \left(\exp{\phi_\beta}\right)^{(2)}\cdot 
\sum_{k\geq 0} \frac{z^k}{k!}\partial^k.\left(\exp{\phi_\alpha}\right)^{(2)}\\
&=z^{(\alpha,\beta)}
\cdot \exp{\phi_{\alpha+\beta}}\sum_{k\geq 0} z^k\;P_{k,\alpha}
\end{align*}
The last expression contains $z^{(\alpha,\beta)}$, which is not a single-valued function for $(\alpha,\beta)\not\in\Z$.\\

We conclude by discussing the definitions above from a representation theoretic point of view. The following construction has appeared in physics in work of Moore-Seilberg \cite{MS89}, and has been made mathematically rigorous in work of Huang \cite{Huang08} in the semisimple case and Huang-Lepowsky-Zhang \cite{HLZ-I}-\cite{HLZ-VIII} in general. 

\begin{theorem}
	Under suitable finiteness conditions (e.g. the so-called $C_2$-cofiniteness), the representations of a vertex operator algebra $\V$ form a braided tensor category. The tensor product $\mathcal{M}\boxtimes\mathcal{N}=\mathcal{L}$ is defined to be the universal object admitting an intertwiner $$\Y:\mathcal{M}\otimes_\C\mathcal{N}\to\mathcal{L}[\log(z)]\{z\}$$	
	The braiding is defined by analytically continuing this intertwiner counterclockwise from~$z$ to $-z$.  In particular if the intertwiner as a function in $z$ is multivalued, then the double braiding is non-trivial. In this case we call the intertwiner non-local.   
\end{theorem}
The Heisenberg algebra does not fulfill all the finiteness properties, but for example there is a semisimple braided tensor category generated by its irreducible representations:
\begin{example}
	For the Heisenberg algebra $\H^r$, the irreducible representations are the spaces $\V_{\lambda}=\H^r e^{\phi_\lambda}$ for all $\lambda\in \C^r$, with the vertex operator $\Y(\partial\phi_\alpha,z)e^{\phi_\lambda}$ above. The tensor product is $\V_{\lambda}\boxtimes \V_{\mu}=\V_{\lambda+\mu}$ with the intertwiner $\Y(\exp{\phi_\lambda})\exp{\phi_\mu}
	=z^{(\lambda,\mu)} \cdot \exp{\phi_{\lambda+\mu}}\cdots$ above. According to the multivalued function $z^{(\lambda,\mu)}$ the braiding is multiplication with~$e^{\pi\i(\lambda,\mu)}$.
\end{example}

An \emph{intertwining algebra} in general can be defined as any algebra inside the tensor category of representations over some vertex algebra $\V$. For any lattice $\Lambda\subset \C^r$, not necessarily integral, Definition \ref{def_latticeVOA} gives an intertwining algebra $\V_\Lambda$ over the Heisenberg vertex algebra, namely the sum of simple $\H^r$-modules $\V_\Lambda=\bigoplus_{\lambda\in \Lambda} V_\lambda$, endowed with the intertwining operators above as a multiplication $\V_\Lambda \boxtimes \V_\Lambda\to \V_\Lambda$. 
More generally, for any $2$-cocycle~$\epsilon:\Lambda\times \Lambda\to \C^\times$ there is an intertwining algebra $\V_\Lambda^\epsilon$ with modified~multiplication. \\

Finally we briefly discuss \emph{lattice vertex algebras}: If $\Lambda$ is an even integral lattice, then we have braidings $e^{\pi\i(\lambda,\mu)}=\pm1$ and trivial double-braidings. Hence for a suitable  choice of $\epsilon$ the intertwining algebra $\V_\Lambda^\epsilon$ becomes a commutative algebra  in the category of representations of the Heisenberg vertex algebra. Then $\V_\Lambda^\epsilon$ itself is again a vertex algebra, a so-called {simple current extension} of $\H^r$. Its representations correspond to the {local modules} over the commutative algebra $\V_\Lambda^\epsilon$ inside the category of $\H^r$-representations:   

\enlargethispage{1.5cm}

\begin{example}\label{exm_latticeVOA}
	For an even integral lattice $\Lambda$, the lattice vertex algebra $\V_\Lambda^\epsilon=\bigoplus_{\lambda\in \Lambda}\V_\lambda$ has irreducible representations $\V_{[\mu]}=\bigoplus_{\lambda\in\mu+ \Lambda}\V_\lambda$ for any coset $[\mu]\in \Lambda^*/\Lambda$ in the dual lattice $\Lambda^*$. The tensor product is $V_{[\lambda]}\boxtimes \V_{[\mu]}=\V_{[\lambda+\mu]}$, using the same intertwiner as for $\H^r$ and the additional factor $\epsilon(\lambda,\mu)$. The braiding is  $e^{\pi\i(\lambda,\mu)}\epsilon(\lambda,\mu)\epsilon(\mu,\lambda)^{-1}$ for some choice of coset representatives, but these choices may cause a nontrivial associator $\omega$. 
	
	Altogether, the category of representations of the lattice vertex algebra $\V_\Lambda^\varepsilon$ is equivalent to the semisimple modular tensor category of vector spaces graded by the finite abelian group $\Lambda^*/\Lambda$, with braiding and associator given by the quadratic form $e^{\pi\i(\lambda,\lambda)}$. 
\end{example}

To study our lattice intertwining algebra $\V_\Lambda$, we could choose an even integral sublattice $L\subset \Lambda$ and work over the lattice vertex algebra $\V_L$, but we find it easier to work directly over $\H^r$. However, in the application to Kazhdan-Lusztig correspondence in Section \ref{sec_KazhdanLusztig} we  works over the a lattice vertex algebra $\V_L$ containing the long screenings.

\subsection{Quantum Field Theory in a Nutshell}\label{subsec_QFT}
	For the curious reader, we briefly sketch the physical motivation and interpretation behind vertex algebras, in particular in the examples of the Heisenberg vertex algebra and the lattice intertwining algebra above.\\ 
	
	To define a model for a physical field, we first define the space of possible field configurations. As an example we may consider the space of all scalar functions $\phi:\Sigma \to \C$ on $\Sigma$, the $d$-dimensional spacetime or a worldsheet. 
	 Then the dynamics of the model is encoded by the choice of a functional on the space of field configurations called \emph{Lagrangian}~$L[\phi]$. For example, the Lagrangian  $L[\phi]=\int_{\Sigma}\mathrm{d}x^d\;\frac{1}{2}\Vert\nabla\phi\Vert^2$ describes a \emph{free scalar field}, without any interaction or external~force.\\
	
	In classical physics, only those field configurations $\phi$ appear ``in nature'', for which $L[\phi]$ is minimal. Equivalently, $\phi$ is a solution of a corresponding Euler-Lagrange differential equation. In the example of a free scalar field, we get harmonic functions $\Delta\phi=0$ for the euclidean space $\R^d$ and we get free waves $-\frac{\partial^2}{\partial t^2}\phi+\Delta\phi=0$ for the Minkowski space $\R^{1,d}$.
	
	 The fact that classical mechanics can be defined equivalently by variation problems (starting with the Fermat principle) or by differential equations was a major achievement of the 18th century mathematical physics. Both approaches generalizes to field theories, as described above, and remain present in quantum mechanics and quantum field theory. \\
	
	Quantum field theory, in the path integral formalism, replaces the minimality~condition of $L[\phi]$ by turning the space of field configurations into a probability space, with probability (more precisely: amplitude) of $\phi$ proportional to $e^{-\frac{\i}{\hbar} L[\phi]}$. For $\hbar\to 0$ only those $\phi$ with minimal $L[\phi]$ make notable contributions, so we recover the classical limit. 
	
	The interesting outputs are expectation values $\langle F\rangle$ of functionals $F[\phi]$. These are the path integrals of $F[\phi]$ over all fields $\phi$ with weight $e^{-\frac{\i}{\hbar} L[\phi]}$, which are mathematically ill-defined, but tremendously useful. Of particular interest are the expectation values of products of $\partial^k\phi$ evaluated at different points $z_1,\ldots,z_n\in \Sigma$, the \emph{$n$-point correlators}. For example, for a free scalar field the path integral can be computed using convolution products of Gauß integrals, and the $2$-point correlators are, after a Fourier transformation
	$$\left\langle\;\vphantom{X^X} \partial\phi(\vec{x_1})\;\partial\phi(\vec{x}_2)\;\right\rangle=\frac{1}{\Vert\vec{x}_1-\vec{x}_2\Vert^2},\qquad \vec{x}_1,\vec{x}_2\in\R^d $$

	In dimension $d=2$ we have an unusually large symmetry algebra, the local conformal transformations, encoded in an action of the Virasoro algebra, and this severely restricts the possible \emph{conformal field theories} (CFT). In this case, we identify the sphere $\Sigma=\R^2\cup\{\infty\}$ with the Riemann sphere $\hat{\C}$ by writing $z_i=x_i+\i y_i$. One can now decompose the correlators into holomorphic functions in $z_1,\ldots, z_n$ and in $\bar{z}_1,\ldots,\bar{z}_n$, which are called \emph{conformal blocks} or \emph{chiral $n$-point correlators}, and which are typically multivalued around branch points ${z_i=z_j}$. The two chiral parts are studied individually as chiral conformal field theories, for example as vertex operator algebras. A \emph{vertex operator algebra} consists of a formal space of functionals $\V$, an action of the Virasoro algebra (including $L_{-1}=\partial$, which corresponds to translation) and a vertex operator $Y(a,z)$. The matrix elements of the vertex operator express the chiral $3$-point correlators on the sphere $\Sigma=\hat{\C}$ via
	\begin{align*}
	\hspace{1cm}a^*(Y(b,z)c)
	\;:=\;\left\langle \vphantom{X^X}\;a(\infty)\;b(z)\;c(0)\;\right\rangle
	\end{align*}  
	where $a^*\in \V^*$. Note that conformal symmetry can move any three points on the Riemann sphere to $(z_1,z_2,z_3)=(\infty,z,0)$, so 
	this determines all chiral $3$-point correlators. By composing several $\Y$, the vertex algebra specifies a set of chiral $n$-point correlators on $\hat{\C}$
		$$1^*\left(Y(a_1,z_1)\cdots Y(a_n,z_n)1\right)
	=\left\langle \vphantom{X^X}\;a_1(z_1)\;\cdots\; a_n(z_n)\;\right\rangle$$

	For example, the Heisenberg vertex algebra $\H^r$ corresponds to the free chiral field $\phi:\Sigma\to \C^r$ with components $\phi_{\alpha_i}$. The space of functionals $\V$ consists of polynomials in derivations of $\phi$ at some point $z$, the previously introduced $\partial^k\phi_{\alpha_i}(z)$. The chiral correlator  
	$$\left\langle \;\vphantom{X^X} \partial\phi_\alpha(z_1)\;\partial\phi_\beta(z_2)\;\right\rangle=\frac{(\alpha,\beta)}{(z_1-z_2)^2}$$
	rewritten as $\langle 1(\infty) \partial\phi(z)\partial\phi(0)\rangle=(\alpha,\beta)z^{-2}$, translates into the first term of the vertex operator $\Y(\partial\phi_\alpha)\partial\phi_\beta=(\alpha,\beta)z^{-2}\cdot 1$ plus higher terms. In a similar way, the lattice vertex algebra $\V_\Lambda$ corresponds to the free chiral field that is compactified on a lattice $\phi:\Sigma\to \C^r/\Lambda$; this allows for additional fields $\exp{\phi_\lambda}$ with prescribed monodromy $\lambda\in \Lambda$.\\
	
By sewing together spheres, one can not only obtain chiral $n$-point correlators on the sphere, but on any Riemann surface $\Sigma$ of genus $g$. 	From a mathematical perspective, this is a fascinating feature of CFTs and vertex algebras. On one hand, already the representation theory of the vertex algebra can access topological data, which has triggered the discovery, that in fact any abstract modular tensor category can be used to construct topological invariants, mapping class group actions and often a full 3-dimensional topological field theory \cite{Tur,KL01}. On the other hand, the representations of a vertex algebra have by construction additional analytical objects attached to them. For example, the chiral correlator for $\Sigma$ a torus ($g=1,n=0$) coincides with the graded dimension of the representation, and from the physical interpretation it is evident that this should depend on the complex structure of the torus. Accordingly, it can be proven \cite{Zhu96} that the graded dimensions of vertex algebra representations are vector-valued modular forms under those action of the torus mapping class group $\mathrm{SL}_2(\Z)$ that is already visible on the topological field theory. For the Heisenberg algebra, we obtain eta functions, for the lattice vertex algebra we obtain theta functions of the respective lattice. \\    	

\newpage
\subsection{Mode- and Residue-Operators}\label{subsec_Residues}

We fix any element $a\in \V$ in a vertex algebra or intertwining algebra $\V$, we expand the operator  
$$\Y(a,z)=\sum_{m\in\C} \Y(a)_m z^{-1-m}
\qquad \Y(a)_m:\V\to \V$$
and we call the linear maps $\Y(a)_m$ \emph{mode operators}. They can also be expressed in terms of the Hopf pairing as follows, if $\langle a,b\rangle_m$ denotes the $z^m$-coefficient with $m\in\Z$ or  $m\in\C$:

\begin{align*}
&\Y(a)_mb
=\sum_{k\geq 0}\langle a^{(1)},b^{(1)}\rangle_{-1-m-k}
\;b^{(2)}\;\frac{1}{k!}\partial^k.a^{(2)}
\end{align*}

 For vertex algebras, the mode operators often give rise to interesting actions of infinite-dimensional Lie algebras. We give some examples: The mode operators ${L_m:=\Y(T)_{m+1}}$ of the conformal element $T$ of any vertex operator algebra give by definition rise to an action of the Virasoro algebra. For the Heisenberg vertex algebra, the mode operators $a_m:=\Y(\partial\phi)_m$ give rise to the action of the Heisenberg Lie algebra. The Wess-Zumino-Witten models are defined with mode operators from a given affine Lie algebra~$\hat{\g}$.


 The $0$-mode operator $\Y(a)_{0}$, or residue, plays a special role, since the commutator formula shows that this is a vertex algebra derivation, for example we have $L_{-1}=\partial$. In Section~\ref{subsec_TrivialLevel} we exemplary study the lattice vertex algebras of a root lattice, where $0$-mode operators give rise to an action of a Lie algebra $\g$. A main point of our article is that for intertwining algebras the $0$-mode operators should be replaced by the following more complicated notion, and this will ultimately give rise to quantum group actions: 

\begin{definition}\label{def_res}
We define a formal \emph{residue} $\Res(f)$ of a function $f\in\C\{z\}$ by setting
$$\mathrm{Res}(z^m)
  \;:=\; \begin{cases}
		0 & m\in\Z \backslash\{-1\}\\
		1 & m= -1\\
		\frac{1}{2\pi\mathrm i\;(m+1)}
		\left(e^{2\pi\mathrm{i}\;(m+1)}-1\right),\quad
		& m\not\in\Z \\             
		\end{cases}$$ 

{Geometrically}, this is the integral of $f(z)$ along the unique lift  of the circle $S^1$ of radius~$1$ to a path in the multivalued covering, starting at the principal branch.
\end{definition}
		
\begin{definition}
  For any given element $a\in \V$ we define the operator
  $\resY(a)$
  
  $$\resY{(a)}b:=\res{\Y{(a,z)b}}=\sum_{m\in\C} \res{z^{-1-m}}\;\Y(a)_mb$$
  If  $\Y(a,z)b$ contains $z$-powers $z^m$ with exponents in a single coset $m+\Z$, then explicitly
  $$
  \resY{(a)}\;b=\begin{dcases*}
		\hspace{4.0cm}\Y(a)_{0}b 
		& for $m\in\Z$\\
		\frac{e^{2\pi\i m}-1}{2\pi\mathrm \i}
		\sum_{k\in\Z}
		\frac{1}{m+k+1}\Y(a)_{-1-m-k}b
		& for $m\not\in\Z$              
		\end{dcases*}
		$$ 
\end{definition}~\\
\begin{remark}\label{rem_illdefines}
  The operator $\resY(a)$ is an infinite linear combination of mode operators $Y(a)_m$, each of which shifts the grading, hence $\resY(a)$ is a well-defined map from $\V$ to the algebraic closure of the graded vector space $\bar{\V}$. For compositions $\resY(a_1)\!\cdots\resY(a_n)$ we encounter infinite series in the coefficients and thus we have to check convergence. For the screening operators from $\V_\Lambda$, in which we are interested later, these series are explicit generalized hypergeometric series, given in Definition \ref{def_quantummonodromynumber}, and we will directly prove their convergence under additional assumptions on the parameters in Lemma \ref{lm_ConditionalConvergence}.\\
  
  \noindent
  A conceptually cleaner course of action in the future would be either:
  \begin{itemize}
   \item Along the lines of \cite{CKLW18}, consider a vertex algebra with norm\footnote{We thank Ingo Runkel for pointing this out.}. Then with  methods similar to Lemma \ref{lm_ConditionalConvergence} one should prove that $\resY(a)$ for $a\in\V$ with sufficiently small norm is a well-defined operator, and bound its operator norm.
   \item Consider $\resY(a)$ as an operator-valued analytic function depending on the radius parameter. This idea will appears in the proof of our main result, but will not be studied more systematically. It would then be very interesting to perform analytic continuations as a function in the radius. For our lattice intertwining  algebra, this analytic continuation will be a generalized hypergeometric function.\\
  \end{itemize}
\end{remark}

\subsection{Screening Operators}\label{subsec_Screenings}

\noindent
We now introduce residue operators $\resY(a)$ 
associated to certain elements $a\in\V_\Lambda$ in the lattice intertwining algebra:

\begin{definition}
 The residue operator of $\phi_\alpha$ defines the  
 linear operator\footnote{This Cyrillic symbol is called \emph{yeri}}
 $\yer_{\alpha}:\V_\Lambda\to {\V}_\Lambda$ 
 $$\yer_\alpha v:=\resY{(\partial\phi_\alpha)v}
 =\Y(\partial\phi_\alpha)_{0}v
 =\langle\partial\phi_\alpha,v^{(1)}\rangle_{-1}\; v^{(2)}$$
\end{definition}
\begin{lemma}\label{lm_ScalarProperties}
 The following properties can be directly computed from the last expression:
 \begin{enumerate}[a)]
  \item $\yer_{\alpha}$ is a {derivation} of the intertwining algebra $\V_\Lambda$ and the algebra $\V_\Lambda$. Explicitly:
    $$\hspace{.2cm}\yer_{\alpha}\;u\exp{\phi_\beta}=(\alpha,\beta)\;u\exp{\phi_\beta}, \quad \;u\in\H^r$$
  \item The map $\alpha\mapsto \yer_{\alpha}$ gives an action of the additive group $\Lambda$ on $\V_\Lambda$, dual to the  
  $\Lambda$-grading. 
  $$\yer_\alpha\yer_\beta=\yer_{\alpha+\beta}\qquad \yer_0=\id$$
 \end{enumerate}
\end{lemma}

\begin{definition}[Screening operator]\label{def_screening}
  The residue operator of $\exp{\phi_\alpha}$ defines the linear operator\footnote{This Cyrillic symbol is called \emph{ze} with a voiced s} ${\zem_\alpha:\V_\Lambda\to \overline{\V}_\Lambda}$, which we will call \emph{screening operator} in the following: 
  $$\zem_\alpha v:=\resY{(\exp{\phi_\alpha})v}$$
  If $(\alpha,\beta)\in \Z$ then $\zem_\alpha$ simplifies to the respective mode operator $\Y(\exp{\phi_\alpha})_0$, explicitly:
  \begin{align*}
  \zem_\alpha \; u\exp{\phi_\beta}
  &=\sum_{k\geq 0} \langle \exp{\phi_\alpha}, 
  u^{(1)}\rangle_{-k-(\alpha,\beta)-1}\; u^{(2)}\;P_{\alpha,k}\;\exp{ \phi_{\alpha+\beta} }  
  \end{align*}
\end{definition}
  Some authors use the term \emph{screening operator} for any $\resY(a)$, others require screening operators to additionally preserve a given choice of a Virasoro action.
\begin{example} Applying the definition (or for $(\alpha,\beta)\in \Z$ the previous formula) gives 
  \begin{align*}
    \zem_0\;v
    &=0 \quad \text{for any $v\in\V_\Lambda$}\\
    \zem_\alpha\;1 
    &= 0\\
    \zem_\alpha\;\partial\phi_\beta 
    &=-(\alpha,\beta)\exp{\phi_\alpha}\\
    \zem_\alpha\;\partial\phi_\beta\partial\phi_\gamma
    &=\left(\vphantom{x^2}-(\alpha,\beta)\partial\phi_\gamma-\partial\phi_\beta 
(\alpha,\gamma)
    +(\alpha,\beta)(\alpha,\gamma)\cdot\partial\phi_\alpha 
\right)\exp{\phi_\alpha}\\
    \zem_\alpha\;\exp{\phi_\beta}
    &=\begin{cases}
    \hspace{4.2cm}0,
    & \mbox{if } (\alpha,\beta)\in\N_0\\
    \hspace{4.2cm} P_{\alpha,k}\; \exp{\phi_{\alpha+\beta}},\qquad
     &\mbox{if }(\alpha,\beta)\in\!-\N,\;\; k:=-(\alpha,\beta)-1\\
      \displaystyle\sum_{k\geq 0} 
     \frac{e^{2\pi\i(k+(\alpha,\beta)+1)}-1}{2\pi\i(k+(\alpha,\beta)+1)}
     \cdot P_{\alpha,k}\;\exp{\phi_{\alpha+\beta}},
     &\mbox{if }(\alpha,\beta)\not\in\Z
    \end{cases}
  \end{align*}
\end{example}


\begin{lemma}\label{lm_ScreeningProperties}
 The following properties hold for $\zem_\alpha$ applied to a general element $v=u\exp{\phi_\beta}$ for any   $\alpha,\beta\in\Lambda$. For $(\alpha,\beta)\in\Z$ these are familiar properties of local operators.
 \begin{enumerate}[a)]
  \item The screening $\zem_{\alpha}$ shifts the $\Lambda$-grading, mapping the graded components $\V_\beta\to\V_{\alpha+\beta}$. 
  \item We have the following commutator relation with the derivation operator $\partial$:
  $$\hspace{-.2cm}[\partial,\zem_\alpha] v=
  \frac{e^{2\pi\i (\alpha,\beta)}-1}{2\pi\i}\left(\sum_{m\in\Z}\langle \exp{\phi_\alpha},u^{(1)}\rangle_m
    \;u^{(2)}\right)\exp{\phi_\beta}\;{\sum_{k\geq 0}\frac{\partial^k}{k!}.\exp{\phi_\alpha}}
    $$
  In particular if $(\alpha,\beta)\in \Z$, then $\zem_\alpha$ and $\partial$~commute.
  \item For $(\alpha,\beta)\in\Z$ we have the following mutual (anti-)commutator relations:
  \begin{align*}
  [\zem_\alpha,\zem_\beta]_\pm
  &=\begin{cases}
  0,
  & (\alpha,\beta)\in\N_0 \\
 \resY{\left(P_{\alpha,k}\exp{\phi_{\alpha+\beta}}\right)},\quad
  & (\alpha,\beta)\in -\N,\;\;k:=-(\alpha,\beta)-1
  \end{cases} 
  \end{align*}
  with $\pm$ depending on the parity of $(\alpha,\beta)$. In particular for $(\alpha,\beta)=-1$ resp. $(\alpha,\alpha)=2$
  \begin{align*}
  [\zem_{\alpha},\zem_{\beta}]_+
  &=\zem_{\alpha+\beta}\\
  [\zem_{\alpha},\zem_{-\alpha}]
  &=\yer_\alpha   
  \end{align*}
  We remark that if we consider instead of $\V_\Lambda$ the twist $\V^\epsilon_\Lambda$ by a $2$-cocycle $\epsilon(\alpha,\beta)$ as in Example~\ref{exm_latticeVOA}, which is the usual lattice vertex algebras with locality, then the anticommutators in the formula above turn into commutators.
\item We have for arbitrary $\alpha,\beta\in\C^r$, due to the grade shift,
\begin{align*}
[\yer_\alpha,\Y_m(v)]&=(\alpha,\beta)\;\Y_m(v),\\
\text{e.g.}\quad  [\yer_\alpha,\zem_\beta]&=(\alpha,\beta)\;\zem_\beta\hspace{1.8cm}
\end{align*}
\item Let $(\alpha,\beta)\in\Z$ and $u$ be any differential polynomial with $\N_0$-degree $|u|$. Then  
    $$(\alpha,\beta)\geq |u|\quad\Rightarrow\quad \zem_\alpha\;u\exp{\phi_\beta}=0 \hspace{2.2cm}$$
    More generally, the screening $\zem_\alpha$ 
    \emph{decreases} the $\N_0$-degree by  $(\alpha,\beta)+1$. Most extremely
    $$(\alpha,\beta)=|u|-1\quad\Rightarrow\quad\zem_\alpha \;u\exp{\phi_\beta}
    ={\langle \exp{\phi_\alpha},u\rangle_{-|u|}}
    \cdot\;\exp{\phi_{\alpha+\beta}}$$
  \end{enumerate}
\end{lemma}
\begin{proof} We first calculate in general
\begin{align*}
 \Y(\exp{\phi_\alpha})\exp{\phi_\beta}
  &=\langle\exp{\phi_\alpha,\phi_\beta}\rangle\cdot 
  \sum_{k\geq 0} \frac{z^k}{k!}\exp{\phi_\beta}\partial^k \exp{\phi_\alpha}
  = \sum_{k\geq 0}  
z^{(\alpha,\beta)+k}\cdot P_{k,\alpha}\exp{\phi_{\alpha+\beta}}
  \end{align*}
\begin{enumerate}[a)]
 \item This follows from the compatibility of the intertwining operator with the $\Lambda$-grading.
 \item The commutator with $\partial$ can be calculated quite general for the Hopf algebra expression; it matches the translation axiom of vertex operator algebras:
 \begin{align*}
 \Y(a,z)\;\partial.b
    &=\sum_{k\geq 0}\langle a^{(1)},\partial^{(1)}.b^{(1)}\rangle(z)
    \;\partial^{(2)}.b^{(2)}\;\frac{z^k}{k!}\partial^k.a^{(2)}\\ 
    &=\sum_{k\geq 0}\partial^{(3)}.\langle a^{(1)},b^{(1)}\rangle(z)\;\frac{z^k}{k!}
    \;\partial^{(2)}.\left(b^{(2)}\;S(\partial^{(2)}).\partial^k.a^{(2)}\right)\\
    &=\partial.\Y(a,z)b
    \;-\;\frac{\partial}{\partial z}
    \sum_{k\geq 0}\langle a^{(1)},b^{(1)}\rangle(z)\;\frac{z^k}{k!}
    \;\left(b^{(2)}\;\partial^k.a^{(2)}\right)
  \end{align*}
  where $S$ is the Hopf algebra antipode with $S(\partial)=-\partial$. The residue of the second term is zero for integral powers, but for non-integral powers we have
    \begin{align*}
    \res{\frac{\partial}{\partial z}z^m}
    &=m\;\res{z^{m-1}}
    = \frac{e^{2\pi\i m}-1}{2\pi\i}
  \end{align*}
  So the residue in question is
  \begin{align*}
  (\partial\zem_\alpha) u\exp{\phi_\beta}
  \;-\;(\zem_\alpha\partial )u\exp{\phi_\beta}
  &=\frac{e^{2\pi\i m}-1}{2\pi\i}\sum_{k\geq 0}\;\sum_{m\in(\alpha,\beta)+\Z}\langle \exp{\phi_\alpha},u^{(1)}\exp{\phi_\beta}\rangle_m\frac{1}{k!}
    \;\left(u^{(2)}\exp{\phi_\beta}\;\partial^k.\exp{\phi_\alpha}\right) \\
  &= \frac{e^{2\pi\i (\alpha,\beta)}-1}{2\pi\i}\left(\sum_{m\in\Z}\langle \exp{\phi_\alpha},u^{(1)}\rangle_m
    \;u^{(2)}\right)\exp{\phi_\beta}\;\sum_{k\geq 0}\frac{\partial^k}{k!}.\exp{\phi_\alpha} \\
  \end{align*}
 \item This follows from the commutator formula for vertex (super-)algebras  \cite{Kac98} (4.9.8)
  $$\zem_\alpha\zem_\beta=(-1)^{(\alpha,\beta)}\zem_\beta\zem_\alpha
  \;+\;\Y(\Y(\exp{\phi_\alpha})_{0}\exp{\phi_\beta})_{0}$$
  together with the previous formula for $\zem_\alpha{\exp{\phi_\beta}}$.
  It implies in particular for ${(\alpha,\beta)\in \N_0}$ that 
  $\Y(\exp{\phi_\alpha})_0\exp{\phi_\beta}=0$, and for $(\alpha,\beta)=-1$ 
  that $\Y(\exp{\phi_\alpha})_0\exp{\phi_\beta}=\exp{\phi_{\alpha+\beta}}$, 
  which yields in the first formula again a screening operator 
  $\zem_{\alpha+\beta}$.
  Note that one may derive similar formulae for other 
  values $(\alpha,\beta)\in-\N$, involving $\partial$ and the screening 
  operator $\zem_{\alpha+\beta}$.
 \item We proceed similarly as in c). Note that $\partial\phi_\alpha\in \V_0$ and $(0,\beta)\in 2\Z$, so the commutator formula applies, without sign. Then using the 
  explicit form of $\yer_\alpha$ as grading we get:
  \begin{align*}
    \Y(\phi_\alpha)_{0}\Y_m(v)-\Y_m(v)\Y(\phi_\alpha)_{0}
    &=\Y\left(\Y(\phi_\alpha)_{0}v\right)_m
    =\Y\left((\alpha,\beta)v\right)_m
  \end{align*}
 \item Since $\exp{\phi_\alpha}$ is grouplike the expression $\langle 
  \exp{\phi_\alpha},u\rangle$ is multiplicative in $u$, and it is anyways 
  $H$-linear, thus the $z$-dependence is $z^{-|u|}$, with coefficients product of factorials.
  Applying the definition of $\zem_\alpha$ and 
  $\langle\exp{\phi_\alpha},\exp{\phi_\beta}\rangle=z^{(\alpha,\beta)}$ yields:
  \begin{align*}
    \zem_\alpha \;u\exp{\phi_\beta}
    &=\sum_{k\geq 0} \langle \exp{\phi_\alpha}, 
  u^{(1)}\exp{\phi_\beta}\rangle_{-k-1}\;
  u^{(2)}\exp{\phi_\beta}\;\frac{1}{k!}\partial^k\exp{\phi_\alpha}\\
    &=\sum_{k\geq 0} \langle \exp{\phi_\alpha}, 
  u^{(1)}\rangle_{-k-1-(\alpha,\beta)}\;
  u^{(2)}\exp{\phi_\beta}\;\frac{1}{k!}\partial^k\exp{\phi_\alpha}
  \end{align*}
  So for non-vanishing terms we need $|u^{(1)}|=-k-1-(\alpha,\beta)$ together 
  with $k\geq 0$, which is only possible for $(\alpha,\beta)\not\in |u|-\N$. In 
  the extremal case $(\alpha,\beta)=|u|-1$ it can only be achieved for $k=0$ 
  and $|u^{(1)}|=|u|$, which means the summand $u^{(1)}\otimes u^{(2)}=u\otimes 
1$.
\end{enumerate}
\end{proof}

\subsection{Example: Root Lattices of Lie Algebras}\label{subsec_TrivialLevel}

We compute the algebra relations between local screening operators $\zem_{\alpha_i}$ in the root lattice of a Lie algebra $\g$, using the formulae in the previous section for $(\alpha_i,\alpha_j)\in\Z$. We find that these local screening operators fulfill the Lie algebra relations of the Borel part $\g^+$. 

This example already anticipates some aspects of what motivates our article in the non-local setting: If we replace root lattices by rescaled root lattices, then the Lie algebra action is replaced by a Nichols algebra action, due to non-locality. Ultimately, we get an action of the Borel part $u_q(\g)^+$ of a small quantum group.\\

%
%
%
%
%
%

\enlargethispage{1.5cm}
Let $\g$ be a finite-dimensional semisimple Lie algebra of rank $r$ with 
Cartan algebra~$\mathfrak{h}$, let $\Lambda=\Lambda_R$ be its root 
lattice and $\Phi$ its set of roots. Choose  positive simple roots 
$\alpha_1,...,\alpha_r$ and thereby a Borel subalgebra 
$\g^{\geq 0}=\mathfrak{h}\oplus \g^+$ and a set of positive roots 
$\Phi^+\subset\mathfrak{h}^*$.

\begin{example}
Assume for simplicity that $\g$ is simply-laced. Then the screening operators $\zem_{\alpha_i}$, together with the operators $\yer_\lambda$
give a representation of the Borel algebra $\g^{\geq 0}$ on the lattice vertex algebra $\V_\Lambda^\epsilon$ and on all modules.
In fact the vertex operator algebra $\V_\Lambda^\epsilon$ is isomorphic to  the vertex operator algebra associated to the    affine Lie algebra $\hat{\mathfrak{g}}_1$ at level one. This is called the \emph{Frenkel-Kac-Segal construction}, see \cite{Kac98} Sec. 5.6.
\end{example}

\noindent
The proof follows directly from the relations in Lemma \ref{lm_ScreeningProperties}. Set
$$E_\alpha:=\zem_\alpha,\;\alpha\in\Phi^+
\quad\qquad F_\alpha:=\zem_\alpha,\;\alpha\in\Phi^-
\quad\qquad H_\lambda:=\yer_\lambda,\;\lambda\in \mathfrak{h}^* $$

Then the important cases in Lemma \ref{lm_ScreeningProperties} were:
\begin{itemize}
 \item $(\alpha,\beta)=-1$, then 
  $\;\;[\zem_\alpha,\zem_\beta]_+=\zem_{\alpha+\beta}$ 
 \item $(\alpha,\beta)=0$, then 
  $\;\;[\zem_\alpha,\zem_\beta]=0$
   \item $(\alpha,\beta)=1$, then 
  $\;\;[\zem_\alpha,\zem_\beta]_+=0$
 \item $(\alpha,-\alpha)=-2,$ then 
 $[\zem_{\alpha},\zem_{-\alpha}]=\yer_\alpha$
 \item Moreover we have 
$\;\;[\yer_\lambda,\zem_\alpha]=(\lambda,\alpha)\;\zem_\alpha$
 \end{itemize}

 We repeat the remark in Lemma \ref{lm_ScreeningProperties} that the anticommutators turn into commutators if we consider the lattice vertex algebra $\V^\epsilon_\Lambda$ involving the $2$-cocycle $\epsilon$. We remark that if $\g$ is {not simply-laced}, then we encounter more cases and the relation between $E_\alpha,F_\alpha$ fails, so we can only find an action of the Borel algebra $\g^{\geq 0}$. \\

\section{A Non-Local Associativity Formula}\label{sec_Associativity}

The goal of this section is to compare the three expressions $\Y(a_1,z_1)\Y(a_2,z_2)$ and $\Y(a_2,z_2)\Y(a_1,z_1)$ and $\Y(\Y(a_1,z_1-z_2)a_2,z_2)$ for intertwining algebras. For vertex algebras, the notions of \emph{associativity} and \emph{locality} (or commutativity) state that these three functions are equal as meromorphic functions in  $z_1,z_2$, while as series their difference is given by the \emph{commutator formula}, see e.g. \cite{Kac98} Section 4.8 and formula 4.9.8.

We now derive for general intertwining operators constructed from Hopf algebra data in Definition  \ref{def_HopfIntertwiner} a formula that simplifies the concatenation of two intertwining operators (repeating some of our steps from \cite{Len07} Sec 4.2 and 4.3.4). Subsequently, this will yield a formula that simplifies the concatenation of two or more non-local screening operators: 
\enlargethispage{1.5cm}
\begin{align*}
  &\Y(a_1,z_1) \Y(a_2,z_2) v
  =\Y(a_1,z_1) \sum_{k\in\N_0} \langle a_2^{(1)},v^{(1)} \rangle(z_2) \cdot  v^{(2)}
  \cdot z_2^{k}\;\frac{\partial^{k}}{k!}.a_2^{(2)}\\
  &=\sum_{k_1\in\N_0}\sum_{k\in\N_0} \langle a_2^{(1)},v^{(1)} \rangle(z_2)
  \cdot z_2^{k} \cdot 
  \langle a_1^{(1)},
  \big(v^{(2)}\frac{\partial^{k}}{k!}.a_2^{(2)})^{(1)}\big\rangle(z_1) 
  \cdot \big(v^{(2)}\frac{\partial^{k}}{k!}.a_2^{(2)}\big)^{(2)}
  \cdot z_1^{k_1}\frac{\partial^{k_1}}{k_1!}.a_1^{(2)}\\
  &=\sum_{k_1\in\N_0}\sum_{k\in\N_0} \langle a_2^{(1)},v^{(1)} \rangle(z_2)\cdot 
  \langle a_1^{(1)},v^{(2)}\rangle(z_1)\\
  &\cdot 
  \sum_{k_{12}+k_2=k}
   z_2^{k}
  \langle a_1^{(2)},\frac{\partial^{k_{12}}}{k_{12}!}.a_2^{(2)}\rangle(z_1)\cdot z_1^{k_1}
  \cdot v^{(3)}\frac{\partial^{k_2}}{k_2!}.a_2^{(3)}
  \frac{\partial^{k_1}}{k_1!}.a_1^{(3)}\\
  &=\sum_{k_1,k_2\in\N_0} 
  \langle a_2^{(1)},v^{(1)} \rangle(z_2) \cdot 
  \langle a_1^{(1)},v^{(2)}\rangle(z_1)\cdot   z_2^{k_2}z_1^{k_1} \cdot
  v^{(3)}\frac{\partial^{k_2}}{k_2!}.a_2^{(3)} 
  \frac{\partial^{k_1}}{k_1!}.a_1^{(3)}\\
  &\cdot \left(\;\sum_{k_{12}\in\N_0} \frac{z_2^{k_{12}}}{k_{12}!}
  (-\frac{\partial}{\partial 
  z_1})^{k_{12}} \langle  a_1^{(2)},a_2^{(2)}\rangle(z_1) \;\right) \qquad (\star)
\end{align*}
We see in the second-to-last line the two independent intertwining operators $Y(a_2^{(1)},z_2) v^{(1)}$ and  $\Y(a_1^{(1)},z_1) v^{(2)}$, and in the last line a scalar function $(\star)$ in $z_1,z_2$, which contains the interplay between $a_1$ and $a_2$. If we assume for 
simplicity that the Hopf algebra $\V$ is commutative and cocommutative, as $\V_\Lambda$ is, then only this last term may prevent commutativity of $\Y(a_1,z_1),\Y(a_2,z_2)$ and respectively of the mode operators $\Y(a_1)_{n_1},\Y(a_2)_{n_2}$.	

\begin{remark}
	For vertex algebras with integral $z$-powers, we now discuss the role of~$(\star)$ in $\Y(a_1,z_1)\Y(a_2,z_2)$ and $\Y(a_2,z_2)\Y(a_1,z_1)$ in its usual analytic interpretation: 
If we abbreviate the value of the Hopf pairing by  $f(z_1)=\langle a_1,a_2\rangle(z_1)$ and $g(z_2)=\langle a_2,a_1\rangle(z_2)$, then on the level of meromorphic functions we can rewrite both $(\star)$ as a translation:
\begin{align*}
\sum_{k_{12}\in\N_0} \frac{z_2^{k_{12}}}{k_{12}!}
(-\frac{\partial}{\partial z_1})^{k_{12}}\; f(z_1)
  =e^{-z_2\frac{\partial}{\partial z_1}} f(z_1)= f(z_1-z_2)\\
  \sum_{k_{21}\in\N_0} \frac{z_1^{k_{21}}}{k_{21}!}
  (-\frac{\partial}{\partial z_2})^{k_{21}}\; g(z_2)
  =e^{-z_1\frac{\partial}{\partial z_2}} g(z_2)= g(z_2-z_1)
  \end{align*}
Hence as meromorphic functions we find that $\Y(a_1,z_1)\Y(a_2,z_2)=\pm \Y(a_2,z_2)\Y(a_1,z_1)$ holds, if additionally $g(z)=\pm f(-z)$, which means $\langle a_2,a_1\rangle(z)=\pm\langle a_1,a_2\rangle(-z)$.\\

However, on the level of series there is no equality, in fact the expressions $(\star)$ in $\Y(a_1,z_1)\Y(a_2,z_2)$ and $\Y(a_2,z_2)\Y(a_1,z_1)$ have different regions of convergence $|z_1|>|z_2|$ and $|z_2|>|z_1|$. As an example, we now calculate for $f(z)=g(-z)=z^{-1}$ the difference of these two expressions as a formal power series:
\begin{align*}
&\sum_{k_{12}\in \N_0} \frac{z_2^{k_{12}}}{k_{12}!}
(-\frac{\partial}{\partial z_1})^{k_{12}}z_1^{-1}
-\sum_{k_{21}\in\N_0} \frac{z_1^{k_{21}}}{k_{21}!}
(-\frac{\partial}{\partial z_2})^{k_{21}}(-z_2^{-1})\\
&=\sum_{k_{12}\in\N_0} z_2^{k_{12}}z_1^{-1-k_{12}}
+\sum_{k_{21}\in\N_0} z_1^{k_{21}}z_2^{-1-k_{21}}\\
&=\sum_{k\in\Z} z_2^{k}z_1^{-1-k}=:\delta(z_1-z_2)
\end{align*}
with $\delta(z_1-z_2)$ the delta function as a formal power series. Similarly, a polar term $f(z)=z^{-k+1}$ produces a formal $k$-th derivative ${\delta^{(k)}(z_1-z_2)}$. The coefficients of these delta functions in the commutator ${\Y(a_1,z_1)\Y(a_2,z_2)-\Y(a_2,z_2)\Y(a_1,z_1)}$ are the essential non-trivial information in a vertex algebra. In contrast, the functions $f(z)$ appearing in the present article for intertwining operators are multivalued or "non-local".

Physically, the appearance of these delta functions encodes the fact that particles only interact in the same point in space, or that processes only occur if momenta are conserved. 
\end{remark}

We now want to compute the formal residue of such products $\Y(a_1,z_1)\Y(a_2,z_2)$. We again use the notation $\langle-,-\rangle_m$ for the coefficient of $z^m,m\in\C$ in the Hopf pairing, and we may safely write $\sum_{m\in\C}$ because there are only finitely many contributions.

\begin{lemma}\label{lm_fractionalAsociativity}
The associativity formula above implies  the following \emph{composition formula} for residue operators associated to general intertwining operators as in Definition \ref{def_HopfIntertwiner}: 
\begin{align*}
&\resY(a_1)\resY(a_2) v\\
 &=\sum_{m_{1},m_{2},m_{12}}\sum_{k_1,k_2\in\N_0} 
\langle a_1^{(2)},a_1^{(2)}\rangle_{m_{12}}\langle a_2^{(1)},v^{(1)} \rangle_{m_2}
\langle a_1^{(1)},v^{(2)}\rangle_{m_1} \cdot v^{(3)}
\frac{\partial^{k_2}}{k_2!}.a_2^{(3)}\frac{\partial^{k_1}}{k_1!}.a_1^{(3)}\\
  &\cdot \underbrace{\sum_{k_{12}\in\N_0}\res{z_2^{m_2+k_2}
  \cdot z_2^{k_{12}}} 
  \res{z_1^{m_1+k_1}\cdot \frac{1}{k_{12}!}(-\frac{\partial}{\partial 
  z_1})^{k_{12}}z_1^{m_{12}}}}_{=:\;\Fp(m_1+k_1,m_2+k_2,m_{12})\;\in\C}
\end{align*}
and with similar calculations an expression for the other bracketing order
\begin{align*}
&\resY\left(\resY(a_1)a_2\right)v\\
 &=\sum_{m_{1},m_{2},m_{12}}\sum_{k_1,k_2\in\N_0} 
 \langle a_1^{(1)},a_2^{(1)}\rangle_{m_{12}}\langle a_2^{(2)},v^{(1)} \rangle_{m_2}
 \langle a_1^{(2)},v^{(2)}\rangle_{m_1} \cdot v^{(3)}
 \frac{\partial^{k_2}}{k_2!}.a_2^{(3)}\frac{\partial^{k_1}}{k_1!}.a_1^{(3)}\\
  &\cdot\underbrace{\sum_{k_{12}\in\N_0} 
  \res{z_{12}^{m_{12}}\cdot z_{12}^{k_{12}}} 
  \res{z_2^{m_2+k_2}\cdot 
  \frac{1}{k_{12}!}(+\frac{\partial}{\partial z_2})^{k_{12}}z_2^{m_{1}+k_1}}
  }_{=:\;\Fm(m_2+k_2,m_{12},m_1+k_1)\;\in\C}\\
\end{align*}
\end{lemma}
\noindent For $n$-fold compositions we respectively get by induction the following formula
\enlargethispage{0cm}
\begin{theorem}\label{thm_associativity}
\begin{align*}
  &\left(\prod_{i=1}^n \resY(a_i)\right)\; v\\
  &=\hspace{-.2cm}\sum_{(m_i,m_{ij})_{i,j}} \;\sum_{(k_i)_i \in \N_0^n}\;
  \prod_{1\leq i \leq n}\langle a_i^{(1)},v^{(i)}\rangle_{m_i}
  \prod_{1\leq i<j \leq n} \langle a_i^{(n-j+2)},a_j^{(n-i+1)}\rangle_{m_{ij}}
  \cdot v^{(n+1)}\prod_{i=n}^1 \frac{\partial^{k_i}}{k_i!}a_i^{(n+1)}\\
  &\cdot\underbrace{\sum_{(k_{ij})_{ij}\in \N_0^{n\choose 2}}\prod_{i} \res{z_i^{(m_i+k_i)+\sum_{i<j}(m_{ij}-k_{ij})+\sum_{j<i}k_{ji}}}\prod_{i<j} (- 1)^{k_{ij}}{m_{ij} \choose k_{ij}}}_{
  =:\;\Fp((m_i+k_i,m_{ij})_{ij})\;\in\C}
\end{align*} 
where the sum is taken over indices $k_{ij}\in \N_0$ for all $1\leq i<j\leq n$, and where the product $\prod_{i=n}^1$ means to be taken left-to-right in order of decreasing index. The expression depends on complex numbers $\Fp((m_i+k_i,m_{ij})_{ij})$, depending on $m_i,k_i,m_{ij}\in\C$ for $i\leq i<j\leq n$. We call them \emph{ quantum monodromy numbers} and study them in the next section.
\end{theorem}

\newpage

\section{The Quantum Symmetrizer Formula}\label{sec_QuantumSymmetrizer}

\noindent
\subsection{The Quantum Monodromy Numbers \texorpdfstring{$\Fpm$}{F} as Series }\label{subsec_QuantumMonodromyNumbers}

We now study in detail the structure constants $\Fpm((m_i,m_{ij})_{ij})\in\C$ that appear in Lemma \ref{lm_fractionalAsociativity} and Theorem~\ref{thm_associativity}. We first repeat the definition as series:

\begin{definition}\label{def_quantummonodromynumber}
	The \emph{quantum monodromy number} $\Fpm((m_i,m_{ij})_{ij})$ for a set of complex parameters $(m_i,m_{ij})_{ij}$ with ${1\leq i<j\leq n}$ is defined by 
	\begin{align*}
	\Fpm((m_i,m_{ij})_{ij})
	&:=\sum_{(k_{ij})_{ij}\in \N_0^{n\choose 2}}
	\prod_{i} \res{z_i^{m_i+\sum_{i<j}(m_{ij}-k_{ij})+\sum_{j<i}k_{ji}}}\prod_{i<j} (\pm 1)^{k_{ij}}{m_{ij} \choose k_{ij}}
	\end{align*}
	Assuming non-integrality $m_i+\sum_{i<j}m_{ij}\not\in \Z$ this is more explicitly:
	$$\Fpm((m_i,m_{ij})_{ij})=\sum_{(k_{ij})_{ij}\in \N_0^{n\choose 2}}
	\prod_{i}\frac{(e^{2\pi\i(m_i+\sum_{i<j} m_{ij})}-1)/2\pi\i}{
		1+m_i+\sum_{i<j}(m_{ij}-k_{ij})+\sum_{j<i}k_{ji}}
	\prod_{i<j} (\pm 1)^{k_{ij}}{m_{ij} \choose k_{ij}}$$
\end{definition}

Convergence of this series depending on $(m_{ij})_{ij}$ is rather subtle. A convergence condition that is sufficient for our purposes will be derived in Lemma \ref{lm_ConditionalConvergence} as an application of our main theorem. Now we prove the following crude convergence condition, which holds for an open subset of values for $(m_i,m_{ij})_{ij}$ and hence suffices for analytic~continuation:

\begin{lemma}\label{lm_firstconvergence}
	 If all $\Re(m_{ij})>0$, then the series $\Fpm((m_i,m_{ij})_{ij})$ is absolutely convergent.
\end{lemma}
\begin{proof}
	A crude bound for the series $\Fpm((m_i,m_{ij})_{ij})$ is  
	$$\sim\sum_{(k_{ij})_{ij}\in \N_0^{n\choose 2}}\prod_{i<j} \left\vert{m_{ij} \choose k_{ij}}\right\vert 
	=\prod_{i<j} \left(\sum_{k\geq 0}\left\vert{m_{ij} \choose k}\right\vert\right)$$
	To estimate each factor $\sum_{k\geq 0}\left|{m \choose k}\right|$ 
	we use the asymptotic ${m\choose k}\approx \frac{(-1)^k}{\Gamma(-m)}k^{-m-1}$ for $k\to \infty$ and fixed $m$, so the sum converges for  $\Re(m_{ij})>0$.
\end{proof}

If all $m_i,m_{ij}$ are equal, then the Selberg integral formula in Example~\ref{exm_Selberg}, together with our quantum symmetrizer formula in Theorem \ref{thm_QuantumSymmetrizer}, will give an expression for $\Fp$ in terms of Gamma functions, see Section \ref{subsec_EqualWeyl}. More generally, if all $m_i,m_{ij}$ are equal modulo $\Z,2\Z$, then there is an expression in terms of Jack polynomials. In this special case, quantum monodromy numbers already appears in our context in \cite{TW}.  

In the following we give an explicit formula for $n=1,2$. In general the author was unable to derive a simpler expression for $\Fpm$. An idea of the present article is to prove relations between the quantum monodromy numbers without actually calculating them.

\begin{example}\label{exm_F2series}
	For $n=2$ the series converges absolutely for $m_{12}>-2$, because the sum is taken over a single variable $k=k_{12}$ and the summand is of magnitude $k_{12}^{-2}k_{12}^{-1-m_{12}}$. For $m_{12}>-1$ we perform partial fraction decomposition and obtain an explicit formula:
	\begin{align*}
	\Fpm(m_1,m_2,m_{12})
	&=\frac{e^{2\pi\i\;m_2}-1}{2\pi\i}\cdot 
	\frac{e^{2\pi\i\;(m_1+m_{12})}-1}{2\pi\i}
	\sum_{k\in\N_0}\frac{(\pm 1)^k{m_{12} \choose k}}{(m_2+k+1)(m_1+m_{12}-k+1)} \\
	&=\frac{e^{2\pi\i\;m_2}-1}{2\pi\i}\frac{e^{2\pi\i\;(m_1+m_{12})}-1}{2\pi\i}\frac{1}{m_1+m_2+m_{12}+2}\\
	&\cdot \left(\sum_{k\in\N_0}\frac{(\pm 1)^k{m_{12} \choose k}}{m_2+k+1}
	-\sum_{k\in\N_0}\frac{(\pm 1)^k{m_{12} \choose k}}{-m_1-m_{12}+k-1} \right)\\
	\Fp(m_1,m_2,m_{12})&=\frac{e^{2\pi\i\;m_2}-1}{2\pi\i}\frac{e^{2\pi\i\;(m_1+m_{12})}-1}{2\pi\i}\frac{1}{m_1+m_2+m_{12}+2}\\
	&\cdot \left(\Beta(m_2+1,m_{12}+1)-\Beta(-m_1-m_{12}-1,m_{12}+1)
	\vphantom{\frac{\sin\pi m_1}{\sin\pi(m_1+m_{12})}}\right)\\
	&=\frac{e^{2\pi\i m_2}-1}{2\pi\i}\frac{e^{2\pi\i m_1+2\pi\i m_{12}}-1}{2\pi\i}
	\;\frac{1}{m_1+m_2+m_{12}+2}\\
	&\cdot \left(\Beta(m_2+1,m_{12}+1)+\frac{\sin\pi m_1}{\sin\pi(m_1+m_{12})}\Beta(m_1+1,m_{12}+1)\right)
	\end{align*}
	where $\Beta(x,y)=\frac{\Gamma(x)\Gamma(y)}{\Gamma(x+y)}$ is the Euler Beta function and in the last equation we use the Euler reflection formula. Alternatively, we recognize the series $\Fpm(m_1,m_2,m_{12})$ as the generalized hypergeometric function $_2F_3$ at $z=\mp1$, which is its boundary of convergence.\\
	
	The formula for $\Fpm(m_1,m_2,m_{12})$ has liftable singularities for certain integral parameter values, according to the cases in Definition \ref{def_res}. 
	We now list these cases:
	
	For $m_2,\;m_1+m_{12}\in \Z$ we get binomial coefficients, and in this sense and in this role the quantum monodromy numbers appear implicitly in standard vertex algebra theory.
	\begin{align*}
	&\Fpm(m_1,m_2,m_{12})\\
	&=\begin{cases}
	(\pm 1)^{m_2+1}{m_{12}\choose -m_2-1} & \text{if }m_1+m_2+m_{12}+2=0,\;-m_2-1\geq 0\\
	0,
	&\text{else}
	\end{cases}\\
	&=\begin{cases}
	(\pm 1)^{m_1+m_{12}+1}{m_{12}\choose m_1+m_{12}+1}
	&\text{if }m_1+m_2+m_{12}+2=0,\;m_1+m_{12}+1\geq 0\\
	0,
	&\text{else}
	\end{cases}
	\end{align*}
	\enlargethispage{1cm}
	This last equality appears implicitly in proofs of locality and the commutator formula (such as Theorem 4.6 in \cite{Kac98}). In other proofs emphasizing contour integrals (e.g. in \cite{FB04}) three different contour integrals over $z_1^{m_1}z_2^{m_2}z_{12}^{m_{12}}$ are compared, which relates to the contour integrals in Section \ref{subsec_Fintegral} in the case of single-valued functions. The equality itself will not generalize to our multi-valued setting. This relates to the fact, that we will no longer get commutator relations and Lie algebras.

	For $m_2\in\Z,\;m_1+m_{12}\not\in \Z$ and for $m_2\not\in\Z,\;m_1+m_{12}\in \Z$ we get respectively
	\begin{align*}
	&\Fpm(m_1,m_2,m_{12})\\
	&=\begin{cases}
	\frac{{(e^{2\pi\i\;(m_1+m_{12})}-1)}/{2\pi\i}}{m_1+m_2+m_{12}+2}(\pm 1)^{m_2+1}{m_{12}\choose -m_2-1},\quad
	&\text{if } -m_2-1\geq 0\\
	0,
	&\text{else}
	\end{cases}\\
	&\Fpm (m_1,m_2,m_{12})\\
	&=\begin{cases}
	\frac{{(e^{2\pi\i\;m_2}-1)}/{2\pi\i}}{m_1+m_2+m_{12}+2}(\pm 1)^{m_1+m_{12}+1}{m_{12}\choose 
		m_1+m_{12}+1},\quad
	&\text{if }m_1+m_2+m_{12}+2\geq 0 \\
	0,
	&\text{else}
	\end{cases}\\
	\end{align*}
\end{example}

\subsection{A Regularity Condition: Subpolar \texorpdfstring{$(m_{ij})_{ij}$}{mij}}\label{subsec_smallness}

The following somewhat subtle condition will appear frequently to ensure convergence of $n$-fold integrals and series. 

\begin{definition}\label{def_smallnessF}
	Let $I=\{1,\ldots,n\}$. A set of complex parameters $m_{ij}$ with $1\leq i<j\leq n$ is called \emph{subpolar}, if for any subset $J\subset I$ with $|J|\geq 2$ the following inequality holds:
	$$\sum_{i<j,\;i,j\in J} \mathfrak{Re}(m_{ij}) > -|J|+1$$
	A set of complex parameters $(m_{ij})_{ij}$ is called \emph{subpolar on intervals},~if the previous inequality holds for any subset  $J\subset I$ with $|J|>2$ that is an interval $J=[k,l]=\{k,\ldots,l\}$.
\end{definition}
\noindent
For our purposes, the condition can be guaranteed in the following setting: 
\begin{lemma}\label{lm_smallness}
	Suppose $\alpha_1,\ldots, \alpha_n$ are elements in a positive-definite euclidean vectorspace and denote the norm by $\Vert v\Vert^2=(v,v)$. 
	Then $m_{ij}:=(\alpha_i,\alpha_j)$ is subpolar, if all $\Vert\alpha_i\Vert^2<1$. It is moreover subpolar, if $\Vert\alpha_i\Vert^2\leq 1$ and $\alpha_i\neq -\alpha_j$ for all $i,j\in I$.  
\end{lemma}
\begin{proof}
	For every subset $J\subset I$ the assumptions of the Lemma hold for $J$, so without loss of generality we assume $J=\{1,\ldots, n\}$. Then 
	\begin{align*}
	\sum_{i<j} m_{ij} 
	& = \frac{1}{2} \sum_{i,j} m_{ij}-\frac{1}{2} \sum_{i} m_{ii}
	= \frac{1}{2} \Vert\sum_i\alpha_i\Vert^2-\frac{1}{2} \sum_{i} \Vert\alpha_i\Vert^2
	> 0-n/2 \geq -n+1
	\end{align*}
	If all $\Vert\alpha_i\Vert^2=1$ then we have to ensure that the inequality holds strictly: For $n>2$ the last inequality holds strictly. For $n=2$ the second-to-last inequality holds strictly unless $\Vert\sum_i\alpha_i\Vert^2=0$, which only happens for $\alpha_1=-\alpha_2$. This case is indeed not subpolar.
\end{proof}

\noindent
If $(m_{ij})_{ij}$ is subpolar, then this has the following crucial implication:
\begin{lemma}\label{lm_Christian2}
	Suppose that $m_{ij}\in\C$ for $1\leq i<j\leq n$ is subpolar on intervals, then	 $\vert\prod_{1\leq i<j\leq n} (z_i-z_j)^{m_{ij}}\vert$, for any chosen branch of the logarithm, is bounded by a linear combination of terms of the form $\prod_{1\leq i<n} |z_i-z_{i+1}|^{\hat{{m}}_{i,i+1}}$, with real numbers $\hat{{m}}_{i,i+1}>-1$.
\end{lemma}

This lemma follows inductively from the following combinatorial construction,\footnote{Many thanks to Christian Reiher for providing this construction.} which explicitly gives the different values of $\hat{{m}}_{i,i+1}$ appearing in the bound.

\begin{lemma}\label{lm_Christian1}
	Let $m_{ij}\in\R$ for $1\leq i<j\leq n$ be subpolar on intervals \ref{def_smallnessF}. We define
	\begin{align*}
	P&:=\{j\in [3, n]\mid m_{1j}\geq 0\}\\
	N&:=\{j\in [3, n]\mid m_{1j}< 0\}
	\end{align*}
	Then for any $K\subset P$ there exist $m_{1j}^K\in \R$ for all $j\in N$ with the following properties:
	\begin{enumerate}[a)]
		\item $m_{1j}\leq m_{1j}^K\leq 0$ for all $j\in N$.
		\item The following numbers $\hat{m}_{ij}\in\R$ for  $2\leq i,j\leq n$ are again subpolar on intervals
		$$\hat{m}_{ij}^K:=\begin{cases}
		m_{ij}, &\text{for }i\geq 3 \\
		m_{2j}, &\text{for }j\in K \\
		m_{2j}+m_{1j}, \quad &\text{for }j\in P\backslash K \\
		m_{2j}+m_{1j}^K, \quad&\text{for }j\in N
		\end{cases}$$
		\item $\hat{m}_{12}^K>-1$ holds for the quantity 
		$$\hat{m}_{12}^K:=\sum_{1\leq i<j\leq n} m_{ij}-\sum_{2\leq i<j\leq n} \hat{m}_{ij}^K
		\;=\;m_{12}+\sum_{j\in K} m_{1j}+\sum_{j\in N} (m_{1j}-m_{1j}^K)$$
	\end{enumerate}
\end{lemma}
\begin{proof}
	By assumption $m_{ij}$ is subpolar on intervals, so for any interval $[r,s]$ we have 
	$$f(r,s)> r-s\qquad \text{ for }\quad f(r,s):=\sum_{r\leq i<j\leq s} m_{ij}$$
	Choose $\epsilon>0$ with $f(r,s)> r-s+\epsilon$ for all $r,s$.\\

	We construct $m_{1j}^K$ recursively, so that a) and b) hold: For each index  $j\in N$, suppose that $m_{1k}^K$ is already defined for all $k>j$ with $k\in N$, then define $m_{1j}^K\geq m_{1j}$ as the minimal value, such that 
	$$\forall s\in[j,n]\qquad 
	f(2,s)
	+\sum_{k\in [3,s]\cap (P\backslash K)} m_{1k} 
	+\sum_{k\in [j,s]\cap N} m_{1k}^K 
	\geq 2-s+\epsilon\qquad (\star)$$
	
	We prove property a): Suppose on the contrary that this minimal value $m_{1j}^K$ is positive (and $m_{1k}^K$ nonpositive for all $k>j$ with $k\in N$), then there must exist an $s\in[j,n]$ with 
	$$f(2,s)
	+\sum_{k\in [3,s]\cap (P\backslash K)} m_{1k} 
	+\sum_{k\in [j+1,s]\cap N} m_{1k}^K 
	<(2-s)+\epsilon$$
	But $f(2,s)>2-s+\epsilon$ and all $m_{1k},k\in P$ are nonnegative, so $[j+1,s]\cap N$ is nonempty and for the smallest element $j'\in [j+1,s]\cap N$ the defining property $(\star)$ of $m_{1j'}^K$ is contradicted.\\
	
	We prove property b): Set $\hat{f}(r,s)=\sum_{r\leq i<j\leq s} \hat{m}_{ij}^K$ for $2\leq r<s\leq n$. Then for $r>2$ we have $\hat{f}(r,s)=f(r,s)>r-s$. For $r=2$ we have:
	$$\hat{f}(2,s)=f(2,s)+\sum_{k\in [3,s]\cap(P\backslash K)} m_{1k}+\sum_{k\in [3,s]\cap N} m_{1k}^K$$
	For $[3,s]\cap N$ nonempty again the defining property  $(\star)$ of $m_{1k}^K$, for $k$ the minimum of $[3,s]\cap N$, ensures $\hat{f}(2,s)\geq (2-s)+\epsilon$ and for $[3,s]\cap N$ empty $\hat{f}(2,s)\geq f(2,s)\geq 2-s+\epsilon$.\\
	
	We prove property c): Assume on the contrary that 
	$$\hat{m}_{12}^K=m_{12}+\sum_{j\in K} m_{1j}+\sum_{j\in N} (m_{1j}-m_{1j}^K) \leq -1$$
	Since $m_{12}>-1$ and all $m_{1j},j\in P$ are nonnegative we can choose the largest index $j\in N$ with $m_{1j}-m_{1j}^K<0$. Since we defined $m_{1j}^K\geq m_{1j}$ to be the minimal choice fulfilling $(\star)$, one of the inequalities holds strictly, i.e., there exists an $s\in[j,n]$ with 
	$$f(2,s)
	+\sum_{k\in [3,s]\cap (P\backslash K)} m_{1k} 
	+\sum_{k\in [j,s]\cap N} m_{1k}^K 
	=(2-s)+\epsilon\qquad (\star\star)$$
	Claim: For all $j'<j$ with $j'\in N$ holds $m_{1j'}^K=0$, otherwise $(\star)$ fails for $\hat{m}_{1j'}^K$.
	\newline
	Using this claim and our choice of $j$, we finally conclude:  
	\begin{align*}
	\hat{m}_{12}^K 
	&>\hat{m}_{12}^K+(-f(1,s)+(1-s)+\epsilon) \\
	&=\hat{m}_{12}^K+\left(- f(2,s)+(2-s)+\epsilon\right) -\sum_{k\in [2,s]} m_{1k}-1\\
	&=\left(m_{12}+\sum_{k\in K} m_{1k}+\sum_{k\in [3,s]\cap N} (m_{1k}-m_{1k}^K)
	+\underbrace{\sum_{k\in [s+1,n]\cap N} (m_{1k}-m_{1k}^K)}_{=0 \text{ by choice of $j$ and $s\geq j$}}\right)\\
	(\star\star)\quad&+\left(\sum_{k\in [3,s]\cap (P\backslash K)} m_{1k} +\sum_{k\in [j,s]\cap N} m_{1k}^K\right)\\
	&-\left(m_{12}+\sum_{k\in [3,s]\cap K} m_{1k}+\sum_{k\in [3,s]\cap (P\backslash K)} m_{1k}+\sum_{k\in [3,s]\cap N} m_{1k}\right)-1 \\
	&=\underbrace{\sum_{k\in [s+1,n]\cap K} m_{1k}}_{\geq 0}+\underbrace{\sum_{k\in [3,j)\cap N} (-m_{1k}^K)}_{=0\text{ by claim}}-1 \quad \geq -1 \quad \qedhere
	\end{align*}
\end{proof}
\newpage
\noindent
Using this combinatorial construction we find the inductive bound as follows:

\begin{lemma}\label{lm_Christian2step}
	Let $m_{ij}\in\R$ for $1\leq i<j\leq n$ be subpolar on intervals and $n\geq 3$. Then
	$$\prod_{1\leq i<j\leq n} |z_i-z_j|^{m_{ij}} 
	\leq C\sum_{K\subset P} |z_1-z_2|^{\hat{m}_{12}^K} \prod_{2\leq i<j\leq n} |z_i-z_j|^{\hat{m}_{ij}^K}  $$ 
	where $\hat{m}_{ij}^K,\hat{m}_{12}^K$ and $P\subset I$ were given in Lemma \ref{lm_Christian1}, and we recall that the construction was made in such a way that $\hat{m}_{ij}^K$ is again subpolar on intervals  and that $\hat{m}_{12}^K>-1$.
\end{lemma}	
\begin{proof}
	The previous Lemma \ref{lm_Christian1} provides for $m_{ij}$, subpolar on intervals, and for any choice $K\subset P$ a set of numbers $m_{1j}^K,j\in N$ and a new set of parameters $\hat{m}_{ij}^K$, again subpolar on intervals by \ref{lm_Christian1} b). We now get an estimate for $\prod_{1\leq i<j\leq n} |z_i-z_j|^{m_{ij}}$: \\
	\noindent
	For $j\in P$ we estimate by the following additive decomposition
	\begin{align*}|z_1-z_j|^{m_{1j}}
	&\leq 2^{m_{1j}}\left(\;\frac{|z_1-z_2|}{2}+\frac{|z_2-z_j|}{2} \;\right)^{m_{1j}}
	\leq 2^{m_{1j}}\left(|z_1-z_2|^{m_{1j}}+|z_2-z_j|^{m_{1j}}\vphantom{\frac{z_1-z_j}{2}}\right)
	\end{align*}
	Multiplying out their product over all $j\in P$ gives a bound for this product by $2^{|P|}$ summands, which we index by subsets $K\subset P$.
	\begin{align*}
	\prod_{j\in P} |z_1-z_j|^{m_{1j}}
	&\leq 2^{\sum_{j\in P}m_{1j}}\sum_{K\subset P} \left(\prod_{j\in K} |z_1-z_2|^{m_{1j}}\prod_{j\in P\backslash K}|z_2-z_j|^{m_{1j}}\right) \\
	\intertext{For $j\in N$ we estimate by the following multiplicative decomposition}
	|z_1-z_j|^{m_{1j}}
	&\leq \mathrm{max}(z_1-z_2,z_2-z_j)^{m_{1j}}
	\leq |z_1-z_2|^{m_{1j}-m_{1j}^K} |z_2-z_j|^{m_{1j}^K} \\
	\intertext{where we choose depending on $K$ the decomposition $m_{1j}=m_{1j}^K+(m_{1j}-m_{1j}^K)$ provided by Lemma \ref{lm_Christian1}, where $m_{1j}^K,\;m_{1j}-m_{1j}^K\leq 0$ by Lemma \ref{lm_Christian1} a). Combining $P$ and $N$ gives}
	\prod_{j\in [2,n]} |z_1-z_j|^{m_{1j}}
	&=|z_1-z_2|^{m_{12}}\prod_{j\in P}|z_1-z_j|^{m_{1j}}\prod_{j\in N}|z_1-z_j|^{m_{1j}}\\
	&\leq 2^{\sum_{j\in P}m_{1j}}\sum_{K\subset P} |z_1-z_2|^{\hat{m}_{12}^K}
	\prod_{j\in P\backslash K}|z_2-z_j|^{m_{1j}}\prod_{j\in N}|z_2-z_j|^{m_{1j}^K} \\
	\intertext{for $\hat{m}^K_{12}=m_{12}+\sum_{j\in K} m_{1j}+\sum_{j\in N} (m_{1j}-m_{1j}^K)$ with $\hat{m}^K_{12}>-1$ by Lemma \ref{lm_Christian1} c).}
	\prod_{1\leq i<j\leq n} |z_i-z_j|^{m_{ij}} 
	&\leq 2^{\sum_{j\in P}m_{1j}}\sum_{K\subset P} |z_1-z_2|^{\hat{m}_{12}^K} \prod_{2\leq i<j\leq n} |z_i-z_j|^{\hat{m}_{ij}^K} \qedhere
	\end{align*}
\end{proof}	
\noindent
By induction, this concludes the proof of Lemma \ref{lm_smallness} for $m_{ij}\in\R$. For $m_{ij}\in\C$ we use
$$\vert\prod_{1\leq i<j\leq n} (z_i-z_j)^{m_{ij}}\vert\leq \prod_{1\leq i<j\leq n} \vert z_i-z_j\vert^{\mathfrak{Re}(m_{ij})} \cdot (2\pi)^{-\mathfrak{Im}(m_{ij})} $$
and note that $(\mathfrak{Re}(m_{ij}))_{ij}$ is subpolar on intervals iff $(m_{ij})_{ij}$ is subpolar on intervals.\qed

\subsection{The Generalized Selberg Integrals \texorpdfstring{$\Sel$}{Sel}}\label{subsec_Selberg}

\begin{definition}\label{def_Selberg}
	We define the \emph{generalized Selberg integral} for a set of complex para\-meters $(m_i,\bar{m}_i,m_{ij})_{ij}$ with $1\leq i<j\leq n$  by the following $n$-fold integral:
	\begin{align*}
	\Sel((m_i,\bar{m}_i,m_{ij})_{ij})
	&:=\idotsint\displaylimits_{1> z_1>\ldots >z_{n}> 0}
	\prod_i z_i^{m_i}\prod_i (1-z_i)^{\bar{m}_i}\prod_{i<j} (z_i-z_j)^{m_{ij}} \;\d z_1\cdots \d z_{n}
	\end{align*}
\end{definition}
\noindent
\begin{lemma}\label{lm_SelbergConvergence}
	The integral $\Sel((m_i,\bar{m}_i,m_{ij})_{i,j})$ is absolutely convergent, provided
	\begin{itemize}
	\item $m_{ij}$ is subpolar on intervals: For all intervals $J\subset I$ with $|J|\geq 2$ holds  
		$$\sum_{i<j,\;i,j\in J} m_{ij} > -|J|+1$$
	\item For the intervals $J=[k,\ldots, n]$ holds moreover 
		$$\sum_{i<j,\;i,j\in J} m_{ij} + \sum_{i\in J} m_i > -|J|$$
	\item For the intervals $J=[1,\ldots,k]$ holds moreover
	$$\sum_{i<j,\;i,j\in J} m_{ij} + \sum_{i\in J} \bar{m}_i > -|J|$$
\end{itemize}	
Moreover, for fixed real $(m_{ij})_{ij}$ and real $m_i\to \infty$ we have the asymptotic bound
$$\Sel((m_i,0,m_{ij})_{i,j})
\leq C' \left(\sum_{1\leq i\leq n} m_i\right)^{-n+\sum_{1\leq i<j\leq j}m_{ij}}$$
\end{lemma}
\begin{proof}
	This is a direct consequence of Lemma \ref{lm_Christian2}, if we introduce additional variables $z_0$ and $z_{n+1}$, which we then set to $1$ and $0$. Hence we now consider a new index set $I^{ext}=\{0,1,\ldots,n,n+1\}$ and define
	$$m_{ij}^{ext}=\begin{cases}
		m_{ij},\quad &0< i<j< n+1 \\
		m_i,\qquad   &0<i<j=n+1 \\
		\bar{m}_{j},\quad &0=i<j<n+1 \\
		-|J|-1-\sum_{1\leq i<j\leq n} m_{ij} - \sum_{1\leq i\leq n} m_i - \sum_{1\leq j\leq n} \bar{m}_j,
		\qquad &0=i,\;\;\;j=n+1 
	\end{cases}$$
	Then the assumptions on $(m_{ij})_{ij}$ are equivalent to $(m_{ij}^{ext})_{ij}$ being subpolar. By Lemma~\ref{lm_Christian2}
	the function $\prod_{0\leq i<j\leq n+1} (z_i-z_j)^{m_{ij}}$ is bounded by sums of terms $\prod_{0\leq i<n+1} |z_i-z_{i+1}|^{\hat{{m}}^{ext}_{i,i+1}}$, where all $\hat{{m}}_{i,i+1}^{ext}>-1$. In particular, after setting $z_0=1,z_{n+1}=0$ the integrand is bound by sums of terms $|1-z_1|^{\hat{{m}}^{ext}_{0,1}}\left(\prod_{1\leq i<n}|z_i-z_{i+1}|^{\hat{{m}}^{ext}_{i,i+1}}\right)|z_n|^{\hat{{m}}^{ext}_{n,n+1}}$. Now, we evaluate the rightmost integral of the bound by substituting $z_n'=z_n/z_{n-1}$, then:
	
	\begin{align*}
	\int_0^{z_{n-1}}\mathrm{d}z_{n}\;(z_{n-1}-z_{n})^{\hat{{m}}^{ext}_{n-1,n}} z_n^{\hat{{m}}^{ext}_{n,n+1}}
	&=z_{n-1}^{1+\hat{{m}}^{ext}_{n-1,n}+\hat{{m}}^{ext}_{n,n+1}}\int_0^{1}\mathrm{d}z_{n}'\;(1-z_{n}')^{\hat{{m}}^{ext}_{n-1,n}} z_n'^{\hat{{m}}^{ext}_{n,n+1}}\\
	&=z_{n-1}^{1+\hat{{m}}^{ext}_{n-1,n}+\hat{{m}}^{ext}_{n,n+1}}\cdot \Beta(\hat{{m}}^{ext}_{n-1,n}+1,\hat{{m}}^{ext}_{n,n+1}+1)
	\end{align*}
	where the Beta integral is convergent for $\hat{{m}}^{ext}_{n-1,n},\hat{{m}}^{ext}_{n,n+1}>-1$. Now we proceed inductively, as the new rightmost variable $z_{n-1}$ has again exponent ${1+\hat{{m}}^{ext}_{n-1,n}+\hat{{m}}^{ext}_{n,n+1}>-1}$.\\
	
	We finally prove the explicit bound for real parameters $(m_i,0,m_{ij})_{ij}$: In the induction in Lemma \ref{lm_Christian1} the term $m_i=m_{i,n+1}^{ext}$ enters only in $\hat{{m}}^{ext}_{n,n+1}$, and for all $m_i>0$ we have $\hat{{m}}^{ext}_{n,n+1}=\sum_{i\in K} m_i\leq \sum_{1\leq i\leq n} m_i$. We also observe that $\sum_{1\leq i<n}\hat{{m}}^{ext}_{i,i+1}\geq\sum_{1\leq i<j\leq n} m_{ij}$. Now, the previous paragraph gives an explicit bound with a Beta function, which for $x$ fixed and $y\to \infty$, is asymptotically $B(x+1,y+1)\sim y^{-1-x}$, so overall we have
 	$\sim \prod_{1\leq i<n}(\hat{{m}}^{ext}_{n,n+1})^{-1-\hat{{m}}^{ext}_{i,i+1}}$. Using the two inequalities above proves the asserted bound.
\end{proof}

\begin{example}\label{exm_Selberg}
	For $n=1$ the generalized Selberg integral is the Euler Beta integral
	\begin{align*}
	\Sel(m_1,\bar{m}_1)&=\Beta(m_1+1,\bar{m}_1+1)\\
	\Sel(m_1,m_2,0,0,m_{12})&=\frac{1}{2+m_1+m_2+m_{12}}  \Beta(m_2+1,m_{12}+1)
	\end{align*}
	More generally, substitution of $z_i'=z_i/z_1$ for $i\neq 1$ and integration over $z_1$ gives
	$$\Sel((m_i,0,m_{ij})_{ij})=\frac{1}{n+\sum_i m_i+\sum_{i<j}m_{ij}} \Sel((m_i,{m}_{1i},m_{ij})_{i,j\neq 1})$$
	For equal parameters $m_i=\alpha-1,\;\bar{m}_i=\beta-1,\;m_{ij}=2\gamma$ the Selberg integral \cite{Sel44} gives 
	$$\Sel(\alpha-1,\beta-1,2\gamma)
	=\frac{1}{n!}\prod_{j=0}^{n-1} \frac{\Gamma(\alpha+j\gamma)\Gamma(\beta+j\gamma)
		\Gamma(1+(j+1)\gamma)}{\Gamma(\alpha+\beta+(n+j-1)\gamma)\Gamma(1+\gamma)}$$
	For all $m_{ij}$ equal and all $m_i,\bar{m}_i$ equal modulo $\Z$, the symmetrization of the generalized Selberg integral can be computed as follows, which was used in our context in  \cite{TW}:
	For any partition with up to $n$ parts
  $\lambda=(\lambda_1,\ldots,\lambda_n)$  and for  $P_{\lambda}^{(1/\gamma)}(z_1,\ldots, z_n)$ the \emph{Jack~polynomial}, we have \emph{Kadell's integral}, see \cite{Kad97} and \cite{Mac95} Sec. VI/10:
	
	\enlargethispage{.5cm}
	\begin{align*}
	&\idotsint\displaylimits_{1> z_1>\ldots >z_{n}> 0} 
	P_{\lambda}^{(1/\gamma)}(z_1,\ldots, z_n)\prod_i z_i^{\alpha-1}\prod_i (1-z_i)^{\beta-1}\prod_{i<j} (z_i-z_j)^{2\gamma} \;\d z_1\cdots \d z_{n}\\
	&=\prod_{1\leq i<j\leq n}\frac{\Gamma((j-i+1)\gamma+\lambda_i-\lambda_j)}
			{\Gamma((j-i)\gamma+\lambda_i-\lambda_j)} 
			\prod_{j=0}^{n-1}
			\frac{\Gamma(\alpha+j\gamma+\lambda_{n-j})\Gamma(\beta+j\gamma)}
			{\Gamma(\alpha+\beta+(n+j-1)\gamma+\lambda_{n-j})}\\
			\end{align*}

\end{example}


\begin{problem}
	Find a full analytic continuation for $\Sel$ and determine the set of poles. 
\end{problem}

\subsection{The Reduced Quantum Monodromy Numbers \texorpdfstring{$\rFpm$}{tildeF}}\label{subsec_redF}

\begin{definition}\label{def_reducedcontourintegral}
	The \emph{reduced quantum monodromy number} $\rFpm((m_i,m_{ij})_{ij})$ for complex parameters $(m_i,m_{ij})_{ij}$ with $1\leq i<j\leq n$ is defined by the following integral:
	$$\rFpm((m_i,m_{ij})_{ij}):=\frac{1}{(2\pi\i)^n}\int_{\triangle}
	\prod_i z_i^{m_i}\prod_{i<j} (z_i\pm z_j)^{m_{ij}}\;\d z_1\cdots \d z_n$$ 
	where we integrate over the following image of an $n$-simplex in the multivalued covering
	$$\triangle:=\{(e^{2\pi\i t_1},\ldots,e^{2\pi\i t_n}) \mid 0<t_1<\cdots <t_n<1\}$$
\end{definition}
\noindent
As a direct consequence of Lemma \ref{lm_Christian2} we find
\begin{lemma}\label{lm_rFpConvergence}
	The integral $\rFp((m_i,m_{ij})_{ij})$ converges absolutely, provided that $(m_{ij})_{ij}$ is subpolar on intervals as in  Definition \ref{def_smallnessF}. 	
\end{lemma}
\begin{lemma}\label{lm_rFSelberg}
The integral $\rFp((m_i,m_{ij})_{ij})$ can be expressed in terms of Selberg integrals
\begin{align*}
\rFp((m_i,m_{ij})_{ij})
&:=\frac{1}{(2\pi\i)^n}\sum_{k=0}^n (-1)^{n-k}\left(\prod_{i=n-k+1}^n e^{2\pi\i\;m_i}\right)\sum_{\eta\in \S_{n-k,\overline{k}}}\left(\prod_{i<j,\;\eta(i)>\eta(j)}e^{\pi\i\;m_{ij}}\right)\\ 
&\cdot 
\Sel((m_{\eta^{-1}(i)},0,m_{\eta^{-1}(i)\eta^{-1}(j)})_{ij})
\end{align*}
where we define the following slight variation of familiar $(n-k,k)$-\emph{shuffles} 
$$\S_{n-k,\overline{k}}:=\{\eta\in \S_n \mid \forall_{i<j\leq n-k}\;\eta(i)<\eta(j)\text{ and } \forall_{n-k<i<j}\;\eta(i)>\eta(j)\}
\qquad |\S_{n-k,\overline{k}}|={n \choose k}$$
and where we assume that the assumptions in Lemma \ref{lm_SelbergConvergence} on $(m_{i},m_{ij})_{ij}$ hold for all such permutations of the index set $I$, so that $\Sel((m_{\eta^{-1}(i)},0,m_{\eta^{-1}(i)\eta^{-1}(j)})_{ij})$ exists.	
\end{lemma}
\begin{proof}
	We deform the integration contour $\triangle$ to the real line and a small circle around $z=0$. Thereby our integration domain decomposes into $2^n=\sum_{k=0}^n {n \choose k}$ pieces, namely $n+1$ cases, respectively, with $t_1,\ldots, t_{n-k}\approx 0$ and $t_{n-k+1},\ldots, t_n\approx 1$, so on the real axis we have 
	$z_1>\cdots > z_{n-k}$ and $z_{n-k+1}<\cdots < z_n$ and in each case all ${n \choose k}$ subcases how the two subsets are shuffled against each other.
	\begin{center}
		\begin{center}
			\includegraphics[scale=.74]{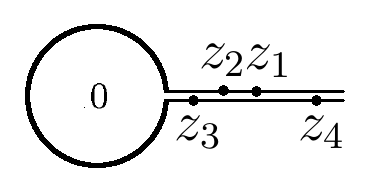}
		\end{center}
		$\qquad \text{Example: }k=2,\;\eta=(1234)\in\S_{2,\overline{2}}$
	\end{center}
	The assumptions of Lemma \ref{lm_SelbergConvergence} ensure that the integration over the small circle vanishes in the limit. For each $k\in\{1,\ldots, n\}$ and $\eta\in \S_{n-k,\overline{k}}$ the integration over the domain and branch discussed above $1>z_{\eta^{-1}(1)}>\cdots >z_{\eta^{-1}(n)}>0$ gives each a Selberg integral: 

	\begin{align*}
	&\rFp((m_i,m_{ij})_{ij})((m_i,m_{ij})_{ij})\\
	&=\sum_{k=0}^n (-1)^{n-k} s'(k)\sum_{\eta\in \S_{n-k,\overline{k}}}s''(k)\cdot 
	\frac{1}{(2\pi\i)^n} \Sel((m_{\eta^{-1}(i)},0,m_{\eta^{-1}(i)\eta^{-1}(j)})_{ij})
	\end{align*}
	Here the factor $(-1)^{n-k}$ is from reversing the direction of integration, the branch of the multivalued function $\prod_i z_i^{m_i}$ gives a factor $s'(k)=\prod_{i=n-k+1}^n e^{2\pi\i\;m_i}$ and the branch of the multivalued function $\prod_{i<j} (z_i-z_j)^{m_{ij}}$ gives a factor $s''(\eta)=\prod_{i<j,\;\eta(i)>\eta(j)}e^{\pi\i\;m_{ij}}$.
\end{proof}

\begin{example}\label{exm_rF2}
    For $n=2$ this formula reads explicitly:
		\begin{align*}
		&\rFp(m_1,m_2,m_{12})\\
		&=\frac{1}{(2\pi\i)^2}\left(\vphantom{x^{x^{x^x}}}
		e^{2\pi\i(m_1+m_2)+\pi\i m_{12}}\;\Sel(m_2,m_1,0,0,m_{12})
		-e^{2\pi\i m_2}\;\Sel(m_1,m_2,0,0,m_{12})\right.\\
		&\left.\vphantom{x^{x^{x^x}}}\hspace{2cm}-e^{2\pi\i m_2+\pi\i m_{12}}\;\Sel(m_2,m_1,0,0,m_{12})
		+e^{0}\;\Sel(m_1,m_2,0,0,m_{12})\right)\\
		&=\frac{1}{(2\pi\i)^2}\left(1-e^{2\pi\i m_2}\right)
		\frac{\Beta(m_2+1,m_{12}+1)}{m_1+m_2+m_{12}+2}\\
		&-\frac{1}{(2\pi\i)^2}e^{2\pi\i m_2+\pi\i m_{12}}\left(1-e^{2\pi\i m_1}\right)
		\frac{\Beta(m_1+1,m_{12}+1)}{m_1+m_2+m_{12}+2}
		\end{align*}
	We observe that, as expected, the poles for $m_i\in -\N$ in the Selberg integrals (here the Beta functions) are no longer poles of $\rFp(m_1,m_2,m_{12})$. In the formula we see that these poles are suppressed by the exponential prefactors resp. by canceling poles. 
\end{example}

\subsection{The Quantum Monodromy Numbers \texorpdfstring{$\Fpm$}{F} as Integrals}\label{subsec_Fintegral}

We now give an analytic continuation of the series $\Fpm((m_i,m_{ij})_{ij})$ in Definition \ref{def_quantummonodromynumber} in terms of contour integrals: For real numbers $r_1>r_2>\ldots>r_n>0$ consider the following unique lift of the $n$-torus with radii $r_1,\ldots, r_n$ to the multivalued cover  
$$\Box_{r_1,\ldots r_n}:=\left\{(r_1e^{2\pi\i t_1},\ldots,r_ne^{2\pi\i t_n}) \mid 0<t_1,\cdots,t_n<1\right\}$$
and consider the following contour integral, which is well-defined for all values of $(m_i,m_{ij})_{ij}$  
$$\Fpm^{r_1,\ldots,r_n}((m_i,m_{ij})_{ij}):=\frac{1}{(2\pi\i)^n}\int_{\Box_{r_1,\ldots r_n}}
\prod_i z_i^{m_i}\prod_{i<j} (z_i\pm z_j)^{m_{ij}}\;\d z_1\cdots \d z_n$$ 	

\enlargethispage{1cm}
\begin{definition}\label{def_contourintegral}
	We define $\Fpm((m_i,m_{ij})_{ij})$ as the limit of  $\Fpm^{r_1,\ldots,r_n}((m_i,m_{ij})_{ij})$ to $r_1=\cdots=r_n=1$. With some caution one could write 
	$$\Fpm((m_i,m_{ij})_{ij})=\frac{1}{(2\pi\i)^n}\int_{\Box}
	\prod_i z_i^{m_i}\prod_{i<j} (z_i\pm z_j)^{m_{ij}}\;\d z_1\cdots \d z_n$$ 
	for a certain lift $\Box$ of the $n$-torus with all radii equal to $1$ to the multivalued covering, on which the integrand has singularities with branch points. 
\end{definition}
The limit certainly exists if all $\Re(m_i),\Re(m_{ij})\geq 0$. It will be a direct consequence of our quantum symmetrizer formula \ref{thm_QuantumSymmetrizer} that $\Fp((m_i,m_{ij})_{ij})$ exists if $(m_{ij})_{ij}$ is subpolar.

\begin{lemma}\label{lm_contourequalsseries}
	Let us assume that all $\Re(m_i),\Re(m_{ij})\geq 0$, where we already established that both the limit of the integral \ref{def_contourintegral} and of the series \ref{def_quantummonodromynumber} are convergent. Then the limit is equal, so these are different expressions for the same analytic function  $\Fpm((m_i,m_{ij})_{ij})$.
\end{lemma}
\begin{proof}
	We first expand for $r_1>\ldots>r_n>0$ the integral $\Fpm^{r_1,\ldots,r_n}((m_i,m_{ij})_{ij})$ into an absolutely convergent series, by expanding 
	$(1\pm z_j/z_i)^{m_{ij}}=\sum_{k_{ij}\in \N_0} (\pm 1)^{k_{ij}}{m_{ij} \choose k_{ij}}$ and integrating it summand-wise:
	\begin{align*}
	\hspace{-.2cm}\sum_{(k_{ij})_{i<j}}
	\prod_i r_i^{m_i+\sum_{i<j} (m_{ij}-k_{ij})+\sum_{j<i}k_{ji}}
	\frac{(e^{2\pi\i(m_i+\sum_{i<j} m_{ij})}-1)/2\pi\i}{1+\sum_{i<j}(m_{ij}-k_{ij})+\sum_{j<i}k_{ji}} 
	\prod_{i<j} (\pm 1)^{k_{ij}}{m_{ij}\choose k_{ij}}
	\end{align*}
	Now by Abel's theorem, the limit of $\Fpm^{r_1,\ldots,r_n}((m_i,m_{ij})_{ij})$ to the point  $r_1=\cdots=r_n=1$ on the boundary of convergence is equal to the evaluation of this series at this point, which is the initial series in Definition \ref{def_quantummonodromynumber}, provided this series converges. Note that again it is sufficient to restrict ourselves to the dense subset of parameters $m_i+\sum_{i<j}m_{ij}\not\in \Z$, where we have a comfortable presentation for the series in Definition \ref{def_quantummonodromynumber}.\end{proof}

\begin{problem}
	Find a full analytic continuation for $\rFp,\Fp$, possibly by combining all $\Box_{r_{\sigma(1)},\ldots,r_{\sigma(n)}}$ to a closed manifold, similar to the Pochhammer contour for $n=2$.
\end{problem}

\subsection{The Main Theorem}\label{subsec_Main}

Recall the contour integral in Definition \ref{def_contourintegral}:
$$\Fp((m_i,m_{ij})_{ij})=\frac{1}{(2\pi\i)^n}\int_{\Box}
\prod_i z_i^{m_i}\prod_{i<j} (z_i- z_j)^{m_{ij}}\;\d z_1\cdots \d z_n$$

In contrast to $\Fm$, this expression is \emph{not} symmetric in the variables $z_i$, not even up to factors $e^{\pi\i m_{ij}}$, because $\Box$ contains singularities with branch points and hence involves choices. However, we now prove a substantially more subtle and interesting algebraic symmetry property for quantum monodromy numbers $\Fp$: Define the diagonal braiding 
$$q_{ij}=q_{ji}=e^{\pi\i m_{ij}},\text{ for }0\leq i<j\leq n$$
and recall from Definition \ref{def_QuantumSymmetrizer} the quantum symmetrizer $\sym_n$ and from Lemma \ref{lm_braidingfactor} the explicit formula for $\sym_{n,q}$ for diagonal braidings in terms of the braiding factor
$$q(\sigma):=\prod_{i<j,\;\sigma(j)>\sigma(i)} q_{ij}=e^{\pi\i \sum_{i<j,\;\sigma(j)>\sigma(i)} m_{ij}},\qquad \forall\sigma\in\S_n$$

\begin{maintheorem}\label{thm_QuantumSymmetrizer}
For complex parameters $(m_i,m_{ij})_{ij}$, $0\leq i<j\leq n$, which are subpolar according to Definition~\ref{def_smallnessF}, so that the contour integral $\rFp$ in Definition \ref{def_contourintegral} 
is well-defined, we have the following analytic quantum symmetrizer formula 
		$$
		\Fp((m_i,m_{ij})_{ij})=\sum_{\sigma\in \S_n} q(\sigma)\;\rFp((m_{\sigma^{-1}(i)},m_{\sigma^{-1}(i)\sigma^{-1}(j)})_{ij})=\sym_{n,q}\;\rFp((m_i,m_{ij})_{ij})
		$$
\end{maintheorem}
\begin{proof}
	In Lemma \ref{lm_contourequalsseries} we have defined $\Fp$ as the limit for $r_1=\cdots=r_n=1$ of a contour integral $\Fp^{r_1,\ldots,r_n}$ and proven that this is an analytic continuation of the series in Definition \ref{def_quantummonodromynumber}. Similarly we can write the contour integral $\rFp$ as a limit of $\rFp^{r_1,\ldots,r_n}$. So it is sufficient to prove the following claim for $r_1>r_2>\cdots >r_n>0$.
	$$\Fp^{r_1,\ldots,r_n}((m_i,m_{ij})_{ij})=\sum_{\sigma\in \S_n} q(\sigma)\;\rFp^{r_{\sigma^{-1}(1)},\ldots,r_{\sigma^{-1}(n)}}((m_{\sigma^{-1}(i)},m_{\sigma^{-1}(i)\sigma^{-1}(j)})_{ij})$$
	for the contour integrals
	\begin{align*}
	\Fp^{r_1,\ldots,r_n}((m_i,m_{ij})_{ij})&=\frac{1}{(2\pi\i)^n}\int_{\Box_{r_1,\ldots,r_n}}
	\prod_i z_i^{m_i}\prod_{i<j} (z_i- z_j)^{m_{ij}}\;\d z_1\cdots \d z_n\\
	\rFp^{r_1,\ldots,r_n}((m_i,m_{ij})_{ij})&=\frac{1}{(2\pi\i)^n}\int_{\triangle_{r_1,\ldots,r_n}}
	\prod_i z_i^{m_i}\prod_{i<j} (z_i-z_j)^{m_{ij}}\;\d z_1\cdots \d z_n
	\end{align*}
	over the integration domains 
	\begin{align*}
	\Box_{r_1,\ldots r_n}&:=\{(r_1e^{2\pi\i t_1},\ldots,r_ne^{2\pi\i t_n}) \mid 0<t_1,\cdots,t_n<1\}\\
	\triangle_{r_1,\ldots r_n}&:=\{(r_1e^{2\pi\i t_1},\ldots,r_ne^{2\pi\i t_n}) \mid 0<t_1<\cdots <t_n<1\}
	\end{align*}
	The main observation for our proof is that we can decompose the rectangular integration domain $\Box$ into $n!$ simplicial integration domains $\triangle$, up to the zero-set consisting of hyperplanes $t_i=t_j$ (the so-called braid arrangement):
	\begin{align*}
	\Box_{r_1,\ldots r_n} 
	&=\bigcup_{\sigma\in\S_n} \sigma\triangle_{r_1,\ldots r_n}\\
	&=\bigcup_{\sigma\in\S_n} \{(r_{\sigma^{-1}(1)}e^{2\pi\i t_{1}},\ldots,r_{\sigma^{-1}(n)}e^{2\pi\i t_{n}}) \mid 0<t_{\sigma(1)}<\cdots <t_{\sigma(n)}<1\}
	\end{align*} 
	Correspondingly, we have a decomposition of the contour integral 
	$$\Fp^{r_1,\ldots,r_n}((m_i,m_{ij})_{ij})
	=\sum_{\sigma\in\S_n}\frac{1}{(2\pi\i)^n}\int_{\sigma\triangle^{r_1,\ldots,r_n}}
	\prod_i z_i^{m_i}\prod_{i<j} (z_i-z_j)^{m_{ij}}\;\d z_1\cdots \d z_n$$
	We finally rewrite the integrals over $\sigma\triangle^{r_1,\ldots,r_n}$ in terms of integrals over $\triangle^{r_{\sigma^{-1}(1)},\ldots,r_{\sigma^{-1}(n)}}$ by permuting the integration variables $z_1,\ldots,z_n$, which causes, due to the choices of branches, the asserted braiding factor and concludes the proof
	\begin{align*} 
	&=\sum_{\sigma\in\S_n}e^{\pi\i \sum_{i<j,\;\sigma(j)>\sigma(i)} m_{ij}}\frac{1}{(2\pi\i)^n}\int_{\triangle^{r_{\sigma^{-1}(1)},\ldots,r_{\sigma^{-1}(n)}}}
	\prod_i z_i^{m_i}\prod_{i<j} (z_i-z_j)^{m_{ij}}\;\d z_1\cdots \d z_n\\
	&=\sum_{\sigma\in\S_n} q(\sigma)\; \rFp^{r_{\sigma^{-1}(1)},\ldots,r_{\sigma^{-1}(n)}}((m_i,m_{ij})_{ij})
	\end{align*}
	\enlargethispage{-1cm}
	An explicit way to see this braiding factor is to write the integral in polar coordinates:
	\begin{align*}
	&\Fp^{r_1,\ldots,r_n}((m_i,m_{ij})_{ij})\\
	&=\frac{1}{(2\pi\i)^n}\int_{[0,2\pi]^n}
	\prod_i \i e^{\i t_i}\d t_i \;
	\prod_i e^{\i t_i(m_i+\sum_{i<j} m_{ij})}r_i^{m_i+\sum_{i<j} m_{ij}}
	\prod_{i<j} e^{\i \theta_{ij}m_{ij}} r_{ij}^{m_{ij}}\\
	&=\frac{\prod_i r_i^{m_i+\sum_{i<j} m_{ij}}}{(2\pi)^n}\int_{[0,2\pi]^n}
	 \d t_1\cdots \d t_n \;
	e^{\i (\sum_i t_i(1+m_i) +\sum_{i<j} (t_i+\theta_{ij}) m_{ij})} 
	\prod_{i<j} r_{ij}^{m_{ij}}
	\end{align*}
	with the obvious geometric functions $r_{ij}(t_i,t_j)\in \R^+$ and $\theta_{ij}(t_i,t_j)\in [-\frac{\pi}{2},\frac{\pi}{2}]$
	\begin{align*}
	r_{ij}(t_i,t_j)&:=\sqrt{1+\left(\frac{r_j}{r_i}\right)^2 -2\frac{r_j}{r_i}\cos(t_j-t_i)}\\
	\theta_{ij}(t_i,t_j)
	&:=\tan^{-1}\left(\frac{1\sin(0)-\frac{r_j}{r_i}\sin(t_j-t_i))}{
		1\cos(0)-\frac{r_j}{r_i}\cos(t_j-t_i))}\right)
	=-\tan^{-1}\left(\frac{\sin(t_j-t_i))}{\frac{r_i}{r_j}-\cos(t_j-t_i))}\right)
\intertext{that give the polar coordinates of $z_i-z_j$.	In the limit $r_1=\cdots=r_n=1$ we have}
	r_{ij}(t_i,t_j)
	&=\sqrt{2-2\cos(t_j-t_i)}\\
	\theta_{ij}(t_i,t_j)
	&=-\tan^{-1}\left(1-\frac{\sin(t_j-t_i)}{\cos(t_j-t_i)}\right)
	=\begin{cases}
	\frac{t_j-t_i}{2}-\frac{\pi}{2},&\quad\text{for }t_j-t_i >0\\
	\frac{t_j-t_i}{2}+\frac{\pi}{2},&\quad\text{for }t_j-t_i <0\\
	\end{cases}\\
	\Fp((m_i,m_{ij})_{ij})
	&=
	\frac{1}{(2\pi)^n}
	\int_{[0,2\pi]^n} 
	\prod_i \d t_i \; 
	e^{\i(\sum_i t_i(1+m_i) + \sum_{i<j} \frac{t_i+t_j}{2} m_{ij})} \\
	&\cdot \prod_{i<j} \sqrt{2-2\cos(t_j-t_i)}^{\;m_{ij}}\prod_{i<j} e^{-\i\frac{\pi}{2} \sgn(t_j-t_i) m_{ij}}
	\end{align*}
	where we introduced the sign-function $\sgn(t_2-t_1)=\pm 1$. Only this last factor depends on the order of the $t_i$, and indeed the ratio of this factor for $0<t_{\sigma(1)}<\cdots <t_{\sigma(n)}<1$ in $\sigma\triangle$ and the factor for the standard ordering $0<t_{1}<\cdots <t_{n}<1$ in $\triangle$ is $q(\sigma)$. 
%
%
\end{proof}	
	
The main consequence of the quantum symmetrizer formula is that $\Fp((m_i,m_{ij})_{ij})$ has zeroes according to elements in the kernel of the quantum symmetrizer.	
		
\begin{example}\label{exm_F2} For $n=2$ the quantum symmetrizer formula reads:
\begin{align*}
  \Fp(m_1,m_2,m_{12})
  &=\rFp(m_1,m_2,m_{12})+e^{\pi\i m_{12}}\cdot \rFp(m_2,m_1,m_{12})
  \end{align*}
  Explicitly, this links the formulae in Example \ref{exm_F2series} for $\Fp$ and \ref{exm_rF2} for $\rFp$ as follows:
  \enlargethispage{0.5cm}
  \begin{align*}
  &\frac{e^{2\pi\i m_2}-1}{2\pi\i}\frac{e^{2\pi\i m_1+2\pi\i m_{12}}-1}{2\pi\i}
  \;\frac{1}{m_1+m_2+m_{12}+2}\cdot\\
  &\cdot \left(\Beta(m_2+1,m_{12}+1)+\frac{\sin\pi m_1}{\sin\pi(m_1+m_{12})}\Beta(m_1+1,m_{12}+1)\right)\\
  &=\frac{1}{(2\pi\i)^2}\left(1-e^{2\pi\i m_2}\right)
  \frac{\Beta(m_2+1,m_{12}+1)}{m_1+m_2+m_{12}+2}\\
  &-\frac{1}{(2\pi\i)^2}e^{2\pi\i m_2+\pi\i m_{12}}\left(1-e^{2\pi\i m_1}\right)
  \frac{\Beta(m_1+1,m_{12}+1)}{m_1+m_2+m_{12}+2}\\
  &+e^{\pi\i m_{12}}\cdot \frac{1}{(2\pi\i)^2}\left(1-e^{2\pi\i m_1}\right)
  \frac{\Beta(m_1+1,m_{12}+1)}{m_1+m_2+m_{12}+2}\\
  &-e^{\pi\i m_{12}}\cdot \frac{1}{(2\pi\i)^2}e^{2\pi\i m_1+\pi\i m_{12}}\left(1-e^{2\pi\i m_2}\right)
  \frac{\Beta(m_2+1,m_{12}+1)}{m_1+m_2+m_{12}+2}
\end{align*}

\enlargethispage{0.8cm}
In this case, two classes of zeroes are visible directly. These zeroes become very transparent from the quantum symmetrizer formula, which holds for subpolar $m_{12}>-1$:
\begin{itemize}
	\item For $m_{12}$ an odd (and positive) integer and $m_1=m_2$ we have $\Fp(m_1,m_2,m_{12})=0$, because the sine-fraction is $-1$ and hence the Beta functions cancel. But if $m_{12}$ is a negative integer, the zero is counteracted by the pole in the Gamma function. 
	
	This zero translates to a Nichols algebra relation $\zem_{\alpha}^2=0$ if $(\alpha,\alpha)$ is an odd (and positive) integer. The exceptionally non-zero value for non-subpolar $(m_{ij})_{ij}$ encodes an extension of the Nichols algebra, see Subsection \ref{subsec_EqualLiouville}.\\
	
	\item For  $m_{12}$ an even or odd (and positive) integer and $m_1,m_2$ equal modulo $\Z,2\Z$, we have $\Fp(m_1,m_2,m_{12})-(-1)^{m_{12}} \Fp(m_2,m_1,m_{12})=0$. But if $m_{12}$ an negative integer, the zero is again counteracted by the pole in the Gamma function.
	
	This zero translates to a Nichols algebra relation ${\zem_{\alpha_1}\zem_{\alpha_2}-(-1)^{m_{12}}\zem_{\alpha_2}\zem_{\alpha_1}=0}$ if $(\alpha,\alpha)$ is an even or odd (and positive) integer. The exceptionally non-zero value for non-subpolar $(m_{ij})_{ij}$ encodes an extension of the Nichols algebra, for example by a Lie algebra with brackets given by the local commutator formula in terms of the poles. For example, this appears for the affine Lie algebra in Subsection~\ref{subsec_TrivialLevel}. Note that for local operators the Nichols algebra itself is just a polynomial ring. \\
	
	\item Another zero depends on $m_i$ and will appear in Subsection \ref{subsec_EqualWeyl}: If $m_1=m_2$ and $m_1+m_2+m_{12}\in 2\N_0$, then $\Fp(m_1,m_2,m_{12})=0$. But if $m_1+m_2+m_{12}=-2$, the zero is counteracted by the pole of the denominator, even if $m_{ij}$ is subpolar. 
	
	This zero does not come from the quantum symmetrizer, but from the formula in Lemma \ref{lm_rFSelberg} for $\rFp$ in terms of $\Sel$. It does not translate into a Nichols algebra relation, but into a statement that on a given module (corresponding to $m_i$) certain elements in  the Nichols algebra act trivially. It precisely corresponds to a Weyl reflection and the exceptional non-zero value encodes an extension of an irreducible Nichols algebra module. This is discussed in Subsection \ref{subsec_EqualWeyl}.
\end{itemize}
\end{example}

More technically, our main theorem and the asymptotic bound for Selberg integrals in Lemma \ref{lm_SelbergConvergence} give the same asymptotic bound for $\Fp((m_i+k_{i},m_{ij})_{ij})$ for $k_i\to \infty$. This is sufficient to prove conditional convergence of the series in Definition \ref{def_quantummonodromynumber} for $\Fp$:

\begin{lemma}\label{lm_ConditionalConvergence}
	The series $\Fp$ converges conditionally, if for all subsets $J\subset I$  with $|J|\geq 2$ holds the following condition, which is slightly weaker than being subpolar:
	$$\sum_{i<j,\;i,j\in J} m_{ij}>-|J|$$ 
\end{lemma}
\begin{proof}
	We inductively express the ${n\choose 2}$-fold series for $\Fp$ in $n+1$ variables in terms of a single sum and  $\Fp$ in $n$ variables. Then we use the bounds from Lemma \ref{lm_SelbergConvergence}:
	\begin{align*}
	&\Fp((m_i,m_{ij})_{ij})\\
	&=\sum_k\frac{1}{1+m_0+\sum_{0<j}m_{0j}-k}\sum_{\sum_i k_{0i}=k} \Fp((m_i+k_{0i},m_{ij})_{1\leq i<j\leq n})\\
	&\sim \sum_k\frac{1}{1+m_0+\sum_{0<j}m_{0j}-k}\int_{\sum_i k_{0i}=k} \Fp((m_i+k_{0i},m_{ij})_{1\leq i<j\leq n})\\
	&\sim \sum_k\frac{1}{1+m_0+\sum_{0<j}m_{0j}-k} k^{(n-1)-n-\sum_{i<j} m_{ij}}
	\end{align*}
	This series converges whenever $\sum_{i<j} m_{ij}>-1$ which inductively shows the claim.
\end{proof}

\section{Example: The Case of Equal \texorpdfstring{$m_i,m_{ij}$}{m\_i,m\_ij}}\label{sec_Equal}

In this section we study as example the case, where all $m_i$ are equal and all $m_{ij}$ are equal, and we denote in this section 
$$m_{\alpha\alpha}=m_{ij},\hspace{.9cm} m_{\alpha\lambda}=m_i$$
$$\hspace{0.8cm}q_{\alpha\alpha}=e^{\pi\i m_{\alpha\alpha}},\hspace{0.6cm} q_{\alpha\lambda}^2=e^{2\pi\i m_{\alpha\lambda}}$$ 
We always require $q_{\alpha\alpha}\neq 1$ and thus $m_{\alpha\alpha}\not\in2\Z$. In some formulae below we will only require that all $m_i$ are equal modulo $\Z$ and all $m_{ij}$ are equal modulo $2\Z$. In this case we still have equal values for self-braiding $q_{\alpha\alpha}$ and double-braiding $q_{\alpha\lambda}^2=q_{\alpha\lambda}q_{\lambda\alpha}$. We remark that the Nichols algebra essentially depends only on these two quantities.

A typical example where this scenario appears is in the powers of a single screening $(\zem_\alpha)^n$ acting on a module $\V_\lambda$, which then involves quantum monodromy numbers $\Fp((m_{\alpha\lambda}+k_i,m_{\alpha\alpha})_{ij})$ for all $(k_i)_i\in\N_0^n$. To see the Nichols algebra relation, we denote by $\Fp((m_i,m_{ij})_{ij})^{\S_n}$ the ordinary symmetrization of $\Fp((m_i,m_{ij})_{ij})$ under permutation of the index set $I$, as in the proof of Theorem \ref{thm_NicholsAlgebra}. In the following example we see the Nichols algebra relation in rank 1 (Example \ref{exm_rank1}), or how it fails if $(m_{ij})_{ij}$ is not subpolar. 

\subsection{Calculating the Combinatorial Prefactors}\label{subsec_EqualPrefactors}
In this subsection we only require that all $m_i$ are equal modulo $\Z$ and all $m_{ij}$ are equal modulo $2\Z$.
We first review how the quantum symmetrizer formula relates $\Fp$ and $\rFp$ in this case:
 
\begin{fact}\label{fact_refcombinatoric} Assume that all $m_{ij}$ are equal modulo $2\Z$ and denote $q_{\alpha\alpha}=e^{\pi\i m_{ij}}$. Then
	\begin{enumerate}[a)]
		\item As polynomials in a formal variable $q_{\alpha\alpha}$ we have the well-known combinatorial identity
		$$\sum_{\eta\in\S_n} q_{\alpha\alpha}^{\mathrm{length}(\eta)}
		\;=\;\prod_{s=1}^{n} \frac{q_{\alpha\alpha}^s-1}{q_{\alpha\alpha}-1}=:[n]_{q_{\alpha\alpha}}!$$
		\item Then the quantum symmetrizer formula in Theorem \ref{thm_QuantumSymmetrizer} reads, for $(m_{ij})_{ij}$ subpolar, 
		$$\Fp((m_i,m_{ij})_{ij})^{\S_n}
		=[n]_{q_{\alpha\alpha}}! \cdot\rFp((m_i,m_{ij})_{ij})^{\S_n}$$
		\item This implies that for all $n\geq\ord(q_{\alpha\alpha})$ we have, again for $(m_{ij})_{ij}$ subpolar, 
		$$\Fp((m_i,m_{ij})_{ij})^{\S_n}=0$$ 
	\end{enumerate}
\end{fact}
\begin{example}
	In particular our screening operator  $\zem_\alpha$ in a lattice intertwining algebra $\V_\Lambda$ fulfills the Nichols algebra relation $\zem_{\alpha}^n=0$ for $n=\ord(e^{\pi\i(\alpha,\alpha)})$, under the additional assumption $(\alpha,\alpha)>-2/n$, so that the constant $n\times n$ matrix $m_{ij}=(\alpha,\alpha)$ is subpolar.
\end{example}

We now focus on the relationship between $\rFp$ and $\Sel$. This  depends on $(m_i)_i$ and~$q_{\alpha\lambda}^2$, which means we study the action of a power of a  single screening on specific module $\V_\lambda$:

\begin{lemma}\label{lm_refcombinatoric} Assume that all $m_i$ are equal modulo $\Z$ and all $m_{ij}$ are equal modulo $2\Z$ and denote $q_{\alpha\lambda}^2=e^{2\pi\i m_{i}}$ and $q_{\alpha\alpha}=e^{\pi\i m_{ij}}$.
	\begin{enumerate}[a)]
		\item As polynomials in two formal variables $q_{\alpha\alpha},q_{\alpha\lambda}^2$ we have the combinatorial identity:
		$$
		\sum_{k=0}^n (-1)^{n-k} q_{\alpha\lambda}^{2k} 
		\sum_{\eta\in \S_{n-k,\overline{k}}} q_{\alpha\alpha}^{\mathrm{length}(\eta)}
		\;=\; \prod_{s=0}^{n-1} \left(q_{\alpha\alpha}^s q_{\alpha\lambda}^2-1\right)$$
		where we again define the following slight variation to $(n-k,k)$-\emph{shuffles} 
		$$\S_{n-k,\overline{k}}:=\{\eta\in \S_n \mid \forall_{i<j\leq n-k}\;\eta(i)<\eta(j)\text{ and } \forall_{n-k<i<j}\;\eta(i)>\eta(j)\}$$
		\item Then Lemma \ref{lm_rFSelberg} reads, under the assumptions in Lemma \ref{lm_SelbergConvergence}. 
		$$\rFp((m_i,m_{ij})_{ij})^{\S_n}
		=\frac{1}{(2\pi\i)^n}\left(\prod_{s=0}^{n-1} \left(q_{\alpha\alpha}^s q_{\alpha\lambda}^2-1\right)\right)\cdot \Sel(m_i,m_{ij})^{\S_n}$$
		\item This implies that if $m_{\alpha\lambda}+s m_{\alpha\alpha}/2\in\Z$ for some $0<s<n$, then  $$\Fp((m_i,m_{ij})_{ij})^{\S_n}=\rFp((m_i,m_{ij})_{ij})^{\S_n}=0$$
	\end{enumerate}
\end{lemma}
\begin{proof}
	Claim (b) and (c) follows from (a) and Lemma \ref{lm_rFSelberg}.\\
	
	Claim (a) can be proven by induction. We picture the $(n+1)$st point added at the beginning and decompose $\S_{n+1-k,\overline{k}}$ for $k\neq n+1$ as a set according to the number $r$ of the $k$ points which are shuffled in front of this first point $\eta(i)<\eta(1)$. The shuffle is then uniquely determined by the remaining $(n-k,k-r)$-shuffle:
	$$\S_{n+1-k,\overline{k}} = \bigcup_{r=0}^{k} \S_{n-k,\overline{k-r}},\qquad k\neq n+1$$
	\noindent
	Accordingly we have for $k\neq n+1$
	\begin{align*}
	\sum_{\eta\in \S_{n+1-k,\overline{k}}} q_{\alpha\alpha}^{\mathrm{length}(\eta)}
	&=\sum_{r=0}^{k}\; \sum_{\eta\in \S_{n-k,\overline{k-r}}} q_{\alpha\alpha}^{\mathrm{length}(\eta)+(n+(n-1)+\cdots+(n-r+1))}
	\end{align*}
	Hence we can calculate for $n+1$, using the induction assertion for all $n-r$:
	\begin{align*}
	&\sum_{k=0}^{n+1} (-1)^{n+1-k} q_{\alpha\lambda}^{2k} 
	\sum_{\eta\in \S_{n+1-k,\overline{k}}} q_{\alpha\alpha}^{\mathrm{length}(\eta)}\\
	&= q_{\alpha\alpha}^{{ n+1\choose 2}}q_{\alpha\lambda}^{2(n+1)}
	+ \sum_{k=0}^{n} (-1)^{n+1-k} q_{\alpha\lambda}^{2k} 
	\sum_{r=0}^{k}q_{\alpha\alpha}^{{ n+1\choose 2}-{ n-r+1 \choose 2}} \hspace{-.4cm} \sum_{\eta\in \S_{n-k,\overline{k-r}}} q_{\alpha\alpha}^{\mathrm{length}(\eta)}\\
	&= q_{\alpha\alpha}^{{ n+1\choose 2}}q_{\alpha\lambda}^{2(n+1)}
	- \sum_{r=0}^{n}   q_{\alpha\alpha}^{{ n+1\choose 2}-{ n-r+1 \choose 2}} q_{\alpha\lambda}^{2r} 
	\sum_{k=r}^{n}(-1)^{\substack{(n-r)-(k-r)}}(q_{\alpha\lambda}^2)^{k-r}\hspace{-.4cm}\sum_{\eta\in \S_{n-k,\overline{k-r}}} q_{\alpha\alpha}^{\mathrm{length}(\eta)}\\
	&= q_{\alpha\alpha}^{{ n+1\choose 2}}q_{\alpha\lambda}^{2(n+1)}
	- \sum_{r=0}^{n}  q_{\alpha\alpha}^{{ n+1\choose 2}-{ n-r+1 \choose 2}}  q_{\alpha\lambda}^{2r} 
	\prod_{s=0}^{n-r-1} \left(q_{\alpha\alpha}^s q_{\alpha\lambda}^2-1\right)\\  
	\intertext{This polynomial can now be factorized as asserted, again inductively: The first asserted factor $\left(q_{\alpha\lambda}^2-1\right)$ for $s=0$ is contained in all summands except $r=n$, and this summand together with the extra summand at the beginning is also divisible by $\left(q_{\alpha\lambda}^2-1\right)$}
	&=\left(q_{\alpha\lambda}^2-1\right)\cdot
	\left(q_{\alpha\alpha}^{{ n+1\choose 2}-{1\choose 2}}q_{\alpha\lambda}^{2n}
	-\sum_{r=0}^{n-1}  q_{\alpha\alpha}^{{ n+1\choose 2}-{ n-r+1 \choose 2}}  q_{\alpha\lambda}^{2r}
	\prod_{s=1}^{n-r-1} \left(q_{\alpha\alpha}^s q_{\alpha\lambda}^2-1\right)\right)
	\\
\intertext{and continuing this process inductively gives}
	&=\left(q_{\alpha\lambda}^2-1\right)\left(q_{\alpha\alpha}q_{\alpha\lambda}^2-1\right)\cdots \left(q_{\alpha\alpha}^{n-1}q_{\alpha\lambda}^2-1\right)  \qedhere
	\end{align*}
\end{proof}

%
\begin{example}
		In particular the screening operator $\zem_\alpha$ in a lattice intertwining algebra acts on a Heisenberg vertex algebra module $\V_\lambda$ with $(\alpha,\lambda)+(n-1)(\alpha,\alpha)/2\in \Z$ such that $\zem_\alpha^{n}=0$, under the assumptions in Lemma \ref{lm_SelbergConvergence}. The only nonzero contributions to $\zem_\alpha^{n}$ can come from parameters $(m_i,m_{ij})_{ij}$ that correspond to poles of $\Sel(m_i,m_{ij})$.
\end{example}

\subsection{Calculating \texorpdfstring{$\Fp,\rFp,\Sel$}{F,tildeF,Sel}}\label{subsec_EqualCalc}
\enlargethispage{1cm}
We now return to the case where all $m_i$ are equal to some $m_{\alpha\lambda}$ and all $m_{ij}$ are equal to some $m_{\alpha\alpha}$ for $1\leq i<j\leq n$. We have  $(m_{ij})_{ij}$ subpolar iff
\begin{align*}
{n\choose 2}m_{\alpha\alpha} 
&> -n+1\\
m_{\alpha\alpha}
&>-2/n
\end{align*}
We continue to assume $q_{\alpha\alpha}\neq 1$ and thus  $m_{\alpha\alpha}\not\in 2\Z$.
We use the Selberg integral formula in Example \ref{exm_Selberg}, which we slightly rearrange and simplify using $\Gamma(z+1)=z\Gamma(z)$, and then use the combinatorial prefactors from the previous subsection to get $\rFp,\Fp$.
\begin{align*}
\Sel(m_{\alpha\lambda},0,m_{\alpha\alpha})
&=\frac{1}{n!}\prod_{j=0}^{n-1} \frac{\Gamma(m_{\alpha\lambda}+1+jm_{\alpha\alpha}/2)\Gamma(1+jm_{\alpha\alpha}/2)
\Gamma(1+(j+1)m_{\alpha\alpha}/2)}{\Gamma(m_{\alpha\lambda}+2+(n+j-1)m_{\alpha\alpha}/2)\Gamma(1+m_{\alpha\alpha}/2)}\\
&=\frac{1}{n!\;\Gamma(1+m_{\alpha\alpha}/2)^n}\cdot \frac{1}{m_{\alpha\lambda}+1+(n-1)m_{\alpha\alpha}/2}\cdot\Gamma(1+nm_{\alpha\alpha}/2)\\
&\cdot \prod_{j=1}^{n-1} \frac{\Gamma(m_{\alpha\lambda}+1+(j-1)m_{\alpha\alpha}/2)
	\Gamma(1+jm_{\alpha\alpha}/2)^2}{\Gamma(m_{\alpha\lambda}+2+(n+j-1)m_{\alpha\alpha}/2)}\\
\rFp(m_{\alpha\lambda},m_{\alpha\alpha})
&=\frac{1}{(2\pi\i)^n}\frac{1}{n!\;\Gamma(1+m_{\alpha\alpha}/2)^n}\cdot \frac{q_{\alpha\alpha}^{n-1} q_{\alpha\lambda}^2-1}{m_{\alpha\lambda}+1+(n-1)m_{\alpha\alpha}/2}\cdot\Gamma(1+nm_{\alpha\alpha}/2)\\
&\cdot \prod_{j=1}^{n-1} \left(q_{\alpha\alpha}^{j-1} q_{\alpha\lambda}^2-1\right)\frac{\Gamma(m_{\alpha\lambda}+1+(j-1)m_{\alpha\alpha}/2)
	\Gamma(1+jm_{\alpha\alpha}/2)^2}{\Gamma(m_{\alpha\lambda}+2+(n+j-1)m_{\alpha\alpha}/2)}\\
\Fp(m_{\alpha\lambda},m_{\alpha\alpha})
&=\frac{1}{(2\pi\i)^n}\frac{1}{n!\;\Gamma(1+m_{\alpha\alpha}/2)^n}\cdot \frac{q_{\alpha\alpha}^{n-1} q_{\alpha\lambda}^2-1}{m_{\alpha\lambda}+1+(n-1)m_{\alpha\alpha}/2}\cdot[n]_{q_{\alpha\alpha}}\Gamma(1+nm_{\alpha\alpha}/2)\\
&\cdot \prod_{j=1}^{n-1} [j]_{q_{\alpha\alpha}}\left(q_{\alpha\alpha}^{j-1} q_{\alpha\lambda}^2-1\right)\frac{\Gamma(m_{\alpha\lambda}+1+(j-1)m_{\alpha\alpha}/2)
	\Gamma(1+jm_{\alpha\alpha}/2)^2}{\Gamma(m_{\alpha\lambda}+2+(n+j-1)m_{\alpha\alpha}/2)}
\end{align*}
The terms in the product already appear for smaller $n$, so we concentrate on the remaining terms. These terms show two types of poles. 
\begin{itemize}
	\item $m_{\alpha\lambda}+(n-1)m_{\alpha\alpha}/2=-1$. This pole only appears in $\Sel$, accordingly the second assumption in Lemma \ref{lm_SelbergConvergence} is violated. In $\Fp,\rFp$ this pole is suppressed by the factor $q_{\alpha\alpha}^{n-1} q_{\alpha\lambda}^2-1$, but it leads to exceptionally non-zero values of $\Fp,\rFp$ in violation of Lemma \ref{lm_refcombinatoric} (c). We will discuss in Subsection~\ref{subsec_EqualWeyl} how this corresponds to module extensions and reflection operators.
	\item $nm_{\alpha\alpha}/2\in-\N$ for $j=n-1$. This pole appears in $\Sel,\rFp$, accordingly the assumption in Lemmas \ref{lm_SelbergConvergence} and \ref{lm_rFpConvergence} of $(m_{ij})_{ij}$ being subpolar is violated. However, in $\Fp$ this pole is suppressed by the factor $q_{\alpha\alpha}^{n}-1$, but it leads to exceptionally non-zero values of $\Fp$ in violation of Fact \ref{fact_refcombinatoric}, and hence a violation of the Nichols algebra relation $(\zem_\alpha)^n=0$, which we will discuss in Subsection \ref{subsec_EqualLiouville}. 
\end{itemize}
\enlargethispage{-2cm}
If $m_i,m_{ij}$ are only equal modulo $\Z,2\Z$, then the combinatorial prefactors still apply, and the Selberg integral can be evaluated via Jack polynomials \cite{Kad97,Mac95,TW}.

\subsection{The Pole of \texorpdfstring{$\Sel$ at $m_{\alpha\lambda}+(n-1)m_{\alpha\alpha}/2=-1$}{Sel}: Reflection operators}\label{subsec_EqualWeyl}

We again assume equal $m_i=m_{\alpha\lambda},m_{ij}=m_{\alpha\alpha}$ and assume $m_{\alpha\alpha}>-2/n$ so $(m_{ij})_{ij}$ is subpolar. 

We consider in Subsection \ref{subsec_EqualCalc} the first type of pole $m_{\alpha\lambda}+(n-1)m_{\alpha\alpha}/2=-1$ in $\Sel(m_{\alpha\lambda},0,m_{\alpha\alpha})$ and calculate in this case
\begin{align*}
\Fp(m_{\alpha\lambda},m_{\alpha\alpha})
&=\frac{1}{(2\pi\i)^n}\frac{[n]_{q_{\alpha\alpha}}!}{n!\;\Gamma(1+m_{\alpha\alpha}/2)^n}\cdot 2\pi\i\cdot \Gamma(1+nm_{\alpha\alpha}/2)\\
&\cdot \prod_{j=1}^{n-1} \left(e^{2\pi\i(j-n)m_{\alpha\alpha}/2}-1\right)\frac{\Gamma((j-n)m_{\alpha\alpha}/2)
	\Gamma(1+jm_{\alpha\alpha}/2)^2}{\Gamma(1+jm_{\alpha\alpha}/2)}\\
&= e^{-\pi\i\;{n\choose 2}m_{\alpha\alpha}/2}\frac{[n]_{q_{\alpha\alpha}}! \;\Gamma(1+nm_{\alpha\alpha}/2) }{n!\;\;\Gamma(1+m_{\alpha\alpha}/2)^n}
\end{align*}
 where in the first line we inserted the condition and in the critical factor used the residue $\lim_{z\to 0} (e^{2\pi\i z}-1)\Gamma(z)=2\pi\i$, and in the second line we canceled $\Gamma(1+jm_{\alpha\alpha}/2)$, resorted the product and used Euler's reflection formula
 $\Gamma(z)\Gamma(1-z)=\frac{\pi}{\sin(\pi z)}$.
 
 In particular, this explicit formula shows that indeed $\Fp(m_{\alpha\lambda},m_{\alpha\alpha})\neq 0$ if and only if  $n<\ord(q_{\alpha\alpha})$, or if $nm_{\alpha\alpha}/2\in-\N$ but this is not in the subpolar region.\\

 We now observe that $\Sel((m_{\alpha\lambda}+k_i,0,m_{\alpha\alpha})_{ij})$ has, according to Lemma \ref{lm_SelbergConvergence} with $m_{\alpha\lambda}+(n-1)m_{\alpha\alpha}/2=-1,m_{\alpha\alpha}>-2/n$, a pole precisely if all $k_i=0$, otherwise it converges. In turn, $\Fp((m_{\alpha\lambda}+k_i,m_{\alpha\alpha})_{ij})\neq 0$ if and only if all $k_i=0$ and $n<\ord(q_{\alpha\alpha})$. 
 
 For our screenings $\zem_\alpha=\zem_{\exp{\phi_\alpha}}$ this has the consequence, that only this leading term contributes to  $(\zem_{\alpha})^n$ in Theorem \ref{thm_associativity}, for example we get explicitly:
 \begin{align*}
 (\zem_{\alpha})^n\; \exp{\phi_\lambda}
 &=\sum_{(k_i)_i}\sum_{(l_i)_i} \;\Fp((m_{\alpha\lambda}+k_i,m_{\alpha\alpha})_{ij})\cdot
 \prod_{1\leq i \leq n}\underbrace{\langle \exp{\phi_\alpha},\exp{\phi_\lambda}\rangle_{m_{\alpha\lambda}}}_{=1}
 \cdot\;\exp{\phi_\lambda}\prod_{i=n}^1 \frac{\partial^{k_i}}{k_i!}\exp{\phi_\alpha}\\
 &=\Fp(m_{\alpha\lambda},m_{\alpha\alpha})\cdot\exp{\phi_{\lambda+n\alpha}}
 \end{align*}
 This is remarkably simple for a composition of non-local operators which each produce infinitely many terms in all degrees - already this fact might suggest that $(\zem_\alpha)^n$ is local.

 \begin{definition}
	Let $n$ depending on $m_{\alpha\lambda},m_{\alpha\alpha}$ be the smallest positive integer such that 
	$$m_{\alpha\lambda}+(n-1)m_{\alpha\alpha}/2\in \Z$$
	If it exists, then we define the \emph{reflection operator} on $\V_\lambda$ by the screening operator power
	$$\zem_{s_\alpha}:=(\zem_\alpha)^n:\;\V_\lambda\to \V_{\lambda+n\alpha}$$
\end{definition}
\enlargethispage{1.3cm}
 \begin{conjecture}
	The reflection operators $\zem_{s_\alpha}$ are local operators, in the sense that they commute with vertex algebra elements in the kernel of $\zem_{s_\alpha}$ and thus give vertex algebra module homomorphisms over this kernel. In particular for a suitable choice of Virasoro action as in Section \ref{sec_KazhdanLusztig}, the reflection operators are Virasoro module homomorphisms. 
	One could try a proof along the lines of \cite{TW} Thm. 4.17. for rank $1$ and the Virasoro~action.
\end{conjecture}  
For example, in rank $1$ these particular powers appear in Felder's complex \cite{Fel89}. Further interpretation of the proposed definitions come from the Lie-theoretic setting, which will be discussed more thoroughly in Section \ref{sec_KazhdanLusztig}. 
\begin{example}
	Let $-\alpha_i/\sqrt{p}$ be rescaled simple roots of the root lattice of a finite-dimensional semisimple Lie algebra, where $p$ is any natural number that is divisible by ${d_i:=(\alpha_i,\alpha_i)/2}$, so we have $\ord(q_{ii})=\ord(e^{\pi\i (-\alpha_i/\sqrt{p},-\alpha_i/\sqrt{p})})=p/d_i$. 
	
	We rewrite the defining condition for $n$ using the Weyl vector $\rho$, fulfilling $(\alpha_i,\rho)=d_i$, and the Weyl vector of the dual rootsystem $\rho^\vee$, fulfilling $(\alpha_i,\rho^\vee)=1$, for simple roots~$\alpha_i$.
	\begin{align*}
	m_{\alpha\lambda}+(n-1)m_{\alpha\alpha}/2&=-r\in \Z\\
	\Longleftrightarrow\qquad
	(\alpha,\lambda+\rho+\rho^\vee pr)\;&\in\; [0,p)
  	\end{align*}
	For this set of weights, the reflection operator on $-\alpha_i/\sqrt{p}$ gives reflection on the affine hyperplane $n=(\alpha,\lambda+\rho+\rho^\vee pr)=0$, which is the affine hyperplane orthogonal to $\alpha_i$ and through the point $-\rho+\rho^\vee pr$. In this sense, reflection operators implement the Weyl group dot-action corresponding to the respective alcove. 
	
	For example for $\g=\sl_2$ the weights are  $k\alpha/2\sqrt{p}$ for $k\in\Z$, and the defining condition reads $k+1=n+pr$. Below we plot for $p=3$ the weights with  $k=-3,\ldots, 7$ as vertical lines for and the reflection operator $(\zem_{-\alpha/\sqrt{p}})^n$ as an arrow pointing from $k$ to $k-2n$. 
	
	\begin{center}
				\begin{tikzpicture}[
				scale=.6,
				zemlja/.style={thick,->,shorten >=2pt,shorten <=2pt,>=stealth},
				level/.style={thin,densely dashed}
				]
				\draw[zemlja] (-1cm,0cm)--(-1.4cm,0cm) ;
				\draw[zemlja] (-0cm,-1cm)--(-2cm,-1cm) ; 
				\draw[zemlja] (1cm,-2cm)--(-3cm,-2cm) ; 
				
				\draw[zemlja] (2cm,-3cm)--(1.6cm,-3cm) ;
				\draw[zemlja] (3cm,-2cm)--(1cm,-2cm) ; 
				\draw[zemlja] (4cm,-1cm)--(0cm,-1cm) ; 
				
				\draw[zemlja] (5cm,0cm)--(4.6cm,0cm) ;
				\draw[zemlja] (6cm,-1cm)--(4cm,-1cm) ; 
				\draw[zemlja] (7cm,-2cm)--(3cm,-2cm) ; 
				
				\draw[level] (-3cm,1cm)--node[above=1cm,fill=white] {$-3\;$}(-3cm,-4cm) ;  
				\draw[level] (-2cm,1cm)--node[above=1cm,fill=white] {$-2\;$}(-2cm,-4cm) ;  
				\draw[level] (-1cm,1cm)--node[above=1cm,fill=white] {$-1\;$}(-1cm,-4cm)
					node[below=0cm,fill=white]{$-\rho$};
				\draw[level] (0cm,1cm)--node[above=1cm,fill=white] {$0$}(0cm,-4cm) ;  
				\draw[level] (1cm,1cm)--node[above=1cm,fill=white] {$1$}(1cm,-4cm) ; 
				\draw[level] (2cm,1cm)--node[above=1cm,fill=white] {$2$}(2cm,-4cm)
					node[below=0cm,fill=white]{$-\rho+3\rho^\vee$} ; 
				\draw[level] (3cm,1cm)--node[above=1cm,fill=white] {$3$}(3cm,-4cm) ;  
				\draw[level] (4cm,1cm)--node[above=1cm,fill=white] {$4$}(4cm,-4cm) ;  
				\draw[level] (5cm,1cm)--node[above=1cm,fill=white] {$5$}(5cm,-4cm)
					node[below=0cm,fill=white]{$-\rho+6\rho^\vee$} ;   
				\draw[level] (6cm,1cm)--node[above=1cm,fill=white] {$6$}(6cm,-4cm) ;  
				\draw[level] (7cm,1cm)--node[above=1cm,fill=white] {$7$}(7cm,-4cm) ;  
				\end{tikzpicture}
	\end{center}
	
\end{example}
\noindent
For diagonal Nichols algebras, the reflection operators $\zem_{s_i}$ matches Kharchenko's formula
$$q_{\alpha\alpha}^{n-1}q_{\alpha\lambda}q_{\lambda\alpha}=1 \quad\text{ or }\quad q_{\alpha\alpha}^{n}=1$$ 
In that regard one could say that the vanishing terms in Lemma \ref{lm_refcombinatoric} correspond to the Nichols algebra acting on irreducible modules, while poles in the Selberg integral generate extensions of irreducible Nichols algebra modules. In the quantum group examples, the blocks of indecomposable modules appear precisely along a shifted Weyl group orbit.

\subsection{The Poles of \texorpdfstring{$\rFp$ at $m_{\alpha\alpha}=-2/n$}{tildeF}: Violated Nichols Algebra Relations}\label{subsec_EqualLiouville}

For $m_{\alpha\alpha}=-2/n$ we have a simple pole in $\rFp(m_i,m_{\alpha\alpha})$ due to $\Gamma(1+nm_{\alpha\alpha}/2)$. This pole is removed in $\Fp$ by the factor $q_{\alpha\alpha}^n-1$ in the quantum symmetrizer, but the result is exceptionally non-zero, in contrast to the zero that would be expected from the quantum symmetrizer formula, see Subsection \ref{subsec_EqualCalc}. The assumption of $(m_{ij})_{ij}$ being subpolar in our main theorem is not fulfilled. In turn we have a violation of the Nichols algebra relation $\zem_\alpha^n\neq 0$ at $q_{\alpha\alpha}^n=1$. These exceptionally non-zero operators are again remarkable:

\begin{theorem}
	For $m_{\alpha\alpha}=-2/n$ the following holds 
	\begin{enumerate}[a)]
		\item The quantum monodromy numbers $\Fp((m_i)_i,m_{\alpha\alpha})$ only depend on the overall sum $s:=\sum_i m_i$. We observe that our proof provides a similar assertion in general for violated Nichols algebra relations on the boundary of the subpolar region.
		\item Using a) we can calculate explicitly: 
		$$\Fp((m_i)_i,m_{\alpha\alpha})=\frac{C}{s+1},\qquad
		C=\frac{1}{(2\pi\i)^{n-1}}\frac{1}{q-1}\frac{[n-1]_q}{(n-1)!\;^n} 
		\prod_{j=0}^{n-1}\left(q^{j-1-s}-1\right) \frac{\Gamma(1-j/n)^2}{\Gamma(1-1/n)}$$
		\item As a consequence, the respective power of this screening is given by one local screening
		$$(\zem_\alpha)^n=C\cdot \zem_{n\alpha}$$
	\end{enumerate}
\end{theorem}
\begin{proof}
	We remark that assertion a) and b) can in this case also be proven directly using Kadell's integral, where the singular term appears only for partitions $\lambda=(1,1,\ldots,0,\ldots)$. \\

	For {assertion a)} we first observe that for all proper subsets $J\subsetneq I$ the subpolar condition holds. Only for $J=I$ we have $\sum_{i<j} m_{ij}={n\choose 2}(-2/n)=-n+1$, which means the subpolar condition fails with equality. 
	
	We now compare $\Fp((m_i)_i,m_{ij})$ with $\Fp(m_i',m_{ij})$ with evenly distributed $m_i'=s/n$: 
	\begin{align*}
	&\Fp((m_i,m_{ij})_{ij})-\Fp(s/n,(m_{ij})_{ij})
	=\int_{\Box}\mathrm{d}z_1\cdots \mathrm{d}z_n\; 
	\left(\prod_i z_i^{m_i}-\prod_i z_i^{s/n}\right)\prod_{i<j} (z_i-z_j)^{m_{ij}}\\
	&\rFp((m_i,m_{ij})_{ij})-\rFp(s/n,(m_{ij})_{ij})
	=\int_{\triangle}\mathrm{d}z_1\cdots \mathrm{d}z_n\; 
	\left(\prod_i z_i^{m_i}-\prod_i z_i^{s/n}\right)\prod_{i<j} (z_i-z_j)^{m_{ij}}
	\end{align*}
	At points $(z_1,\ldots, z_n)$ where only a subset of coordinates are equal, the function is integrable by our initial observation. On the other hand if all $z_i\to z$ the large bracket now decreases the problematic pole order, so the function is again integrable.
	In consequence, we can again use the idea of our quantum symmetrizer formula and find:
	\begin{align*}
	\Fp((m_i,m_{ij})_{ij})-\Fp(s/n,(m_{ij})_{ij})
	&=[n]_q!(\rFp((m_i,m_{ij})_{ij})-\rFp(s/n,(m_{ij})_{ij}))=0
	\end{align*}
	
	For {assertion b)} we can now simply calculate 
	$\Fp(s/n,-2/n)$ from our formulae in Subsection \ref{subsec_EqualCalc}, with $q:=q_{\alpha\alpha}=e^{-\frac{2\pi\i}{n}}$ and $q_{\alpha\lambda}^2=q^{-s}$, most Gamma functions cancel, and  we use again the residue $\frac{2\pi\i}{q-1}$ for the critical term, here $[n]_{q_{\alpha\alpha}}\Gamma(1+nm_{\alpha\alpha}/2)$:
	\enlargethispage{1cm}
	\begin{align*}
	\Fp(s/n,-2/n)&=\frac{1}{(2\pi\i)^n}\frac{1}{n!\;\Gamma(1-1/n)^n}\cdot \frac{q^{n-1-s}-1}{s/n+1/n}\cdot \frac{2\pi\i}{q-1}\\
	&\cdot \prod_{j=1}^{n-1} [j]_{q} \left(q^{j-1-s}-1\right) \frac{\Gamma(s/n+1-(j-1)n)
		\Gamma(1-j/n)^2}{\Gamma(s/n+2-(n+j-1)/n)}\\
	&=\frac{1}{s+1}\cdot \frac{1}{(2\pi\i)^{n-1}}\frac{1}{q-1}\frac{[n-1]_q}{(n-1)!\;^n} 
	 \prod_{j=0}^{n-1}\left(q^{j-1-s}-1\right) \frac{\Gamma(1-j/n)^2}{\Gamma(1-1/n)}
	\end{align*}
	
	For {assertion c)} we proceed as in the last section: We apply the expression from assertion b) to the associativity formula. This formula simplifies because $a_i=\exp{\phi_\alpha}$ is grouplike and  $\langle a,a\rangle_{m_{\alpha\alpha}}=1$ and all $a_i$ are the same, but the main conclusion holds~generally.
	
	\begin{align*}
	&\left(\prod_{i=1}^n \resY(a)\right)\; v\\
	&=\sum_{(k_i)_i}\sum_{(m_i)_{i}} \;\Fp((m_i+k_i)_i,m_{\alpha\alpha}\cdot
	\prod_{1\leq i \leq n}\langle a,v^{(n-i+1)}\rangle_{m_i}
	\cdot v^{(n+1)}\prod_{i=n}^1 \frac{\partial^{k_i}}{k_i!}a\\
	&=C\sum_{m,k}\frac{1}{m+k+1}
	\sum_{(m_i)_i,\sum_i m_i=m}	
	\;\prod_{1\leq i \leq n}\langle a,v^{(n-i+1)}\rangle_{m_i}
	\sum_{(k_i)_i,\sum_i k_i=k}v^{(n+1)}\prod_{i=n}^1 \frac{\partial^{k_i}}{k_i!}a\\
	&=C\sum_{m,k}\frac{1}{m+k+1}\langle a^n,v^{(1)}\rangle_{m}
	v^{(2)}\prod_{i=n}^1 \frac{\partial^{k}}{k!}a^n\\
	&=\resY(a^n)v
	\end{align*}  
	
\end{proof}
\enlargethispage{.5cm}

\begin{example}[Liouville and Toda case]
	Let $\Lambda=\frac{1}{\sqrt{p}}\Lambda_R$ be the rescaled root lattice with $p$ divisible by all $(\alpha,\alpha)/2$. Take now $p<0$, so $\Lambda$ is a negative definite lattice. Then our result shows that all Nichols algebra truncation relations for short screenings are violated; instead of zero they give long screenings
	$$\left(\zem_{\alpha_i/\sqrt{p}}\right)^{p/d_i}=C_i\cdot \zem_{\alpha_i^\vee\sqrt{p}}$$
	Moreover, since for $p<0$ holds $(\alpha_i\sqrt{p},\alpha_j\sqrt{p})\geq 0$ for $i\neq j$, the long screening commute. \\
	
	We would conjecture that the same is true for non-simple roots, and that the resulting non-vanishing powers give a central subalgebra isomorphic to the algebra of functions on the Borel subgroup of the corresponding dual Lie group, as in the Kac-Procesi-DeConcini quantum group. More generally, we expect an extension of any  Nichols algebra by a large center according to Cartan-type roots \cite{AA} and possibly exceptional Serre relations.\\  

	Physically, these negative-definite cases could be called \emph{non-compact cases} and the kernel of screenings for negative-definite root lattices appear in \emph{Liouville- and Toda-theory} \cite{FF93}, there for generic $q$. The representation theory of the kernel of this infinite-dimensional algebra of screenings becomes a highly infinite category; it is graded by the eigenvalues of the elements $\zem_{\alpha_i^\vee\sqrt{p}}$ in the large center. 
\end{example}	
	
	\newpage
\section{Applications to the Logarithmic Kazhdan-Lusztig Correspondence}\label{sec_KazhdanLusztig}

\subsection{Nichols Algebra Action}\label{subsec_NicholsAlgebraAction}

As a corollary of the previous two sections, we prove:

\begin{theorem}\label{thm_NicholsAlgebra}
  Let $\Lambda\subset \C^r$ be a positive-definite lattice and $\alpha_1,\ldots, \alpha_{r}$ be a fixed basis, which fulfills $\Vert\alpha_i\Vert^2\leq 1$. Then the screening operators $\zem_{\alpha_i}:\V_\Lambda\to \bar{\V}_\Lambda$ in Definition~\ref{def_screening} are composable arbitrarily often and they obey the relations of the diagonal Nichols algebra generated by $\zem_{\alpha_i}$ with braiding matrix 
  $$q_{ij}=e^{\pi\i\;(\alpha_i,\alpha_j)}$$
  \end{theorem}
  \begin{proof}
  The essential idea is now clear: Quantum monodromy numbers $\Fp$ can be written by Theorem \ref{thm_QuantumSymmetrizer} as quantum symmetrizer with diagonal braiding $q_{ij}=e^{\pi\i m_{ij}}$ of some $\rFp$. Thus every formal linear combination of quantum monodromy numbers in the kernel of the quantum symmetrizer vanishes. Since $\V_\Lambda$ is based on a commutative
  and cocommutative Hopf algebra, this will imply that a linear combination of compositions of screening operators $\zem_{\alpha_i}$ vanishes, if it is in the kernel of the quantum symmetrizer. But the Nichols algebra is by definition the quotient of the free algebra by the kernel of this quantum symmetrizer.
  We now make this idea precise: \\
  
  We wish to prove a vanishing statement about monomials $\zem_{\alpha_{i(1)}}\cdots \zem_{\alpha_{i(n)}}$ of total degree $n$ in generators $\zem_{\alpha_1},\ldots,\zem_{\alpha_r}$. Such monomials are formally parametrized by maps $i(-):X\to I$ from the index set $X=\{1,\ldots, n\}$ to $I=\{1,\ldots,r\}$. Linear combinations of monomials in total degree $n$ are parametrized by their coefficient  $c_{i(-)}\in\C$ for each monomial. For a given map $i(-)$ we extend the braiding matrix to $X$ by $q_{i(x),i(y)}=e^{\pi\i(\alpha_{i(x)},\alpha_{i(y)})}$ and denote accordingly the braiding factor $q_{i(-)}(\sigma)$, for all $\sigma\in\S_n$ as in Lemma \ref{lm_braidingfactor}.\\
  
  By Lemma \ref{lm_smallness} the condition on $\Lambda$ and $\Vert\alpha_i\Vert^2\leq 1$ implies for all $\alpha_{i(1)},\ldots, \alpha_{i(n)}$ that $m_{i(x),i(y)}$ is subpolar; note that since we assumed $\alpha_i$ to be a basis, we have $\alpha_{i(x)}\neq -\alpha_{i(y)}$, so the condition $\Vert\alpha_i\Vert^2\leq 1$ with equality suffices.  By Lemma \ref{lm_ConditionalConvergence} this implies not only the well-definedness of the contour integral expression of $\Fp$, but indeed for the series expression of $\Fp$. Hence the composition of screening operators $\zem_{\alpha_{i(1)}}\cdots \zem_{\alpha_{i(n)}}$ in Theorem \ref{thm_associativity} has convergent coefficients. This proves the first claim. \\
  
  Let $\C^I$ be the vector space with formal basis $e_i$ and diagonal braiding $q_{ij}=e^{\pi\i(\alpha_i,\alpha_j)}$, and $\B(\C^I,q)$ the Nichols algebra. Consider in $\B(\C^I,q)$ an arbitrary relation in degree $n$:
  $$0\stackrel{!}{=} e_c:=\sum_{i(-)} c_{i(-)}\cdot \left(e_{i(1)}\otimes \cdots \otimes  e_{i(n)}\right)$$ 
  By the defining property of a Nichols algebra this means that in the tensor algebra 
  \begin{align*}
  \hspace{.5cm}0\stackrel{!}{=}\sym_q(e_c)
  &:=\sum_{\sigma\in\S_n}\sum_{i(-)} q_{i(-)}(\sigma) \cdot c_{i(-)}\cdot \left(e_{i(1)}\otimes \cdots \otimes  e_{i(n)}\right)\\
  &=\sum_{i(-)} \left(\sum_{\sigma\in\S_n} q_{i(-)}(\sigma)c_{i(\sigma(-))}\right)\left(e_{i(1)}\otimes \cdots \otimes  e_{i(n)}\right)
  \end{align*}
  so the large bracket vanishes for all $i(-)$. In this situation, we wish to prove that the corresponding expression in the screening operators vanishes on any  element $v\in \V_\Lambda$.
  $$0\stackrel{?}{=}\sum_{i(-)} c_{i(-)}\left(\prod_{x=1}^n \zem_{\alpha_{i(x)}}\right)v\hspace{5.8cm}$$
  We apply our Associativity Theorem \ref{thm_associativity}. For $a_x$ grouplike and $\langle a_x,a_y\rangle=1\cdot z^{m_{xy}}$ with parameters $m_{xy}=(\alpha_{i(x)},\alpha_{i(y)})$ this theorem reads:
  \begin{align*}
  &\left(\prod_{x=1}^n \resY(a_x)\right)\; v\\
  &=\sum_{(k_x)_x \in \N_0^n}\sum_{(m_x,m_{xy})_{xy}} \;
  \prod_{1\leq x \leq n}\langle a_x^{(1)},v^{(x)}\rangle_{m_x}
  \prod_{1\leq x<y \leq n} \langle a_x^{(n-y+2)},a_y^{(n-x+1)}\rangle_{m_{xy}}\\
  &\cdot v^{(n+1)}\prod_{x=n}^1 \frac{\partial^{k_x}}{k_x!}a_x^{(n+1)}\cdot\Fp((m_x+k_x,m_{xy})_{x,y})\\
  \
  &=\sum_{(k_x)_x \in \N_0^n}\underbrace{\left( 
  \sum_{(m_x)_{x}} \; \prod_{1\leq x \leq n}\langle \exp{\phi_{\alpha_{i(x)}}},v^{(x)}\rangle_{m_{i(x)}}
  \cdot v^{(n+1)}\prod_{x=n}^1 \frac{\partial^{k_x}}{k_x!}\exp{\phi_{\alpha_{i(x)}}}
  \right)}_{=:[k_1,\ldots,k_n]}\\
  &\cdot \sum_{i(-)} c_{i(-)}\;\Fp((m_{i(x)}+k_x,m_{i(x)i(y)})_{x,y})
  \end{align*}
  We now abbreviate the large bracket by $[k_1,\ldots,k_n]\in\V_\Lambda$ for $k_1,\ldots,k_m\in\N_0$. The only property necessary for the proof now is that, due to the commutativity and cocommutativity of the Hopf algebra $\V_\Lambda$, the symbols $[k_1,\ldots,k_n]$ are permutation invariant:
  ${[k_{\sigma(1)},\ldots, k_{\sigma(n)}]=[k_1,\ldots,k_n],\sigma\in \S_n}$. 
  Then we apply the quantum symmetrizer formula to $\Fp$ and prove that the respective expression in screening operators is zero as asserted:
\enlargethispage{0.5cm}
  \begin{align*}
&\sum_{i()} c_{i()}\left(\prod_{x=1}^n \zem_{\alpha_{i(x)}}\right)\; v\\
    &=\hspace{-.3cm}\sum_{(k_x)_x\in \N_0^n}[k_1,\ldots, k_n]\sum_{i(-)} c_{i(-)}\; \Fp((m_{i(x)}+k_x,m_{i(x),i(y)})_{1 \leq x,y\leq n}) \\
    &=\hspace{-.3cm}\sum_{(k_x)_x\in \N_0^n}[k_1,\ldots, k_n]\sum_{i(-)}\sum_{\sigma\in \S_n} c_{i(-)} q_{i(-)}(\sigma)\;\rFp((m_{i(\sigma^{-1}(x))}+k_{\sigma^{-1}(x)},m_{i(\sigma^{-1}(x))i(\sigma^{-1}(y))})_{x,y})\\
 &=\hspace{-.3cm}\sum_{(k_x)_x\in \N_0^n}[k_1,\ldots, k_n]\cdot\sum_{i(-)} 
\sum_{\sigma\in \S_n} c_{i(\sigma(-))} q_{i(-)}(\sigma)\;\rFp((m_{i(x)}+k_{\sigma^{-1}(x)},m_{i(x)i(y)})_{x,y})\\
  &=\hspace{-.3cm}\sum_{(k_x)_x\in \N_0^n}[k_1,\ldots, k_n]\cdot\sum_{i(-)} 
 \cdot\left(\sum_{\sigma\in \S_n} c_{i(\sigma(-))} q_{i(-)}(\sigma)\right)\;\rFp((m_{i(x)}+k_{x},m_{i(x)i(y)})_{x,y})
 =0
    \end{align*}
  We remark that using the symmetry of $[k_1,\ldots,k_n]$ in the last equality is no formality: In the Nichols algebra e.g. the implication $q_{\alpha\alpha}=-1 \Rightarrow (\zem_\alpha)^2=0$ makes use of the fact that we have two equal factors, but in the screening operator power $(\zem_\alpha)^2$ we have a-prior different $m_x=m_{\alpha\lambda}+k_x$, but it turns out to be sufficient that the $k_x$ enter permutation-invariantly. Differently spoken, we have e.g. typically $\Fp(m_{\alpha\lambda}+k_1,m_{\alpha\lambda}+k_2,m_{\alpha\alpha})\neq 0$ for $m_{\alpha\alpha}$ a positive odd integer, but our result states that the symmetrization vanishes:
  \begin{align*}&\Fp(m_{\alpha\lambda}+k_1,m_{\alpha\lambda}+k_2,m_{\alpha\alpha})^{\S_2}\\
  &:=\frac{1}{2}\left(\vphantom{\sum}\Fp(m_{\alpha\lambda}+k_1,m_{\alpha\lambda}+k_2,m_{\alpha\alpha})
  +\Fp(m_{\alpha\lambda}+k_2,m_{\alpha\lambda}+k_1,m_{\alpha\alpha})\right)=0
  \end{align*}
  We also remark that in non-commutative cases, such as $\V_\Lambda^\epsilon$ with nontrivial $2$-cocycle $\epsilon(\alpha_i,\alpha_j)$, the symmetry of $[k_1,\ldots,k_n]$ is replaced by some braided version.. 
  As a result, we encounter a Nichols algebra with a braiding, which is the product of the monodromy braiding $e^{\pi\i(\alpha_i,\alpha_j)}$ and the additional symmetric braiding $\epsilon(\alpha_i,\alpha_j)\epsilon^{-1}(\alpha_j,\alpha_i)$. The crucial quantities for the Nichols algebra $q_{ii},q_{ij}q_{ji}$ are thereby kept invariant, and algebraically the resulting Nichols algebra is a Doi twist of the former Nichols algebra by $\epsilon(\alpha,\beta)$.
\end{proof}
\begin{remark}
  The author would conjecture that the endomorphisms $\zem_{\alpha_i}$ even generate precisely the Nichols algebra, or respectively, that the quantum monodromy numbers fulfill precisely the Nichols algebra relations. This seems to be suggested  by the universal property of the Nichols algebra as the smallest quotient with a Hopf algebra structure, since is seems that the screening operators  admit the  
  coproduct $\Delta(\zem_\alpha)=1\otimes \zem_\alpha+\zem_\alpha\otimes 1$, because it acts as derivation with respect to certain vertex algebra structures.
\end{remark}

As a main application, our result proves a long-standing expectation about the appearance of small quantum groups as screening operators:

\enlargethispage{2cm}
\begin{example}[Quantum groups]
  Let $\g$ be a finite-dimensional semisimple Lie algebra with root lattice $\Lambda_R$ and Killing form $(-,-)$. Choose an even integer $\ell=2p\in\N$, such that all $(\alpha_i,\alpha_i)$ divide $\ell$. Consider the rescaled root lattice $\Lambda:=\frac{1}{\sqrt{p}}\Lambda_R$ with basis $-\alpha_i/\sqrt{p}$, to which we associate in the context of this article the quantities
  $$m_{ij}=  (-\alpha_i/\sqrt{p},-\alpha_j/\sqrt{p})=(\alpha_i,\alpha_j)/p,\qquad 
  q_{ij}=e^{\frac{2\pi\i}{\ell} (\alpha_i,\alpha_j)}$$ 
  This is precisely the braiding matrix underlying the Nichols algebra $u_q(\g)^+$ with $q$ a primitive $\ell$-th root of unity. Thus Theorem \ref{thm_NicholsAlgebra} shows that $E_i\mapsto \zem_{-\alpha_i/\sqrt{p}}$ give a representation of the Nichols algebra $u_q(\g)^+$ on $\V_{\frac{1}{\sqrt{p}}\Lambda_\g}$ resp. on a subspace of the algebraic closure.
  
   In addition, the relations in Lemma \ref{lm_ScalarProperties} show that $K_i\mapsto e^{\pi\i\;\yer_{\alpha_i/\sqrt{p}}}$ extends this to a representation of the quantum Borel subalgebra ${u_q(\g)^\geq}=u_q(\g)^0u_q(\g)^+$.  
\end{example}

\subsection{The Logarithmic Kazhdan-Lusztig Correspondence}\label{subsec_OutlookKL}

In this article we have constructed a graded infinite-dimensional representation of a diagonal Nichols algebra $\B((q_{ij})_{ij})$, in particular of the Borel part $u_q(\g)^+$ of a quantum group with $q$ a root of unity.
As explained in the introduction, we now  discuss a program that has attracted much attention, and whose goal is to construct vertex algebras with the same representation theory as a quantum group. It was the main motivation for this article. Needless to say, the author has a generalization to diagonal and arbitrary Nichols algebras in mind.

The program is often called \emph{Logarithmic Kazhdan-Lusztig correspondence}, because such an equivalence of finite, non-semisimple modular tensor categories of representations would realize a non-semisimple variant of the Kazhdan-Lusztig correspondence \cite{KL94} between an affine Lie algebra and the respective quantum group. One would even expect to be able to construct from the vertex algebra approach an explicit bimodule, with a vertex algebra action on one side and a quantum group action on the other side. 

\begin{definition}
  Let $\g$ be a finite-dimensional semisimple complex Lie algebra and denote the root-, coroot-, and weight-lattice by $\Lambda_R,\Lambda_{R^\vee}$ and $\Lambda_W$, and a choice of simple roots by $\alpha_1,\ldots,\alpha_r$. For any positive integer $\ell=2p$ divisible by all $(\alpha_i,\alpha_i)$ we consider the so-called \emph{short screenings} $\alpha_i^{short}=-\alpha_i/\sqrt{p}$ spanning a~non-integral lattice $\Lambda=\frac{1}{\sqrt{p}}\Lambda_R$. Then in the context of the present article we have $m_{ij}=(\alpha_i,\alpha_j)/p$ and the braiding matrix $q_{ij}=e^{\pi\i(\alpha_i,\alpha_j)/p}$. Thus Theorem \ref{thm_NicholsAlgebra} proves for these choices that the \emph{short screening operators} $\zem_{-\alpha_i/\sqrt{p}}$ give an action of the small quantum group $u_q(\g)^+$ with $q=e^{\frac{2\pi\i}{2p}}$.
  
  Moreover, we consider an associated set of \emph{long screenings} $\alpha_i^{long}=\alpha_i^\vee\sqrt{p}$, which span the even integral lattice $L=\sqrt{p}\Lambda_{R^\vee}$ with dual lattice $L^*=\frac{1}{\sqrt{p}}\Lambda_W$. 
\end{definition}

The long screenings were chosen such that the following conformal structure exists. For simply-laced Lie algebras this was derived in \cite{FT}, and we computed the non-simply-laced formulae in \cite{FL18}, where the dual rootsystem and the additional divisibility condition appears. For diagonal Nichols algebras the short root lattices are computed in \cite{FL19}, but what the long basis should be is not clear yet. The group $L^*/L$ matches the modular tensor category in our factorizable quasi-quantum group construction \cite{GLO18}.

\begin{lemma}[Virasoro structure]\label{fact_Vir}
 There exists a Virasoro algebra action (or conformal structure) on the lattice vertex algebra $\V_{L}$, such that all $\exp{\phi_{\alpha_i^{long}}}$ and all  $\exp{\phi_{\alpha_i^{short}}}$ in respective modules have $L_0$-eigenvalue (or conformal weight) equal to  $1$. The central charge~is
 $$c= \mathrm{rank}(\mathfrak{g}) -12\left(\vphantom{x^{x^{x^x}}}\vert\rho^\vee\vert^2p-2(\rho, \rho^\vee)+\vert\rho\vert^2/p\right) $$
  with $\rho,\rho^\vee$ the Weyl vectors of the rootsystem and the dual rootsystem.
  
  As primary fields with $L_0$-eigenvalue~$1$, the commutator formula implies that the associated {long screening operators} $\zem_{\alpha_i^{long}}$ are Virasoro algebra homomorphisms. The same conclusion does not hold for short screening operators $\zem_{\alpha_i^{short}}$ acting on an arbitrary module, due to non-local effects (see Lemma \ref{lm_ScreeningProperties}), but we expect suitable powers to be Virasoro homomorphisms, namely the reflection operators in Section \ref{subsec_EqualWeyl}.
\end{lemma}

The long screening operators are local screening operators on $\V_L$, so they generate a Lie algebra by the commutator formula. The following conjecture was proven for simply-laced Lie algebras in \cite{FT} with a trick, but it is probably not difficult to generalize:
\begin{conjecture}
  The long screenings $\zem_{\alpha_i^{long}}$ constitute a representation of the Borel algebra $(\g^\vee)^{\geq0}$ with the dual root system on $\V_L$ and on all its modules.
  
  Note that, as in Lemma \ref{lm_ScreeningProperties} and Section \ref{subsec_TrivialLevel}, this will only hold precisely, with commutators instead of anticommutators, if $\V_L$ is endowed with the usual $2$-cocycle $\epsilon(\alpha,\beta)$. This is the reason behind somewhat obscure sign choices in \cite{FT} formula (2.8) and~(2.9).
\end{conjecture}

Now these operators can be used to define subspaces of the lattice vertex algebra $\V_{L}$. This method of presenting a vertex algebra is called \emph{free-field realization}:

\begin{conjecture}
The kernel of long and short screenings in the lattice vertex algebra~$\V_{L}$
$$\mathrm{W}(\g,\ell):=\bigcap_{i} \ker{\zem_{\alpha_i^{short}}} \;\cap\; \bigcap_{i}\ker{\zem_{\alpha_i^{long}}}\;\cap\; \V_0$$
in $\Lambda$-degree $0$ is isomorphic to the $\mathrm{W}$-algebra (Hamiltonian reduction) of to the respective affine Lie algebra, 
This is related to the free field realization of $\mathrm{W}$-algebras in \cite{FF88}.
\end{conjecture}
\begin{definition}
Consider the larger kernel of only the short screenings:
$$\mathcal{W}(\g,\ell):=\bigcap_{i} \ker{\zem_{\alpha_i^{short}}}\;\cap\;\V_L\hspace{3cm}$$
\end{definition}
\begin{conjecture}
 The kernel of short screenings $\mathcal{W}$ is a vertex subalgebra of $\V_{L}$ with an action of the full Lie algebra $\g^\vee$ including the long screenings. As such, $\W$ is  generated~by a pure exponential $\exp{\phi_{-\omega}}$ for the highest roots $\omega$, spanning the adjoint representation of~$\g^\vee$.
\end{conjecture}
\enlargethispage{1cm}
\begin{conjecture}
The \emph{Logarithmic Kazhdan Lusztig correspondence} asserts that for the category of representations of the vertex operator algebra $\mathcal{W}\subset \V_L$ the following holds:
 \begin{itemize}
  \item It is a non-semisimple modular tensor category. The simple modules are the unique simple submodules $\W_{[\lambda]}$ of the \emph{Verma modules} $\V_{[\lambda]}$ for all $[\lambda]\in L^*/L$, where $\V_{[\lambda]}$ is the restriction of the respective lattice vertex algebra module to $\W$.
  \item It is as an abelian category equivalent to the representation category of a version of the respective quantum group  $u_q(\g)$ with  $q=e^{\frac{2\pi\i}{2p}}$ a $2p$-th root of unity. 
  \item It is as modular tensor category equivalent to the representation category of a factorizable quasi-Hopf algebra variant of the small quantum group $u_q(\g)$ at even roots of unity, which we constructed in \cite{GLO18}.
  \item The action of the mapping class group $\mathrm{SL}_2(\Z)$ on the graded characters \cite{Zhu96} and on the quantum group center coincide, see \cite{FGSTmodular}\cite{RT} in the case $\sl_2$.
 \end{itemize}
Moreover, the representation category of $\mathrm{W}$ (the Hamiltonian reduction) should be equivalent to the representation category of Lusztig's infinite quantum group of divided powers, and  
the representation category of $\W\cap \V_0$ (the singlet algebra) should be equivalent to the representation category of the unrolled quantum group, see e.g. \cite{AM14, CGR20, Len19}.  
\end{conjecture}

\begin{conjecture}
In the authors opinion one should try to approach these conjectures over the following intermediate steps, that also maybe make more clear, why the correspondence would hold:
\begin{itemize}
\item Consider the reflection operators $\zem_{s_{\alpha_i}}$ from Section \ref{subsec_EqualWeyl}. If they do indeed commute with the action of the kernel, as conjectured, then they decompose the lattice vertex algebra modules $\V_{[\lambda]}$ over $\W$. More precisely, one should prove that the kernel of all reflection operators $\W_{[\lambda]}$ for all $[\lambda]\in L^*/L$ gives all simple $\W$-modules. This may require a Zhu algebra argument, quite different from arguments below. 
\item \cite{FT} have proven a remarkable character formula for certain modules of $\W$ constructed by cohomologies of vertex algebra bundles over the flag variety. These should be again the  $\W_{[\lambda]}$, but this is not proven. The inclusion $\W_{[\lambda]}^{Len}\subset \W_{[\lambda]}^{FT}$ would follow from extending the action of the long screening to a full $\g^\vee$ action, probably via divided powers of the short screenings. Vice-versa one could argue with the known characters that equality holds. One result of this proof would be that the composition factors of the Verma modules $\V_{[\lambda]}$ are as expected. 
\item One should prove the linkage principle, namely that $\W_{[\lambda]},\W_{[\mu]}$ are in the same block only if they are in an orbit of the Weyl group. Probably this could be deduced from the respective linkage principle for $\mathrm{W}$-algebras in \cite{Arak07}. 
\item \cite{AJS94} have proven the conjecture by Lusztig relating the principal block of the representations of a small quantum group $u_q(\g)$ to the representations of $\g$ in characteristic $p$. Their remarkable proof consists, very roughly, of setting up a machinery of translation- and wall-crossing-functors in both categories and then showing that they have the same combinatorics. After establishing the facts above, it would be intriguing to try to set up the same machinery a third time, now in the category of $\W$-representations. Note that this makes heavy use of~the tensor category structure, which is the most difficult piece of vertex algebra representation theory, but whose existence is assured by \cite{HLZ-VIII} and which has the ability to produce proven projective modules: By  tensoring the Steinberg module.  
\item The equivalence should be made explicit by constructing a projective module over~$\W$ with the action of a quantum group. This was done successfully for $\g=\sl_2$ at $p=2$ in \cite{FGSTsl2} -- however $p=2$ is the only case that  does not involve truely non-local operators and in the general case the quantum group action will not preserve $\V_{\Lambda}$. Instead the screening operators in the present article will have to be taken literally, giving infinite sums. Then, the idea is again to use reflection operators to glue together a bimodule of the form $\bigoplus_{[\lambda]\in L^*/L}\W_{[\lambda]}\otimes P_{[\lambda]}\cong \bigoplus_{[\lambda]\in L^*/L}\mathcal{P}_{[\lambda]}\otimes L_{[\lambda]}$, where $\mathcal{P}_{[\lambda]}$ is the projective cover of $\W_{[\lambda]}$, and where $P_{[\lambda]}$ is the projective cover of the irreducible representation $L_{[\lambda]}$ of $u_q(\g)$. In view of the course of proof in \cite{AJS94}, the case $\sl_2$ seems to be an essential step. 
\end{itemize}
 \end{conjecture}

On the other hand the author would like to extend the original question to pointed Hopf algebras and beyond. The following ultimate goal is probably very hard:
\begin{problem}\label{problem_main}
	Construct for every modular tensor category $\mathcal{C}$ at least one vertex algebra $\W$, such that its representation category is the given modular tensor category $\mathcal{C}$.  	
\end{problem}
We have still a poor understanding of the classification of  semisimple tensor categories. On the other hand, we have the successful Andruskiewitsch-Schneider program \cite{AS}, which classifies Hopf algebras $H$ under the assumption that the coradical (maximal cosemisimple subalgebra) $H_0$ is of a given form, and via the classification of all Nichols algebras over $H_0$. 
For example, if $H_0$ is a abelian group ring, then $H$ is a generalization of a quantum group, involving any of the diagonal Nichols algebras classified in \cite{Heck09}. These correspond to quantum groups, quantum supergroups and some exotic examples.

\begin{problem}
	In this spirit, the following problem might be approachable: 
	For a given class of semisimple modular tensor categories $\mathcal{S}$, and given vertex algebras $\V$ with this representation theory, construct for every modular tensor category $\mathcal{C}$ with semisimple part  $\mathcal{S}$ (in some sense) at least one vertex algebra $\W$ with representation category $\mathcal{C}$.
	
	One could tackle this problem by trying to write $\mathcal{C}$ as the category of Yetter-Drinfeld modules over a Nichols algebra $\B$ inside $\mathcal{S}$. Then one could use the realization of $\mathcal{S}$ via $\V$ to define an action of $\B$ on $\V$ via non-local screening operators, using generalizations of the technology developed in this article. Then one could hope that the kernel of screenings $\W=\V^\B$ has always a representation theory equivalent to $\mathcal{S}$. 
\end{problem}

Compare this to group orbifolds, where the representation category of $\V^G$ is some extension of the representation category of $\V$ by the Drinfeld double of $G$. In this analogy, ``orbi\-folding $\V_L$ by a Nichols algebra $\B$'' gives a quantum group with coradical~$L^*/L$.  

\begin{example}[Pointed Hopf algebras with abelian coradical]
	The case $\mathcal{S}=\mathrm{Vect}_A$ is the setting of \cite{AS}, so we know all respective Hopf algebras in terms of $H_0=\C[A]$ and a known list of diagonal Nichols algebras. One might expect, from the explicit formulae, that the only modular tensor categories do arise for Drinfeld doubles of one Nichols algebra.
	In this case we have a realization of $\mathcal{S}$ via lattice vertex algebras, and the result of this paper realizes any diagonal Nichols algebras as algebras of screenings over these vertex algebras. We calculated a list of corresponding lattices in \cite{FL19}, following \cite{ST}. 
\end{example}
	
 One step further, the author has in \cite{Len12}\cite{Len14} constructed non-diagonal Nichols algebras by  \emph{folding} and \cite{HV} show they are essentially all such Nichols algebras (and probably the only ones leading to modular tensor categories besides Drinfeld doubles).

\begin{conjecture}
  Let $\g$ be a Lie algebra, $\ell=4$ and $\theta$ an automorphism of the Dynkin diagram. Associated to this data we can consider the following three objects
  \begin{itemize}
   \item An orbifold model of the vertex operator algebra $\W(\g,\ell)$ by the automorphism $\theta$.
   \item A kernel of screening operators in the orbifold model of the vertex algebra of the rescaled root lattice of $\g$ by $\theta$. Here, the Nichols algebra braiding is non-diagonal.
   \item A factorizable Hopf algebra $\tilde{u}_q({^\theta}\g)$ associated to the author's non-diagonal Nichols algebra for the folding data $\g,\theta$ and $q$ a primitive $\ell$-th root of unity. 
  \end{itemize}
  Are their representation categories equivalent? (the first two seem to be clearly equal)
\end{conjecture}

Even one step further, one should consider a case beyond quasi Hopf algebras. Most suitable seems the case where $\V$ is an affine Lie algebra $\hat{\g}$ at integer positive level, and the associated semisimple modular tensor categories $\mathcal{S}$. It has additional advantages: On one side the groundstates and possible screenings have clear Lie-theoretic interpretations, and on the other side $\mathcal{S}$  can be studied using semisimplification of quantum groups.

\subsection{Example: \texorpdfstring{$\mathfrak{sl}_2$ at $\ell=4$}{sl2 at ell=4}}\label{subsec_ExampleKL}

We will make parts of the vague and conjectural outlook in the previous section precise on the example $\g=\mathfrak{sl}_2$ where the main conjectures have been proven in \cite{FGSTsl2} for $\ell=4$ and in \cite{TN} for arbitrary $\ell$. The second source, despite being a remarkable paper, might be considered sketchy in some details.\\

In the case $\sl_2$ we have $\Lambda_R=\alpha\Z$ with $(\alpha,\alpha)=2$. For $\ell=2p\in\N$ we consider the lattice vertex algebra $\V_{L}$ with $L=\alpha\sqrt{p}\;\Z$ and $L^*=\alpha/2\sqrt{p}\;\Z$, which has $2p$ modules~$\V_{[k]}$ that are enumerated by cosets $[k]\in \Z_{2p}$ corresponding to lattice cosets $[k\alpha/2\sqrt{p}]\in L^*/L$. \\

One easily finds the conformal element $T=\frac{1}{4} \partial\phi_\alpha\partial\phi_\alpha
-(\sqrt{p}-1/\sqrt{p})\partial^2\phi_{\alpha/2}$, giving a Virasoro algebra action with central charge $13-6(p+\frac{1}{p})$, such that  the two pure exponentials $\exp{-\alpha/\sqrt{p}}$ and $\exp{+\alpha\sqrt{p}}$ have $L_0$-eigenvalue $1$. In general, the $L_0$-eigenvalues of the pure exponentials $\exp{\phi_{k\alpha/2\sqrt{p}}}$ is as follows 
$$h_k:=h_{k+1,1}=\frac{(k-(p-1))^2}{4p}+\frac{c-1}{24}=\frac{k(k-2p+2)}{4p}$$
This is the characteristic parabola centered around the \emph{Steinberg point} ${Q=(p-1)\alpha/2\sqrt{p}}$ for $k=p-1$, visible in the picture below.
This Virasoro action is chosen with respect to  $\exp{-\alpha/\sqrt{p}},\;\exp{+\alpha\sqrt{p}}$, so that it is compatible with the \emph{short and long screening operators}
$$\zem_{-\alpha/\sqrt{p}},\qquad \zem_{+\alpha\sqrt{p}}$$

The groundstates (in the sense of Zhu's algebra) of the module $\V_{[k]}$ are $\exp{\phi_{k\alpha/2\sqrt{p}}}$ for those coset representatives $k\alpha/2\sqrt{p}\in L^*$ with minimal 
distance to the Steinberg point. So the groundstate is $2$-dimensional with $k=-1,2p-1$ for the 
class $[-1]\in \Z_{2p}$ and it is $1$-dimensional with $k\in\{0,\ldots 2p-2\}$ for all other classes $[k]\in \Z_{2p}$.\\

Our results in this article, applied to $m_{11}=(-\alpha/\sqrt{p},-\alpha/\sqrt{p})=2/p$ and the braiding  $q_{11}=e^{\frac{2\pi\i}{p}}$, proves that the Nichols algebra relation $(\zem_{-\alpha/\sqrt{p}})^p=0$ holds and thus the short screening generates $u_q(\mathfrak{sl}_2)^-$. On the other hand, $e^{\pi\i(+\alpha\sqrt{p},+\alpha\sqrt{p})}=e^{2\pi\i p}=1$ and the long screening 
$\zem_{+\alpha\sqrt{p}}$ is local and generates the Borel part $U(\mathfrak{sl}_2)^+$.\\


The kernel of the short screening $\W(\sl_2,\ell)$ is the \emph{triplet algebra} $\mathcal{W}(p)$, previously discovered by \cite{Kausch91}, for which the present construction gives a free-field realization:
$$\mathcal{W}(p)\cong \V_{\sqrt{p}\alpha\Z}\;\cap\; \ker{\zem_{-\alpha/\sqrt{p}}}$$

This vertex operator algebra is known to be generated by a pure exponential together with the $\mathfrak{sl}_2$-action of the long screenings, generating an adjoint representation of $\mathfrak{sl}_2$:
\begin{align*}
 \mathcal{W}^- &:=e^{-\phi_{\sqrt{p}\alpha}}\\
 \mathcal{W}^0\;&:=\zem_{+\alpha\sqrt{p}} e^{-\phi_{\sqrt{p}\alpha}}\\
 \mathcal{W}^+ &:=\zem_{+\alpha\sqrt{p}}\zem_{+\alpha\sqrt{p}} e^{-\phi_{\sqrt{p}\alpha}}\hspace{4cm}
\end{align*}
\begin{example}
For  $p=1$  we get the affine Lie algebra at level one, see Section \ref{subsec_TrivialLevel}:
\begin{align*}
  \mathcal{W}^-&:=\exp{\phi_{-\alpha}}\\
  \mathcal{W}^0\;&:=\partial\phi_\alpha\\
  \mathcal{W}^+&:= -2\exp{+\phi_\alpha}\hspace{6cm}
\end{align*}
The first nontrivial case $p=2$ was studied in \cite{FGSTsl2} in terms of symplectic fermions
\begin{align*}
 \mathcal{W}^-&:=\exp{\phi_{-\alpha\sqrt{2}}}\\
 \mathcal{W}^0\;&:=\frac{1}{6}\left(
\partial\phi_{\alpha\sqrt{2}}\partial\phi_{\alpha\sqrt{2}}
\partial\phi_{\alpha\sqrt{2}}
+3\partial\phi_{\alpha\sqrt{2}}\partial^2\phi_{\alpha\sqrt{2}}
+\partial^3\phi_{\alpha\sqrt{2}}\right) \hspace{-1.5cm}\\
 \mathcal{W}^+&:=8\left(\partial\phi_{\alpha\sqrt{2}}\partial\phi_{\alpha\sqrt{2}}+\partial^2\phi_{\alpha\sqrt{2}}\right)
 \exp{\phi_{+\alpha\sqrt{2}}}
\end{align*}
\end{example}

\noindent
\begin{theorem}[\cite{AM08}] 
The representations category of $\mathcal{W}(p)$ for $\;p\geq 2$ is as follows:
\begin{enumerate}
 \item It has $2p$ irreducible representations, denoted $\Lambda(k+1),\Pi(k+1)$ for ${k=0,\ldots, p-1}$, with
 groundstates\footnote{Semikhatov used the Cyrillic letters  $\Lambda,\Pi$ in order to visualize the $1$- resp. $2$-dim groundstates.} of dimension $1$ resp. $2$ and $L_0$-eigenvalues $h_{k}$ resp. 
 $h_{(p-k-2)+2p}$.
 \item It is non-semisimple. There are intertwining operators, which have logarithmic 
 singularities. There exist modules, where $L_0$ acts 
  non-diagonalizable.
 \item It is $C_2$-cofinite \cite{CF}\cite{AM08}. 
 \end{enumerate}
\end{theorem}

\noindent
We now discuss \emph{Verma modules}, the restriction of the simple $\V_{L}$-modules $\V_{[k]}$ to $\mathcal{W}(p)$:
For the short screenings $\zem_{-\alpha/\sqrt{p}}$ we consider the reflection operators in Subsection~\ref{subsec_EqualWeyl}
$$\left(\zem_{-\alpha/\sqrt{p}}\right)^{n}:\;
\V_{[{(p-1)+n}]}\longrightarrow \V_{[{(p-1)-n}]}$$
For the {Steinberg points} $k=p-1$ and $k=2p-1\equiv -1$ (which is a non-regular weight) the reflection operator maps the module $\V_{[k]}$ to itself, the reflection operator is trivial, and the module stays simple over $\mathcal{W}(p)$. We call these $\mathcal{W}(p)$-modules 
$$\Lambda(p):=\V_{[{p-1}]}
\qquad \Pi(p):=\V_{[{-1}]}$$
All other simple $\V_L$-modules arise in pairs $k=(p-1)\pm n$ around the Steinberg point and each restricts over $\W(p)$ to an indecomposable extensions of two simple $\mathcal{W}(p)$-modules, namely  image-by-kernel of the respective reflection operators~$(\zem_{-\alpha/\sqrt{p}})^{n},(\zem_{-\alpha/\sqrt{p}})^{p-n}$
$$0\to \Lambda(p-n)\to\V_{[(p-1)-n]}\to\Pi(n)\to 0\qquad\hspace{1cm} n=1,\ldots, p-1$$
$$0\to \Pi(n)\to\V_{[(p-1)+n]}\to\Lambda(p-n)\to 0\hspace{4.4cm}$$

One can show this abelian category is equivalent to $u_q(\mathfrak{sl}_2)$ at $\ell=2p$. The $2$ semisimple blocks and $p-1$ non-semisimple blocks are according to the Weyl group dot action of~$\sl_2$.

\begin{remark}[]
 The lattice vertex algebra $\V_{\sqrt{p}\alpha\Z}$ has the graded dimensions
  \begin{align*}
 \ch_{\V_{[{(p-1)\pm n}]}}(q)
 &=\frac{\Theta_{p,n}(q)}{\eta(q)}\\
\intertext{which decompose over $\W(p)$ into the following graded dimensions, see \cite{AM08}}
 \ch_{\Lambda(n)}(q)
     &=\frac{n\;\Theta_{p,p-n}(q)+2p\;\partial\Theta_{p,p-n}}{p\;\eta(q)}\\
 \ch_{\Pi(n)}(q)
     &=\frac{n\;\Theta_{p,n}(q)-2p\;\partial\Theta_{p,n}}{p\;\eta(q)}
 \end{align*}
 where one defines 
 $$\Theta_{p,n}(q):=
 \sum_{j\in \Z}q^{\frac{1}{4p}(n+2pj)^2} 
 \qquad
 \partial\Theta_{p,n}(q):=
 \sum_{j\in \Z}\left(j+\frac{n}{2p}\right)q^{\frac{1}{4p}(n+2pj)^2}$$
 These $2p$ characters together with $p-1$ pseudocharacters, see \cite{CG17}, give a vector-valued modular form with $3p-1$ components. The $\mathrm{SL}_2(\Z)$ action matches the $\mathrm{SL}_2(\Z)$-action on the center of $u_q(\sl_2)$ \cite{FGSTmodular}.  
\end{remark}

To conclude, we summarize the information about $\W(p)$ for $p=2$ in the following picture: 
 Recall that $\V_L$ has $2p=4$ representations $\V_{[j]}=\bigoplus_{k\in[j]}\V_{k\alpha/2\sqrt{p}}$ with $[j]\in\Z_{2p}$. Dots denote basis vectors of the representations of $\V_L$ and fat dots represent the pure exponentials $\exp{\phi_{k\alpha/2\sqrt{p}}}$ on top of each cone-shaped $\H^r$-submodule $\V_{k\alpha/2\sqrt{p}}$.  The $\Y$-axis denotes $L_0$-eigenvalues and the $X$-axis denotes $k$. We draw the vacuum module $\V_{[0]}$ as a union of cones with boundary lines, the module $\V_{[2]}$ as a union of cones without boundary lines, and we only indicate with gray fat dots the Steinberg modules $\V_{[-1]}$ and $\V_{[1]}$ that will remain irreducible over $\W(p)$.
The straight short arrows denote short screenings with $(\zem_{-\alpha/\sqrt{p}})^{2}=0$  
and the bent long arrows denote long screenings giving the $\mathfrak{sl}_2$ action. The union of the gray areas is $\W=\Lambda(1)$, the kernel of the short screening in the vacuum representation $\V_L=\V_{[0]}$, which is the new logarithmic vertex algebra. We see the decompositions of $\V_{[0]}$ into white and grey areas according to 
$\Lambda(1)\to \V_{[0]}\to \Pi(1)$.
\begin{center}
 \includegraphics[scale=.35]{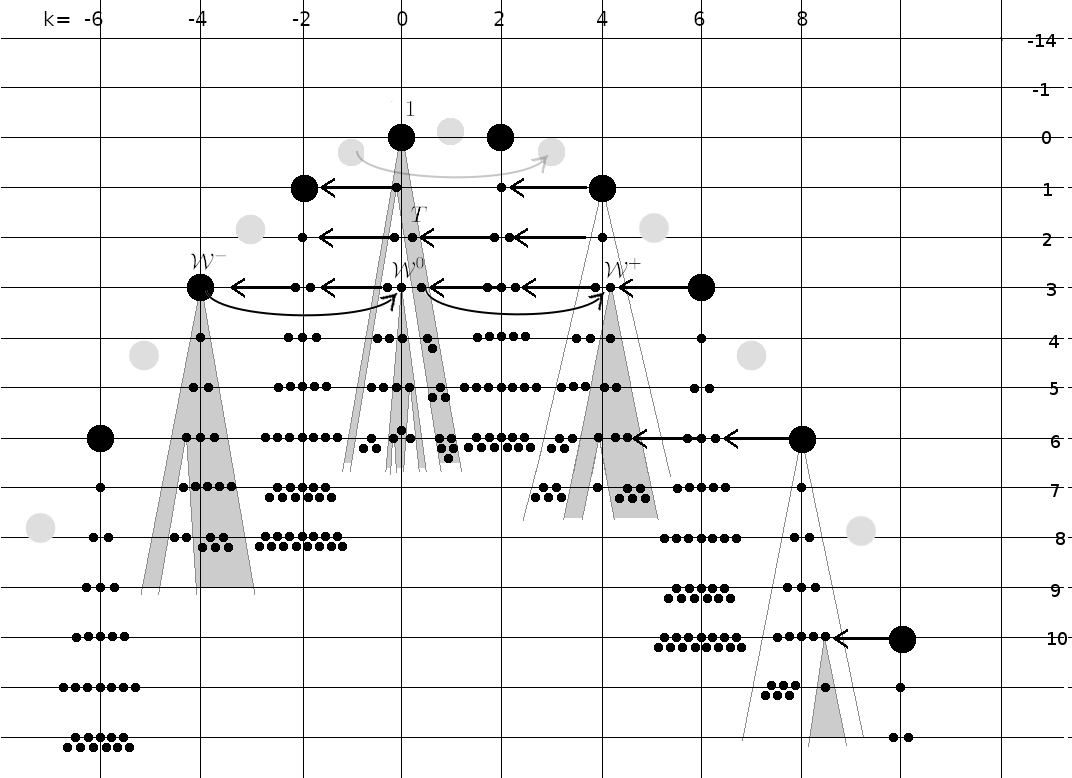}
\end{center}

\enlargethispage{2.5cm}
\begin{acknowledgementX}
I am indebted and grateful to A. Semikhatov, for hospitality in Moscow and 
for countless hours of exposing me to the topic and his ideas, to C. Schweigert, 
whose continuing support and inspiration made all of this possible, and to  B. Feigin, who  
created such beauty, continuously answered my questions and explained the big picture. \\

Also many thanks to Christian Reiher for enduring help on the analysis side and to Yorck Sommerh\"auser and the anonymous referee for carefully reviewing and improving the manuscript.
My visit in Moscow with the DAAD PRIME program was funded by the BMBF Germany and the Marie Curie Action of the EU. I am also grateful for partial support from the DFG Graduiertenkolleg 1670 at the University Hamburg.\\

\noindent 
To my son Jonathan.

\end{acknowledgementX}

\newpage
\renewcommand{\stretch}{\enlargethispage{0cm}}
\providecommand{\href}[2]{#2}

\end{document}